\documentclass[a4paper,11pt,oneside,leqno]{article}

\usepackage{amsmath,amsthm,amsfonts,amssymb,amscd}
\usepackage{mathtools}
\usepackage{pgf,tikz}
\usetikzlibrary{calc}
\usepackage{mathrsfs}
\usepackage{float}
\usepackage{csvsimple}
\usepackage{enumerate}
\usepackage{booktabs}
\usepackage{pgfplots}
\usepackage{pgfplotstable}
\usepackage[noend]{algpseudocode}
\usepackage{graphicx,bm,xcolor}
\usepackage{algorithm}
\usepackage{pdfpages}
\usepackage{algpseudocode}
\usepackage{stmaryrd}
\usetikzlibrary{arrows}
\usepackage{caption}

\usepackage{t1enc}%%%%%%%%%%%%%%%%%%%%%%%%%%%%%%%%%%%%%%%%%%%%%%%%%%%%%%%%%%%%%%%%%%%%%%%%%%%%%%%%
\usepackage{amsmath,amsthm,amsfonts,amssymb,amscd}
\usepackage{graphicx}
\usepackage{mathtools}
\usepackage{pgf,tikz}
\usepackage{mathrsfs}
\usepackage{float}
\usepackage{enumerate}
\usepackage{booktabs}
\usepackage{pgfplots}
\usepackage[noend]{algpseudocode}
\usepackage{graphicx,bm,xcolor}
\usepackage{algorithm}
\usepackage{algpseudocode}
\usepackage{stmaryrd}
\usetikzlibrary{arrows}
\usepackage{caption}
\usepackage[T1]{fontenc}
\usepackage{t1enc}
\usepackage[polish,german,english]{babel}
\usepackage{verbatim} 

\usepackage{url}

\newif\ifMAKEPICS
\MAKEPICSfalse

\ifMAKEPICS
\usepackage[cleanup,subfolder]{gnuplottex}
\usepackage{xparse}

\ExplSyntaxOn
\DeclareExpandableDocumentCommand{\convertlen}{ O{cm} m }
{
	\dim_to_decimal_in_unit:nn { #2 } { 1 #1 } cm
}
\ExplSyntaxOff
\fi

\usepackage[german,english]{babel}
\usepackage{verbatim} 

\usepackage{url}

\newif\ifMAKEPICS
\MAKEPICSfalse

\ifMAKEPICS
\usepackage[cleanup,subfolder]{gnuplottex}
\usepackage{xparse}

\ExplSyntaxOn
\DeclareExpandableDocumentCommand{\convertlen}{ O{cm} m }
{
	\dim_to_decimal_in_unit:nn { #2 } { 1 #1 } cm
}
\ExplSyntaxOff
\fi

\usepackage{authblk}
\usepackage{todonotes}
\usepackage[top=2cm, bottom=2cm, left=2cm, right=2cm]{geometry}
\theoremstyle{plain} % default
\newtheorem{theorem}{Theorem}[section]
\newtheorem{propositon}{Proposition}[section]
\newtheorem{lemma}[theorem]{Lemma}

\newtheorem{remark}[theorem]{Remark}
\newtheorem{definition}[theorem]{Definition}

\newtheorem{assumption}{Assumption}

\theoremstyle{definition} %

\theoremstyle{remark} %

\newcommand{\ignore}[1]{}

\begin{document}

\title{
  Multigoal-oriented optimal control problems \\
  with nonlinear PDE constraints
} 
\author[1,2]{B. Endtmayer}
\author[2]{U. Langer}
\author[3]{I. Neitzel}
\author[4,5]{T. Wick}
\author[6]{W. Wollner}

\affil[1]{Doctoral Program on Computational Mathematics, Johannes Kepler University, Altenbergerstr. 69, A-4040 Linz, Austria}
\affil[2]{Johann Radon Institute for Computational and Applied Mathematics, Austrian Academy of Sciences, Altenbergerstr. 69, A-4040 Linz, Austria}
\affil[3]{Institut f\"ur Numerische Simulation,
  Endenicher Allee 19b,
  53115 Bonn,
  Germany}
\affil[4]{Institut f\"ur Angewandte Mathematik, Leibniz Universit\"at Hannover, Welfengarten 1, 30167 Hannover, Germany}
\affil[5]{Cluster of Excellence PhoenixD (Photonics, Optics, and
  Engineering - Innovation Across Disciplines), Leibniz Universit\"at Hannover, Germany}
\affil[6]{Technische Universit\"at Darmstadt,
  Fachbereich Mathematik,
  Dolivostr. 15,
  64293 Darmstadt, Germany }

\date{}

\maketitle

\begin{abstract}

  In this work, we consider an optimal control problem subject to a nonlinear PDE constraint
 and apply it to the regularized $p$-Laplace equation. To this end, a reduced
  unconstrained optimization problem in terms of the control variable is formulated.
  Based on the reduced approach, we then derive an a posteriori error representation and mesh adaptivity for
  multiple quantities of interest. 
  All quantities are combined to one, and then the dual-weighted residual (DWR) method
  is applied to this combined functional. Furthermore, 
  the estimator allows for balancing the
  discretization error and the nonlinear iteration error. 
  These developments allow us to formulate an adaptive solution strategy, which 
  is finally substantiated 
  via
  several numerical examples.
\end{abstract}

% ##############################################################################%
\section{Introduction}
\label{Introduction}
Optimal control problems with nonlinear PDE constraints have been studied 
for a long time in many works. In particular, employing 
the (regularized) $p$-Laplacian (see e.g.,~\cite{Glowinski1975,DiRu07,Hi2013,ToWi2017}) as a nonlinear constraint
of an optimal control problem was considered 
for instance in~\cite{CaKoLeu16}.

In many applications, however, not the entire solution is of interest, but 
only parts or certain quantities of interest, so-called goal functionals. 
In the past, often a single goal functional was analyzed. However,
it may be of interest to control multiple goal functionals simultaneously 
\cite{HaHou03, Ha08, BruZhuZwie16,EnWi17,KerPruChaLaf2017, PARDO20101953}.
In this paper, these three topics are combined: optimal control, the
regularized $p$-Laplacian as a numerical example of a quasi-linear PDE constraint, 
and multiple goal-oriented a posteriori error estimation.

In the following, we briefly refer to studies that treat parts of the three topics.
Optimal control problems (specifically, a priori estimates and optimality conditions)
with quasi-linear (as the $p$-Laplacian can be classified) 
elliptic 
PDE constraints
were considered in \cite{CaFe89,CaTroe09,CaDah11}.
More recently,
the extension to optimal control with parabolic PDEs
was discussed in~\cite{BoNei18} and~\cite{CasChry18}.

Optimal control problems with (single) goal functionals 
were investigated in~\cite{BeKaRa00,MeidnerVexler:2007,BeBaMeRa07,VeWo08,Wollner:2008,RannacherVexlerWollner:2011}. 
The $p$-Laplacian and a posteriori
error estimates were considered in~\cite{Liu2001,CaKlo03,Creuse2007,CaWeNi06}, and, more specifically,
for goal functional evaluations, we refer to~\cite{Hi2013,RanVi2013,EnLaWi18}.
To estimate goal functionals, we adopt the dual-weighted residual (DWR) method~\cite{BeckerRannacher:1996,BeRa01}
in which an adjoint problem is solved to obtain (local) sensitivity measures
that are used for mesh refinement. 
As is well-known, using a gradient-based
approach for the numerical solution of optimal control problems, 
the same adjoint problem as for the DWR error estimator can be employed. For this reason,
it is natural to combine gradient-based optimization with adjoint-based error
estimation.

We are specifically interested in an extended DWR version
in which the discretization and (linear/nonlinear) iteration error are balanced
\cite{MeiRaVih109,RanVi2013,MalGourVohYou2018}. 
As localization technique we employ integration by parts as done in \cite{BeRa01} 
or, for residual based error estimates, in \cite{Verfuerth:1996a}. 
The extension of \cite{RanVi2013} to multiple goal functionals was recently undertaken 
in \cite{EnLaWi18}.

Three major aims constitute the main contents of this paper: first, the design of a framework for 
goal-oriented error estimation
for optimal control subject to a nonlinear PDE
and balancing 
the discretization and nonlinear iteration error (Section~\ref{sec_DWR}).
From the optimization point of view, we carefully revisit the important 
elements for the DWR estimator for optimization problems. The main result in
this respect 
is the a posteriori error representation for the reduced optimal control system for an abstract problem formulation.
The second aim is the extension to the simultaneous control of multiple goal
functionals (Section~\ref{sec_multigoal}). 
As a third goal, based on our theoretical developments, we carefully
design an adaptive solution algorithm (Section~\ref{sec_algo}). The
performance of our algorithms are investigated in terms 
of the usual quality measures of convergence behavior and effectivity indices in Section~\ref{sec_num_tests}. 
The latter one measures the quality of our proposed error estimator in
comparison to (known) true errors, which are computed on sufficiently refined
meshes.

We summarize the outline of this work as follows: In Section~\ref{sec_problem}, the 
problem setting is introduced. 
Next, in Section~\ref{sec_DWR}, the dual-weighted residual method 
for the reduced optimization problem is formulated. The multi-goal approach 
is then introduced in Section~\ref{sec_multigoal}. Our algorithmic
developments to solve the multiple goal-functional optimal control problem 
are derived in Section~\ref{sec_algo}. In Section~\ref{sec_num_tests}, we 
present several numerical examples that demonstrate the performance of our 
approach. Therein, we study different Tikhonov regularization parameters, we 
perform mesh refinement studies, and consider different goal functionals.
In Section~\ref{sec_conclusions}, we summarize the key outcomes of this work.

\newpage
% ##############################################################################%
\section{The Optimal Control Problem}
\label{sec_problem}
In this section, we define an abstract problem formulation and collect some
properties that we will rely on when deriving the a posteriori error
estimates. 
%It will depend upon the specific PDE to prove these properties, which is, however, not within the scope of this paper.
\subsection{The Abstract Problem Formulation}
Let $U$ and $Q$ be Banach spaces. We would like to find  a control
$\overline{q} \in Q$ and an associated state $\overline{u} \in U$ such that
the pair $(\overline{u}, \overline{q})$ is a local minimizer of some given
cost functional $J(u,q)\colon U\times Q\to\mathbb{R}$,
where $u$ and $q$ have to fulfill the so called state equation $A(u,q)=0$ with nonlinear differential operator $A$ acting between Sobolev spaces. More precisely, the arising PDE-constrained optimization problem reads as follows:

\begin{align}  \label{Problem: Abstract Optimization Problem}
	\begin{split}
		\min_{(u,q)} J(u,q) & \qquad  u \in U, q \in Q,\\
		\text{s.t.  }A(u,q)=0& \qquad \text{ in } V^*,
	\end{split}
\end{align}
for some operator $A\colon U \times Q \mapsto V^*$, where $V^*$ denotes the dual space of some Banach space $V$. For the theoretical findings in this paper, we assume that, for each $q\in Q$, the PDE is uniquely solvable. 
More precisely, we assume the following:
\begin{assumption}\label{Assumption: Existence of the solution operator}
  Let there exist a unique mapping $ S\colon Q \mapsto U$ which is implicitly defined by 
  \begin{equation}\label{Equation: Solution Operator}
    A(S(q),q)=0, \qquad \forall q \in Q.
  \end{equation}
  Moreover, we assume that $S$ is twice continuously Fr{\'e}chet differentiable.
\end{assumption}
Without further mention, we also assume the existence of a at least one global minimizer for Problem \eqref{Problem: Abstract Optimization Problem}. 
For instance, we refer
to \cite{Troe09} for general theorems on existence of solutions for problems with linear and semilinear state equations.
Moreover, let $A$ and $J$ be smooth enough for all operations occurring in the next Section.

With the help of the so called control-to-state mapping $S$, we reformulate~\eqref{Problem: Abstract Optimization Problem} as an unconstrained optimization problem
\begin{equation}\label{eq:reducedproblem}
  \min_q j(q),  \qquad q \in Q,\nonumber
\end{equation}
where $j(q) := J(S(q),q)$. Here, we will also assume sufficient smoothness in order to derive all further estimates.

\subsection{First Order Necessary Optimality Conditions}\label{Section: First Order Necessary Optimality Conditions}
It is clear, that under our implicit smoothness assumptions, the first order necessary optimality conditions for a locally optimal control $\bar q\in Q$ for Problem~\eqref{eq:reducedproblem} are given by 
\begin{equation} \label{Problem: Reduced Primal Problem}
  j'(\overline{q})(\delta q)=0 \qquad\forall \delta q \in Q.
\end{equation}
For completeness and further use, we rewrite these conditions for the
non-reduced formulation with the help of the well-known Lagrange approach. We define the Lagrangian $\mathcal{L}\colon U \times Q \times V \mapsto \mathbb{R}$ for this problem  as follows
\begin{align} \label{Lagrangian}
  \mathcal{L}(u,q,z) := J(u,q)-A(u,q)(z), \qquad \forall u \in U, q \in Q, z \in V.
\end{align}
To shorten notation, we consider the abbreviation  $B'_{\zeta}:=\frac{\partial}{\partial \zeta}B$ for the partial derivatives of some operator $B$.
The first order necessary optimality conditions for~\eqref{Problem: Abstract Optimization Problem} are then given by
\begin{equation} \label{Equation: First Order Opt u}
  \begin{aligned}
    J'_{u}(\bar u,\bar q)(\delta u)-A'_{ u}(\bar u,\bar q)(z)(\delta u)=\mathcal{L}'_{u}(\bar u,\bar q,\bar z) (\delta u)&=0 \qquad \forall \delta u \in U, \\
    J'_{q}(\bar u,\bar q)(\delta q)-A'_{q}(\bar u,\bar q)(\bar z)(\delta q)=\mathcal{L}'_{q}(\bar u,\bar q,\bar z) (\delta q)&=0 \qquad \forall \delta q \in Q, \\
    -A(\bar u,\bar q)(\delta z)=\mathcal{L}'_{z}(\bar u,\bar q,\bar z)(\delta z) &=0  \qquad \forall \delta z \in V.
  \end{aligned}
\end{equation}
Moreover, $\bar u =S\bar q$ denotes the optimal state associated with $\bar q$, and $\bar z=(S'(\bar q))^*J'_u(\bar u,\bar q)$ the associated adjoint state.
In order for the Newton algorithm to work, and for the error estimator we need the following assumption.
\begin{assumption}
  We assume that $A'_{ u}= \mathcal{L}''_{uz}$ is invertible.
\end{assumption}

\subsection{An Example: the Regularized $p$-Laplacian and Tracking-type Cost Functional} \label{Example: p-Laplace}
Let us finish this section by defining $A$ for a concrete example (i.e., a PDE) that
motivates our numerical studies. To this end, a (regularized) $p$-Laplace equation
for $p\neq 2$ is considered, even though, it does not necessarily fit into the theory setting. 
For details, we refer to~\cite{DiRu07,Hi2013,ToWi2017} and the references therein
regarding the (regularization of) the $p$-Laplace equation. 
We consider the following setting: Let $\Omega \subset \mathbb{R}^d$ be open
and bounded with $C^1$ boundary, and let $p \in (\frac{2d}{2+d}, \infty)$. Then
we define 
\[\mathcal{A}_p\colon  W^{1,p}_0(\Omega) \times (W^{1,p}_0(\Omega))^*\mapsto
  (W^{1,p}_0(\Omega))^*,\]  
by the identity
\begin{align}
  \mathcal{A}_p(u,q)(v) :=& \langle{(\varepsilon^2+ |\nabla u|^2)^{\frac{p-2}{2}}\nabla u,\nabla v\rangle}_{(L^p(\Omega))^*\times L^p(\Omega)}-\langle{f+q,v\rangle}_{ (W^{1,p}_0(\Omega))^*\times W^{1,p}_0(\Omega)},\nonumber
\end{align}
for $u,v \in U:= W^{1,p}_0(\Omega),$ $f \in  V^*$, where $\langle\cdot,\cdot\rangle$ is the usual notation for duality pairings. Note that in this example, we have $U = V$. 
\label{sec_2}
Let $u\in U$ be the state, and $q\in Q$, e.g., $Q=L^2(\Omega)$, be the control variable.
Then  our optimal control problem is given by
\begin{align}
  \min_{(u,q)}J(q,u) \qquad u \in U, q \in Q \\\text{s.t. } \quad \mathcal{A}_p(u,q)(v) = 0, \nonumber
\end{align}
with the tracking-type cost functional
\begin{align*}
  J(q,u) = \frac{1}{2}\|u - \bar{u}^d\|_{L^2(\Omega)}^2 + \frac{\alpha}{2}\|q - \bar{q}^d\|_{L^2(\Omega)}^2,
\end{align*}
with $\alpha >0$ and given $f \in U^*$, $\bar{u}^d \in L^2(\Omega)$ and $\bar{q}^d\in L^2(\Omega)$. 

% ##############################################################################%
\section{The Dual Weighted Residual Method for the Reduced System}\label{Section: DWR for Reduced System}
\label{sec_DWR}
We now formulate the DWR method for the reduced optimal control
system and develop a posteriori error estimators. The presentation is kept as general as possible so that the 
extension to multiple goal functionals outlined in Section~\ref{sec_multigoal}
can easily be  incorporated.
Firstly, we briefly outline the important elements of the discretization.

\subsection{Discretization}\label{sec_disc}
The method of choice, which will be used in the numerical examples, is the
finite element method~\cite{Ciarlet:2002:FEM:581834,
  Braess,GrRoSty07}. However, the algorithms presented in this work can also be
adapted to other discretization techniques where adaptivity can be
accomplished, like isogeometric analysis, the virtual element
method, or finite cell methods.
For the spaces $U_h=V_h$, we use continuous tensor product finite elements
$Q_c^{r}$ ;see, for instance,~\cite{Ciarlet:2002:FEM:581834}.
For $Q_h$ we  use discontinuous tensor product finite elements $Q_{DG}^{r}$. 
Let $\mathcal{T}_h$ be a subdivision (triangulation) of the domain $\Omega$ into quadrilateral elements 
such that $\bigcup_{K \in \mathcal{T}_h}\overline{K}=\overline{\Omega}$
and $K \cap K'= \emptyset$ for all $K,K' \in \mathcal{T}_h$ where $K \neq K'$.
Furthermore, let $\psi_K$ be a multilinear  mapping from the reference element $\hat{K}=(0,1)^d$ to the element $K \in \mathcal{T}_h$.
We define the space $Q_{DG}^r$  as
\begin{equation*}
  Q_{DG}^{r}:=\{ v_h \in L^\infty(\Omega): v_{h|K} \in Q_r(K), \,   \forall K \in \mathcal{T}_h\},
\end{equation*}
with $Q_r(K):=\{v_{|\hat{K}}\circ\psi_K^{-1}:\, v(\hat{x})= \prod_{i=1}^{d} (\sum_{\beta=0}^{r} c_{\beta,i}\hat{x}_i^{\beta}),\, c_{\beta,i} \in \mathbb{R}\}$.
The use of these finite dimensional spaces leads to a conforming discretization for Example~\ref{Example: p-Laplace}.
We point out that the conforming discretization is needed in order to keep Theorem~\ref{Theorem: Error Representation for full Optimality Systems} valid.
The discretized abstract model problem reads as follows: Find $u_h \in U_h$ and $q_h \in Q_h$ such that they are a local solution pair of 
\begin{align}  \label{Problem: Abstract Discrete Optimization Problem}
  \begin{split}
    \min_{(u_h,q_h)} J(u_h,q_h) & \qquad  u_h \in U_h, q_h \in Q_h,\\
    \text{s.t } A(u_h,q_h)=0& \qquad \text{ in } V^*_h.
  \end{split}
\end{align}
\begin{assumption}
  There exists a unique discrete mapping $S_h\colon Q_h\mapsto U_h$, which is implicitly defined by
  \begin{equation}\label{Equation: Discrete solution operator}
    A(S_h(q_h),q_h)=0 \qquad \forall q_h \in Q_h.
  \end{equation}
  As for its continuous counterpart, we assume that it is twice
  continuously Fr{\'e}chet differentiable.
\end{assumption}
Using the discrete mapping $S_h$, we can reformulate Problem~\eqref{Problem: Abstract Discrete Optimization Problem} as the unconstrained optimization problem: Find $q_h \in Q_h$ such that it solves 
\begin{equation}
  \min_{q_h}j_h(q_h) \qquad q_h \in Q_h. \nonumber
\end{equation}
Similar to Section~\ref{Section: First Order Necessary Optimality Conditions}, we also provide the discrete version of the first order necessary optimality conditions. If $\bar q_h\in Q_h$ is a local solution, then these conditions are given by 
\begin{equation} \label{Problem: Reduced Discrete Primal Problem}
  j_h'({\bar q_h})(\delta q_h)=0 \qquad\forall\delta q_h \in Q_h.
\end{equation} We will also use the non-reduced formulation with the help of the Lagrange-approach, with 
\begin{align} 
  \mathcal{L}(u_h,q_h,z_h) := J(u_h,q_h)-A(u_h,q_h)(z_h), \qquad \forall u_h \in U_h, q_h \in Q_h, z_h \in V_h.\nonumber
\end{align}
The discrete first order necessary optimality conditions for~\eqref{Problem: Abstract Discrete Optimization Problem} are then given by
\begin{align} 
  \label{Equation: Discrete Adjoint Equation} 
  \begin{aligned}
    J'_{u}(\bar u_h,\bar q_h)(\delta u_h)-A'_{u}(\bar u_h,\bar q_h)(\bar z_h)(\delta u_h)=\mathcal{L}'_{u}(\bar u_h,\bar q_h,\bar z_h) (\delta u_h)&=0 \qquad \forall \delta u_h \in U_h, \\
    J'_{q}(\bar u_h,\bar q_h)(\delta q_h)-A'_{q}(\bar u_h,\bar q_h)(\bar z_h)(\delta q_h)=\mathcal{L}'_{q}(\bar u_h,\bar q_h,\bar z_h) (\delta q_h)&=0 \qquad \forall \delta q_h \in Q_h, \\
    -A(\bar u_h,\bar q_h)(\delta z_h)=\mathcal{L}'_{z}(\bar u_h,\bar q_h,\bar z_h)(\delta z_h) &=0  \qquad \forall \delta z_h \in V_h,
  \end{aligned}
\end{align}
where $\bar u_h=S_h(\bar q_h)$ and $\bar z_h=(S_h'(\bar q_h))^*J'_u(\bar u_h,\bar q_h)$.

\subsection{Error Representation for the Reduced System}
\label{sec_error_reduced_system}
We are now interested in an error estimator for a quantity of interest $I\colon U \times Q \mapsto \mathbb{R}$.
Let $\overline q$ be an optimal control of Problem~\eqref{eq:reducedproblem}
with associated optimal state $\bar u=S\bar (q)$. While we are interested in
$I(\overline{u},\overline{q})$, we can only compute an approximation
$I(\tilde{u}_h,\tilde{q}_h)$ of this value. Note that we assume, for most of
what follows, that $\tilde u_h:=S_h(\tilde q_h)$ is exactly solved  by means
of the solution operator  $S_h$ for the discrete state equation,
cf. Section~\ref{sec_disc}. To estimate this error, we apply the 
%well known dual weighted residual method~\cite{BeRa01} 
previously mentioned DWR method (e.g., \cite{BeRa01})
to the first order optimality conditions of our reduced system. 

Defining $i(q):=I(S(q),q)$ as well as $i_h(q):=I(S_h(q),q)$, the error between $I(S(\bar q),\bar q)$ and $I(S_h(\tilde q_h),\tilde q_h)$ can be split into
\begin{equation}
  I(S(\bar q),\bar q)-I(S_h(\tilde q_h),\tilde q_h)=i(\bar q)-i(\tilde q_h)+i(\tilde q_h)-i_h(\tilde q_h).\nonumber
\end{equation}
Therefore, $i_h$ still corresponds to our "true" quantity of interest, but computed with the discrete solutions $\tilde q_h$ and $S_h(\tilde q_h)$. 
We start by estimating the first part of the error, which actually has a practical relevance: if some approximate control $\tilde q_h$ is computed and applied in a practical situation, then the corresponding physical system will produce a "true" state $\tilde u:=S(\tilde q_h)$ instead of an approximation $\tilde u_h=S_h(\tilde q_h)$.

As a first result, we formulate a theoretical error estimator, where we need the adjoint problem to the first order optimality conditions, which is given by: Find $\overline{p} \in Q$ such that
\begin{equation} \label{Problem: Reduced Adjoint Problem}
  j''(\overline{q})(\delta q,\overline{p})=-i'(\overline{q}) (\delta q) \qquad \forall \delta q \in Q.
\end{equation}
\begin{assumption}
  We assume that~\eqref{Problem: Reduced Adjoint Problem} has a unique solution.
\end{assumption}

\begin{theorem}[Error Representation for Reduced System]\label{Theorem: Error Representation}
  Let us assume that $j \in \mathcal{C}^4(Q ,\mathbb{R})$ and $i \in \mathcal{C}^3(Q ,\mathbb{R})$. 
  If $\overline{q}$ solves~\eqref{eq:reducedproblem} and $\overline{p}$ solves~\eqref{Problem: Reduced Adjoint Problem}
  for $ \overline{q}\in Q$, 
  then, for arbitrary fixed  $\tilde{q}_h \in Q$ and $ \tilde{p}_h \in Q$, we find:
  \begin{align*} 
    \begin{aligned}
      i(\overline{q})-i(\tilde{q}_h)&= \frac{1}{2}\rho(\tilde{q}_h)(\overline{p}-\tilde{p}_h)+\frac{1}{2}\rho^*(\tilde{q}_h,\tilde{p}_h)(\overline{q}-\tilde{q}_h) 
      +\rho (\tilde{q}_h)(\tilde{p}_h) + \mathcal{R}^{(3)},
    \end{aligned}
  \end{align*}
  where
  \begin{align}
    \label{Error Estimator: primal}
    \rho(\tilde{q}_h)(\cdot) &:= j'(\tilde{q}_h)(\cdot), \\ 
    \rho^*(\tilde{q}_h,\tilde{p}_h)(\cdot) &:= i'(\tilde{q}_h)(\cdot)+j''(\tilde{q}_h)(\cdot,\tilde{p}_h), 	\nonumber
  \end{align}
  and the remainder term satisfies
  \begin{equation}
    \begin{aligned}	\label{Error Estimator: Remainderterm}
      \mathcal{R}^{(3)}:=\frac{1}{2}\int_{0}^{1}[i'''(\tilde{q}_h+se)(e,e,e)
      +j''''(\tilde{q}_h+se)(e,e,e,\tilde{p}_h+se^*)
      +3j'''(\tilde{q}_h+se)(e,e,e^*)]s(s-1)\,ds,
    \end{aligned} \nonumber
  \end{equation}
  with $e:=\overline{q}-\tilde{q}_h$ and $e^* :=\overline{p}-\tilde{p}_h$.
  \begin{proof}%
    The proof follows the same idea as in~\cite{RanVi2013,EnLaWi18} but is stated for completeness of presentation.
    Define $e$ and $e^*$ as above and let $\overline{x}$, $x$, $\tilde{x}_h$  be defined as $\overline{x}:=(\overline{q},\overline{p})$, ${x}:=({q},{p})$, $\tilde{x}_h:=(\tilde{q}_h,\tilde{p}_h),$ as well as $m(x):= i(q)+j'(q)(p)$. Furthermore, let $e_x$ be defined $e_x:= \overline{x}-\tilde{x}_h$.
    By the fundamental theorem of calculus as well as the trapezoidal rule, we observe that
    \begin{equation}\label{eq:ansatz}\begin{aligned}
        m(\overline{x})-m(\tilde{x}_h)=& \int\limits_0^1m'(\tilde{x}_h+ s e_x)(e_x) \text{ d}s \\
        =& \frac{1}{2}  \big( m'(\tilde{x}_h)(e_x) d s + m'(\overline{x})(e_x)\big) +\frac{1}{2}\int\limits_0^1m'''(\tilde x_h+se_x)(e_x,e_x,e_x)s(s-1)\text{ d}s.
      \end{aligned}\end{equation}
    By carefully inspecting $\frac{1}{2}\int\limits_0^1m'''(\tilde x_h+se_x)(e_x,e_x,e_x)s(s-1)\text{ d}s$, it follows that it coincides with $\mathcal{R}^{(3)}$. Additionally, we can deduce that
    \begin{equation}
      \label{eq:ansatz2}
      m'(\overline{x})(e_x)=i'(\overline{q})+j''(\overline{q})(e,\overline{p})+j'(\overline{q})(e^*)=0\end{equation}
    due to~\eqref{Problem: Reduced Primal Problem} and~\eqref{Problem: Reduced Adjoint Problem}.
    Combining~\eqref{eq:ansatz} and~\eqref{eq:ansatz2} results in the following identity
    \begin{equation}
      \label{eq:mdiff}
      m(\overline{x})-m(\tilde{x}_h) = \frac{1}{2}  m'(\tilde{x}_h)(e_x) + \mathcal{R}^{(3)}.\end{equation}
    Therefore,  using again~\eqref{Problem: Reduced Primal Problem} as well as~\eqref{Error Estimator: primal}, we get
    \begin{align*}
      i(\overline{q})-i(\tilde{q}_h)=&m(\overline{x})-j'(\overline{q})(\overline{q})-m(\tilde{x}_h)+{j'(\tilde{q}_h)(\tilde{p}_h)}\\
      =&m(\bar x)-m(\tilde x_h)+\rho (\tilde{q}_h)(\tilde{p}_h)\\
      =&\frac{1}{2}  m'(\tilde{x}_h)(e_x) + \mathcal{R}^{(3)}+\rho (\tilde{q}_h)(\tilde{p}_h),
    \end{align*}
    where we have applied~\eqref{eq:mdiff}.
    This proves the theorem after verifying that $m'(\tilde{x}_h)(e_x)= \rho(\tilde{q}_h)(\overline{p}-\tilde{p}_h)+\rho^*(\tilde{q}_h,\tilde{p}_h)(\overline{q}-\tilde{q}_h) $.
  \end{proof}
\end{theorem}
\begin{remark}
  One objective of this representation, in addition to the fact that for instance $\tilde u=S(\tilde q_h)$ is not readily available exactly, is to obtain indicators for local adaptivity. By inspecting the primal part of the error estimator $\rho(\tilde{q}_h)(\overline{p}-\tilde{p}_h)$, we observe that 
  \begin{align*}
    \rho(\tilde{q}_h)(\overline{p}-\tilde{p}_h)=-j'(\tilde{q}_h)(\overline{p}-\tilde{p}_h)=- J'_q(S(\tilde{q}_h),\tilde{q}_h)(\overline{p}-\tilde{p}_h)- J'_u(S(\tilde{q}_h),\tilde{q}_h)(S'(\tilde{q}_h)(\overline{p}-\tilde{p}_h)).
  \end{align*}
  Since it is not clear how to localize $S'(\tilde{q}_h)(\overline{p}-\tilde{p}_h)$, we do not follow this path to compute the error indicators, but prove a localizable error estimator in a
  similar fashion in Theorem \ref{Theorem: Error Representation for full Optimality Systems}, which makes use of \eqref{Equation: First Order Opt u} as well.
%   A similar argument holds for the second part of the error estimator $\rho^*(\tilde{q}_h,\tilde{p}_h)(\overline{q}-\tilde{q}_h) $.
\end{remark}

For another idea, we consider the adjoint problem  to  the first order optimality conditions for the Lagrangian defined in~\eqref{Lagrangian}: Find $(\overline{v},\overline{p_2},\overline{y} )\in U \times Q \times V$ such that
\begin{align} \label{Problem:  Adjoint Optimality System}
  \begin{pmatrix}
    \mathcal{L}''_{uu} &\mathcal{L}''_{uq} &\mathcal{L}''_{uz} \\
    \mathcal{L}''_{qu} &\mathcal{L}''_{qq} &\mathcal{L}''_{qz} \\
    \mathcal{L}''_{zu} &\mathcal{L}''_{zq} &0
  \end{pmatrix} 		\begin{pmatrix}
    \overline{v}\\
    \overline{p_2}\\
    \overline{y}
  \end{pmatrix} = -
  \begin{pmatrix}
    I'_u\\
    I'_q\\
    0
  \end{pmatrix}\qquad \text{ in }  U^* \times Q^* \times V^*,
\end{align}
where the argument in the partial derivatives is always given by $(\overline{u},\overline{q},\overline{z})$.
\begin{assumption}
  We assume that~\eqref{Problem:  Adjoint Optimality System} has a unique solution.
\end{assumption}
In order to obtain the variables $\overline{v}$ and $\overline{y}$ with the help of the solution of the reduced adjoint problem~\eqref{Problem: Reduced Adjoint Problem}, the following lemma is useful.
\begin{lemma} \label{Lemma: p=p}
  If $\bar q\in Q$ with associated state $\bar u=S(\bar q)$ is a
  local solution of~\eqref{eq:reducedproblem}, and $\bar p$
  solves~\eqref{Problem: Reduced Adjoint Problem}, then $\bar
  v=S'(\bar q)\bar p$, $\bar p_2=\bar p$, and $\bar y$ given
  by~\eqref{Proof: y= ...} solve~\eqref{Problem:  Adjoint Optimality System}.
  \begin{proof}
    Let $\mathfrak{p} \in Q $ be arbitrary.
    Using the definition of the reduced functionals, we obtain
    \begin{align}\label{Equation: Proof q=q j''}
      \begin{aligned}
        j''(\overline{q})(\mathfrak{p})=&\;J''_{uu}(S(\overline{q}),\overline{q})(S'(\overline{q})\mathfrak{p})\circ
        S'(\overline{q})+J''_{uq}(S(\overline{q}),\overline{q})(\mathfrak{p})\circ
        S'(\overline{q})\\&+J''_{qu}(S(\overline{q}),\overline{q})(S'(\overline{q})\mathfrak{p})+J''_{qq}(S(\overline{q}),\overline{q})(\mathfrak{p})+J'_u(S(\overline{q}),\overline{q})\circ
        S''(\overline{q})(\mathfrak{p}),
      \end{aligned}
    \end{align}
    and 
    \begin{align*}
      \begin{aligned}
        i'(\overline{q}) = I'_{u}(S(\overline{q}),\overline{q})S'(\overline{q})+I'_{q}(S(\overline{q}),\overline{q}).
      \end{aligned}
    \end{align*}
    Furthermore, with the definition of the solution operator, we obtain from~\eqref{Equation: Solution Operator} that
    \begin{align}\label{Equation: Proof q=q 1. 00}
      \begin{aligned}
        0=&\;A'_{u}(S(\overline{q}),\overline{q})(S'(\overline{q})\mathfrak{p})+A'_{q}(S(\overline{q}),\overline{q})(\mathfrak{p})=\mathcal{L}''_{zu}(S(\overline{q}),\overline{q},\overline{z})(S'(\overline{q})\mathfrak{p})+\mathcal{L}''_{zq}(S(\overline{q}),\overline{q},\overline{z})(\mathfrak{p})
      \end{aligned}
    \end{align}
    and
    \begin{align}\label{Equation: Proof q=q 1. 0}
      \begin{aligned}
        0=&\;A''_{uu}(S(\overline{q}),\overline{q})(S'(\overline{q})\mathfrak{p})\circ
        S'(\overline{q})+A''_{uq}(S(\overline{q}),\overline{q})(\mathfrak{p})\circ
        S'(\overline{q})\\&+A''_{qu}(S(\overline{q}),\overline{q})(S'(\overline{q})\mathfrak{p})
        +A''_{qq}(S(\overline{q}),\overline{q})(\mathfrak{p})+A'_u(S(\overline{q}),\overline{q})\circ
        S''(\overline{q})(\mathfrak{p}).
      \end{aligned}
    \end{align}
    By subtracting~\eqref{Equation: Proof q=q
      1. 0} from~\eqref{Equation: Proof q=q j''}, it follows that
    \begin{align}
      \begin{aligned}
        j''(\overline{q})(\mathfrak{p})=&\;
        j''(\overline{q})(\mathfrak{p})-0\\
        =&\;\mathcal{L}''_{uu}(S(\overline{q}),\overline{q},\overline{z})(S'(\overline{q})\mathfrak{p})\circ
        S'(\overline{q})+\mathcal{L}''_{uq}(S(\overline{q}),\overline{q},\overline{z})(\mathfrak{p})\circ
        S'(\overline{q})\\&+\mathcal{L}''_{qu}(S(\overline{q}),\overline{q},\overline{z})(S'(\overline{q})\mathfrak{p})
        +\mathcal{L}''_{qq}(S(\overline{q}),\overline{q},\overline{z})(\mathfrak{p})+\mathcal{L}'_u(S(\overline{q}),\overline{q},\overline{z})\circ
        S''(\overline{q})(\mathfrak{p}).
      \end{aligned}\nonumber
    \end{align}
    Further, from~\eqref{Equation:
      Proof q=q 1. 00}
    we get
    \[
      S'(\overline{q}) = -[\mathcal{L}''_{zu}(S(\overline{q}),\overline{q})]^{-1}\mathcal{L}''_{zq}(S(\overline{q}),\overline{q})).
    \]
    Thus $\overline{p}_2 =
    \overline{p}$ and
    $\overline{v}=S'(\overline{q})\overline{p}$
    satisfy the third line in~\eqref{Problem:  Adjoint Optimality System}.
    
    To proceed, we note that $\overline{q}$,
    $\overline{u}=S(\overline{q})$ and
    $\overline{z}$ solves~\eqref{Equation: First
      Order Opt u}, thus we have that $\mathcal{L}'_u(S(\overline{q}),\overline{q},\overline{z})=0$.
    This leads to 
    \begin{align}
      \label{Proof: q=q SSLuu+SLuq +Luq+lqq}
      \begin{aligned}
        j''(\overline{q})(\overline{p}_2)
        =&\;\mathcal{L}''_{uu}(S(\overline{q}),\overline{q},\overline{z})(\overline{v})\circ
        S'(\overline{q})+\mathcal{L}''_{uq}(S(\overline{q}),\overline{q},\overline{z})(\overline{p}_2)\circ
        S'(\overline{q})\\&+\mathcal{L}''_{qu}(S(\overline{q}),\overline{q},\overline{z})(\overline{v})+\mathcal{L}''_{qq}(S(\overline{q}),\overline{q},\overline{z})(\overline{p}_2)
      \end{aligned}
    \end{align}
    Now, we define $\overline{y}$ by the first
    line of~\eqref{Problem:
      Adjoint Optimality System}, we get
    \begin{align}\label{Proof: y= ...}
      \begin{aligned}
        \mathcal{L}''_{uz}(S(\overline{q}),\overline{q}))(\overline{y})=-\big(\mathcal{L}''_{uu}(S(\overline{q}),\overline{q},\overline{z})(\overline{v})+\mathcal{L}''_{qu}(S(\overline{q}),\overline{q},\overline{z})(\overline{p}_2)+I'_{u}(S(\overline{q}),\overline{q})\big).
      \end{aligned}
    \end{align}
    With this, we can rewrite~\eqref{Proof: q=q SSLuu+SLuq +Luq+lqq} as
    \begin{align*}                                                                                                        \begin{aligned}
        j''(\overline{q})(\overline{p}_2)
        =&\;\Bigl(\mathcal{L}''_{uu}(S(\overline{q}),\overline{q},\overline{z})(\overline{v})+\mathcal{L}''_{uq}(S(\overline{q}),\overline{q},\overline{z})(\overline{p}_2)\Bigr)\circ
        S'(\overline{q})\\&+\mathcal{L}''_{qu}(S(\overline{q}),\overline{q},\overline{z})(\overline{v})+\mathcal{L}''_{qq}(S(\overline{q}),\overline{q},\overline{z})(\overline{p}_2)\\
        =&\;
        -\Bigl(\mathcal{L}''_{uz}(S(\overline{q}),\overline{q}))(\overline{y})
        + I'_{u}(S(\overline{q}),\overline{q})\Bigr)\circ
        S'(\overline{q})\\&+\mathcal{L}''_{qu}(S(\overline{q}),\overline{q},\overline{z})(\overline{v})+\mathcal{L}''_{qq}(S(\overline{q}),\overline{q},\overline{z})(\overline{p}_2).                               \end{aligned}
    \end{align*}
    Now, we can use the definition of $\overline{p}$,
    $\mathcal{L}''_{uz} = (\mathcal{L}''_{zu})*, \mathcal{L}''_{zq} =
    (\mathcal{L}''_{qz})*$, 
    the formula for $S'(\overline{q})$ and the
    representation of $i'(\overline{q})$ to get
    \begin{align*}
      -I'_{u}(S(\overline{q}),\overline{q})S'(\overline{q})&-I'_{q}(S(\overline{q}),\overline{q})
      =-i'(\overline{q})\\
      &=j''(\overline{q})(\overline{p}_2)\\
      &= -I'_{u}(S(\overline{q}),\overline{q})\circ S'(\overline{q})
      - S'(\overline{q})^*
      \mathcal{L}''_{uz}(S(\overline{q}),\overline{q}))(\overline{y})
      \\
      &\;\;\;\;+\mathcal{L}''_{qu}(S(\overline{q}),\overline{q},\overline{z})(\overline{v})+\mathcal{L}''_{qq}(S(\overline{q}),\overline{q},\overline{z})(\overline{p}_2)\\
      &= -I'_{u}(S(\overline{q}),\overline{q})\circ S'(\overline{q})
      +
      \mathcal{L}''_{qz}(S(\overline{q}),\overline{q}))(\overline{y})
      \\
      &\;\;\;\;+\mathcal{L}''_{qu}(S(\overline{q}),\overline{q},\overline{z})(\overline{v})+\mathcal{L}''_{qq}(S(\overline{q}),\overline{q},\overline{z})(\overline{p}_2)
    \end{align*}
    and the second line in~\eqref{Problem:  Adjoint Optimality System} follows.
  \end{proof}
\end{lemma}
Lemma~\ref{Lemma: p=p} allows to obtain $\overline{p}=\bar p_2$ by solving the reduced adjoint equation~\eqref{Problem: Reduced Adjoint Problem}. Then, $\overline{v}$  can be computed by solving the tangent equation 
\begin{align*}
  \mathcal{L}''_{zu}(\overline{u},\overline{q},\overline{z})( \cdot ,\overline{v}) +&\mathcal{L}''_{zq} (\overline{u},\overline{q},\overline{z})( \cdot ,\overline{p})=0,
\end{align*}
which is the last row of~\eqref{Problem:  Adjoint Optimality System}. Using this solution, we can deduce $\overline{y}$ from the first row of~\eqref{Problem:  Adjoint Optimality System}.

An analogue to~\eqref{Problem:  Adjoint Optimality System} on the discrete level is given by: Find $(\tilde v_h,\tilde p_h,\tilde y_h )\in U_h \times Q_h \times V_h$
such that
\begin{align} \label{Problem: Discrete Adjoint Optimality System}
  \begin{pmatrix}
    \mathcal{L}''_{uu} &\mathcal{L}''_{uq} &\mathcal{L}''_{uz} \\
    \mathcal{L}''_{qu} &\mathcal{L}''_{qq} &\mathcal{L}''_{qz} \\
    \mathcal{L}''_{zu} &\mathcal{L}''_{zq} &0
  \end{pmatrix} 		\begin{pmatrix}
    \tilde v_h\\
    \tilde p_h\\
    \tilde y_h
  \end{pmatrix} = -
  \begin{pmatrix}
    I'_u\\
    I'_q\\
    0
  \end{pmatrix},
\end{align}
where the arguments in the partial derivatives are given by $(\tilde u_h,\tilde q_h,\tilde z_h)$.
\begin{remark}\label{rem:dualvariables}
  If~\eqref{Problem: Discrete Adjoint Optimality System} is considered at the linearization point $\bar q_h,\bar u_h,\bar z_h$, then Lemma~\ref{Lemma: p=p} holds also true for the discrete problem, i.e.  if $\overline{p}_h \in Q_h$ solves
  \begin{equation} \label{Problem: Reduced Discrete Adjoint Problem}
    j_h''(\overline{q}_h)(\delta q_h,\overline{p}_h)=-i_h'(\overline{q}_h) (\delta q_h) \qquad \forall \delta q_h \in Q_h,
  \end{equation} then $\tilde p_h=\bar p_h$. This can be shown by the same proof replacing $S$ by $S_h$.\end{remark}
Similar as explained above, the variables $\tilde v_h$ and $\tilde y_h$ can be deduced from the knowledge of $\bar p_h$ and the discrete version of Lemma~\ref{Lemma: p=p}.

\begin{theorem}[Localizable Error Representation for Reduced System]\label{Theorem: Error Representation for full Optimality Systems}
  Let us assume that $j \in \mathcal{C}^4(Q ,\mathbb{R})$ and $i \in \mathcal{C}^3(Q ,\mathbb{R})$. 
  Let $\overline{q}$ be a local solution
  of~\eqref{eq:reducedproblem}, with
  $\overline{\xi}=(\overline{u},\overline{q},\overline{z})$ the
  corresponding KKT-triplet given by~\eqref{Equation: First
    Order Opt u}, and let the triple $\overline{\xi}^*=(\bar
  v,\overline{p},\bar y)\in U\times Q\times V$
  solve~\eqref{Problem:  Adjoint Optimality System}. Moreover,
  let $\tilde{q}_h \in Q_h$ be an arbitrary fixed discrete
  control, and let $\tilde{\xi}_h^* = (\tilde p_h,\tilde
  v_h,\tilde y_h)$ be the solution
  to~\eqref{Equation: Discrete Adjoint Equation} and the first
  and last row of~\eqref{Problem: Discrete Adjoint Optimality
    System} at the linearization point $\tilde{\xi}_h = (\tilde u_h, \tilde
  q_h,\tilde z_h)$ with $\tilde u_h=S_h(\tilde q_h)$ and $\tilde
  z_h=S'_h(\tilde q_h)^*J'_u(\tilde u_h,\tilde q_h)$.
  Then we have the error representation 
  \begin{align} \label{Equation: Error Representation for full Optimality Systems}
    \begin{aligned}
      i(\overline{q})-i_h(\tilde{q}_h)=&\; \frac{1}{2}\big[\rho_u(\tilde{\xi}_h,\tilde{\xi}_h^*)(\overline{y}-\tilde{y}_h) + \rho_z(\tilde{\xi}_h,\tilde{\xi}_h^*)(\overline{v}-\tilde{v}_h) +\rho_q(\tilde{\xi}_h,\tilde{\xi}_h^*)(\overline{p}-\tilde{p}_h)\\			
      &\;\;\;\;\;+\rho_v(\tilde{\xi}_h,\tilde{\xi}_h^*)(\overline{z}-\tilde{z}_h) + \rho_y(\tilde{\xi}_h,\tilde{\xi}_h^*)(\overline{u}-\tilde{u}_h)	  +\rho_p(\tilde{\xi}_h,\tilde{\xi}_h^*)(\overline{q}-\tilde{q}_h)\big]\\
      &-j'_h (\tilde{q}_h)(\tilde{p}_h) + \tilde{\mathcal{R}}^{(3)},
    \end{aligned}
  \end{align}
  where
  \begin{align*}
    \rho_u(\tilde{\xi}_h,\tilde{\xi}_h^*)(\cdot):=&\mathcal{L}'_{z}(\tilde{\xi}_h)(\cdot),  \\
    \rho_q(\tilde{\xi}_h,\tilde{\xi}_h^*)(\cdot):=&\mathcal{L}'_{q}(\tilde{\xi}_h)(\cdot),  \\
    \rho_z(\tilde{\xi}_h,\tilde{\xi}_h^*)(\cdot):=&\mathcal{L}'_{u}(\tilde{\xi}_h)(\cdot) , \\
    \rho_v(\tilde{\xi}_h,\tilde{\xi}_h^*)(\cdot):=&\mathcal{L}''_{zu}(\tilde{\xi}_h)(\cdot,\tilde{v}_h)+\mathcal{L}''_{zq}(\tilde{\xi}_h)(\cdot,\tilde{p}_h),  \\
    \rho_p(\tilde{\xi}_h,\tilde{\xi}_h^*)(\cdot):=&I'_{q} (\tilde{u}_h,\tilde{q}_h)(\cdot)+ \mathcal{L}''_{uq}(\tilde{\xi}_h)(\tilde{v}_h,\cdot)+\mathcal{L}''_{qq}(\tilde{\xi}_h)(\tilde{p}_h,\cdot)+\mathcal{L}''_{zq}(\tilde{\xi}_h)(\tilde{y}_h,\cdot) , \\
    \rho_y(\tilde{\xi}_h,\tilde{\xi}_h^*)(\cdot):=&I'_{u} (\tilde{u}_h,\tilde{q}_h)(\cdot)+ \mathcal{L}''_{uu}(\tilde{\xi}_h)(\tilde{v}_h,\cdot)+\mathcal{L}''_{qu}(\tilde{\xi}_h)(\tilde{p}_h,\cdot)+\mathcal{L}''_{zu}(\tilde{\xi}_h)(\tilde{y}_h,\cdot) , 
  \end{align*}
  and the remainder term $$\tilde{\mathcal{R}}^{(3)}=\frac{1}{2}\int\limits_{0}^{1}[I'''(\tilde{\xi}+s\tilde{e}_\xi)(\tilde{e}_\xi,\tilde{e}_\xi,\tilde{e}_\xi)
  +\mathcal{L}''''(\tilde{\xi}+s\tilde{e}_\xi)(\tilde{e}_\xi,\tilde{e}_\xi,\tilde{e}_\xi,\tilde{p}+s\tilde{e}_\xi^*)
  +3\mathcal{L}'''(\tilde{\xi}+s\tilde{e}_\xi)(\tilde{e}_\xi,\tilde{e}_\xi,\tilde{e}_\xi^*)]s(s-1)\,ds,$$ 
  with $\tilde{e}_\xi = \overline{\xi}-\tilde{\xi}_h$, $\tilde{e}_\xi^*=\overline{\xi}^*-\tilde{\xi}_h^*$.
  \begin{proof}
    %%%%%%%%%%%%%%%%%%%%%%%%%%%%%%%%%%%%%%%%%%%%%%%%%%%%%%%%%%%%%%%%%%%%%%%%%%%%%%%% 
    The proof follows a similar structure as the proof of Theorem~\ref{Theorem: Error Representation}.
    Let
    $\overline{x}:=(\overline{\xi},\overline{\xi}^*)$,
    $\tilde{x}_h:=(\tilde{\xi}_h,\tilde{\xi}^*_h)$. For
    $x = (\xi,\xi^*) = (u,q,z,\xi^*)$ we define 
    $\mathcal{M}(x):= I(\xi
    )+\mathcal{L}'(\xi)(\xi^*) = I(u,q)+\mathcal{L}'(\xi)(\xi^*)$.
    It holds that
    \begin{align} \label{Equation: Trapez of M}
      \begin{aligned}
        \mathcal{M}(\overline{x})-\mathcal{M}(\tilde{x}_h)=& \int\limits_0^1\mathcal{M}'(\tilde{x}_h+ s e_x)(e_x) \text{ d}s \\
        =& \frac{1}{2}  \big( \mathcal{M}'(\tilde{x}_h)(e_x) d s + \mathcal{M}'(\overline{x})(e_x)\big) +\frac{1}{2}\int\limits_0^1\mathcal{M}'''(\tilde x_h+se_x)(e_x,e_x,e_x)s(s-1)\text{ d}s
      \end{aligned}
    \end{align}
    where $e_x = \overline{x}-\tilde{x}_h$.			
    By carefully inspecting $\mathcal{M}'''(\tilde x_h+se_x)(e_x,e_x,e_x)$ it follows that
    \begin{align*}
      \mathcal{M}'''(\tilde{x}_h+se_x)(e_x,e_x,e_x)
      &=(\mathcal{M}'''_{\xi \xi \xi } +
      3\mathcal{M}'''_{\xi \xi \xi^*}
      +3\mathcal{M}'''_{\xi \xi^* \xi^*}
      +\mathcal{M}'''_{\xi^* \xi^* \xi^*
      })(\tilde{x}_h+se_x)(e_x,e_x,e_x)\\
      &=I'''(\tilde{\xi }_h+s\tilde{e}_{\xi})(\tilde{e}_{\xi},\tilde{e}_{\xi},\tilde{e}_{\xi})
      +\mathcal{L}''''(\tilde{\xi
      }_h+s\tilde{e}_{\xi})(\tilde{e}_{\xi},\tilde{e}_{\xi},\tilde{e}_{\xi},\tilde{\xi}^*_h+s\tilde{e}_{\xi}^*)\\
      &\;\;\;\;+3\mathcal{L}'''(\tilde{\xi }_h +s(\tilde{e}_{\xi}))(\tilde{e}_{\xi},\tilde{e}_{\xi},\tilde{e}_{\xi}^*)
    \end{align*}
    since $\mathcal{M}''_{\xi^* \xi^* }=0$ and
    $\mathcal{M}'_{\xi^*} (\tilde{x}_h+se_x)(e_x)=\mathcal{L}'(\tilde{\xi }_h
    +s(\tilde{e}_{\xi}))(\tilde{e}_{\xi}^*)$.
    Thus,~\eqref{Equation: Trapez of M} gives  
    \begin{align*}
      \begin{aligned}
        \mathcal{M}(\overline{x})-\mathcal{M}(\tilde{x}_h)
        =& \frac{1}{2}  \big( \mathcal{M}'(\tilde{x}_h)(e_x) d s + \mathcal{M}'(\overline{x})(e_x)\big) 
        + \tilde{\mathcal{R}}^{(3)}.
      \end{aligned} 
    \end{align*}	
    For the part $\mathcal{M}'(\overline{x})(e_x)$ of~\eqref{Equation: Trapez of M}, we can deduce that
    $$\mathcal{M}'(\overline{x})(e_x)=I'(\overline{\xi })+\mathcal{L}''(\overline{\xi })(\tilde{e}_{\xi},\overline{\xi }^*)+\mathcal{L}'(\overline{\xi })(\tilde{e}_{\xi}^*)=0,$$
    since $\overline{\xi}^*$ solves \eqref{Problem:  Adjoint Optimality System} and $\overline{\xi}$ solves \eqref{Equation: First Order Opt u}.
    Finally, relation~\eqref{Equation: Trapez of M} reduces to  the following identity
    $$\mathcal{M}(\overline{x})-\mathcal{M}(\tilde{x}_h) = \frac{1}{2}  \mathcal{M}'(\tilde{x}_h)(e_x) + \tilde{\mathcal{R}}^{(3)}.$$
    Therefore, we get
    \begin{align*}
      I(\overline{\xi })-I(\tilde{\xi }_h)=\mathcal{M}(\overline{x})-\mathcal{L}'(\overline{\xi })(\overline{\xi })-\mathcal{M}(\tilde{x}_h)+\mathcal{L}'(\tilde{\xi }_h)(\tilde{\xi}^*_h)= \frac{1}{2}  \mathcal{M}'(\tilde{x}_h)(e_x) + \tilde{\mathcal{R}}^{(3)}+\mathcal{L}'(\tilde{\xi }_h)(\tilde{\xi}^*_h).
    \end{align*}
    Furthermore, we can deduce that $\mathcal{M}'(\tilde{x}_h)(e_x)= \mathcal{L}'(\tilde{\xi}_h)(\tilde{e}_{\xi}^*) 
    + I'(\tilde{\xi}_h)(\tilde{e}_\xi)+\mathcal{L}''(\tilde{\xi}_h)(\tilde{e}_\xi,\tilde{\xi}^*_h)$.
    Gathering the results from above, we obtain,
    noting that $\tilde{\xi}_h =
    (\tilde{u}_h,\tilde{q}_h,\tilde{z}_h) = (S_h(\tilde{q}_h),\tilde{q}_h,\tilde{z}_h)$
    \begin{align}\label{Proof: Representation hat}
      \begin{aligned}
        i(\overline{q}) -
        i_h(\tilde{q}_h)
        &= 
        I(\overline{u},\overline{q})-I(\tilde{u}_h,\tilde{q}_h)\\
        &= \frac{1}{2}\big[\mathcal{L}'(\tilde{\xi}_h)(\tilde{e}^*) 
        + I'(\tilde{\xi}_h)(\tilde{e}_\xi)+\mathcal{L}''(\tilde{\xi}_h)(\tilde{e}_\xi,\tilde{\xi}^*_h)
        \big]
        +\mathcal{L}'(\tilde{\xi}_h)(\tilde{\xi}_h^*)
        + \tilde{\mathcal{R}}^{(3)}.
      \end{aligned}
    \end{align}			
    Straightforward
    calculations show
    \begin{equation}
      \label{Equation: Full Error Estimator Proof Primal Part}
      \mathcal{L}'(\tilde{\xi}_h)(\tilde{e}^*)=\rho_u(\tilde{\xi}_h,\tilde{\xi}_h^*)(\overline{y}-\tilde{y}_h) + \rho_z(\tilde{\xi}_h,\tilde{\xi}_h^*)(\overline{v}-\tilde{v}_h)+\rho_q(\tilde{\xi}_h,\tilde{\xi}_h^*)(\overline{p}-\tilde{p}_h),
    \end{equation}
    and 
    \begin{equation}
      \label{Equation: Full Error Estimator Proof Adjoint Part}
      I'(\tilde{\xi}_h)(\tilde{e}_\xi)+\mathcal{L}''(\tilde{\xi}_h)(\tilde{e}_\xi,\tilde{\xi}_h^*)=\rho_v(\tilde{\xi}_h,\tilde{\xi}_h^*)(\overline{z}-\tilde{z}_h) + \rho_y(\tilde{\xi}_h,\tilde{\xi}_h^*)(\overline{u}-\tilde{u}_h)	  +\rho_p(\tilde{\xi}_h,\tilde{\xi}_h^*)(\overline{q}-\tilde{q}_h).
    \end{equation}
  \end{proof}
\end{theorem}
Let us end this section with  some further observations.
 	\begin{remark}
 		Note that if $\tilde q_h=\bar q_h$, then $(\tilde v_h,\tilde p_h=\bar p_h, \tilde y_h)$ in fact solve~\eqref{Equation: Discrete Adjoint Equation}, cf. Remark~\ref{rem:dualvariables}, and consequently $j'_h (\tilde{q}_h)(\tilde{p}_h)=0$.
 	\end{remark}

\begin{remark}
  From numerical experiments for the regularized $p$-Laplacian computed
  in~\cite{EnLaWi18b}, we can deduce that $\mathcal{R}^{(3)}$ can be neglected
  on sufficiently refined meshes.
\end{remark}
An identity also observed in \cite{Wo10}, is the following:
\begin{propositon}
  If $I =J$ and $j''(q)$ is injective, then we have $(\overline{v},\overline{p},\overline{y})=(0,0,\overline{z})$.
  \begin{proof}
    Since $J$ is the cost functional and $( \overline{u},\overline{q})$ is a local minimizer of our optimization problem the first order necessary condition is given by $j'(\overline{q})=0$. Therefore the adjoint equation reads as 
    $$j''(\overline{q}) \overline{p} =-i'(\overline{q})= -j'(\overline{q})=0.$$
    If $j''(\overline{q})$ is injective, then $\overline{p}=0$.
    From the  tangent equation
    $$	\mathcal{L}''_{zu}(\overline{u},\overline{q},\overline{z})( \cdot ,\overline{v}) +\mathcal{L}''_{zq} (\overline{u},\overline{q},\overline{z})( \cdot ,\overline{p})=0,$$
    we can deduce that $\overline{v}=0$.
    Finally the optimality system  reduces to 
    $$\mathcal{L}''_{zu} (\overline{u},\overline{q},\overline{z})(\overline{y},\cdot) + I'_u(\cdot)=0, \quad \text{ and } \quad  \mathcal{L}''_{zq} (\overline{u},\overline{q},\overline{z})(\overline{y},\cdot) + I'_q(\cdot)=0.$$
    From this follows that $\overline{y}=\overline{z}$, which completes the proof.
  \end{proof}
\end{propositon}

\subsection{The Parts of the Error Estimator}\label{Subsection: The error estimator parts}
We now briefly discuss the two main parts of the error estimator:
\[
  \eta_{h,k}^{(2)}:= \eta_k + \eta_h^{(2)},
\]
where the first part refers to the iteration error, and the second term
denotes the discretization error to be defined in the following. We recall that $\eta_{h,k}^{(2)}$ is designed to
estimate $i(\overline{q})-i_h(\tilde{q}_h)$ given in 
\eqref{Equation: Error Representation for full Optimality Systems}.

\paragraph{The iteration error estimator}
The iteration error estimator 
\[
  \eta_k:=-j'_h (\tilde{q}_h)(\tilde{p}_h)
\]
can be used as stopping rule for the nonlinear solver like for Newton's 
method as in~\cite{RanVi2013, EnLaWi18,EnLaWi18} 
and Algorithm~\ref{inexat_newton_algorithm_for_multiple_goal_functionals}
presented in Section~\ref{sec_algo}.

\paragraph{The discretization error estimator}
Of course the exact solution of the optimal control problem in formula~\eqref{Equation: Error Representation for full Optimality Systems} are not
known. They can either be replaced by a (patch-wise) higher order polynomial 
interpolation
or by approximations on enriched spaces 
\cite{BeRa01,BaRa03}.

The discretization error estimator using the solutions $(u_h^{(2)},q_h^{(2)},z_h^{(2)})$ and $(v_h^{(2)},p_h^{(2)},y_h^{(2)})$  on enriched spaces reads as 
	\begin{align}\label{Equation: Computable  error estimator parts}
		\eta_h^{(2)}:=\frac{1}{2}\big[&\rho_u(\tilde{\xi}_h,\tilde{\xi}_h^*)(y_h^{(2)}-\tilde{y}_h) + \rho_z(\tilde{\xi}_h,\tilde{\xi}_h^*)(v_h^{(2)}-\tilde{v}_h) +\rho_q(\tilde{\xi}_h,\tilde{\xi}_h^*)(p_h^{(2)}-\tilde{p}_h)\\			
		+&\rho_v(\tilde{\xi}_h,\tilde{\xi}_h^*)(z_h^{(2)}-\tilde{z}_h) + \rho_y(\tilde{\xi}_h,\tilde{\xi}_h^*)(u_h^{(2)}-\tilde{u}_h)	  +\rho_p(\tilde{\xi}_h,\tilde{\xi}_h^*)(q_h^{(2)}-\tilde{q}_h)\big]. \nonumber
		\end{align}
		The replacement is justified if a strengthened saturation
                assumption is fulfilled as shown in~\cite{EnLaWi18b} for
                both the nonlinear state equation and the goal functionals. 

We briefly recall that the localization can be performed in three ways:
classical 
integration by parts yielding the strong problem formulation~\cite{BeRa01}, a filtering
approach employing the weak problem formulation~\cite{BraackErn02}, or 
a partition-of-unity using again the weak form of the problem~\cite{RiWi15_dwr}. All three techniques are analyzed (theoretically and computationally) with respect to
their effectivity in~\cite{RiWi15_dwr}. 
In the theoretical analysis, a discrete version of Lemma~\ref{Lemma: p=p} is
necessary to justify that $(v_h^{(2)},p_h^{(2)},y_h^{(2)})$ is indeed a
solution in the enriched spaces.

%%%%%%%%%%%%%%%%%%%%%%%%%%%%%%%%%%%%%%%%%%%%%%%%%%%%%%%%%%%%%%%%%%%%%%%%%% 
\section{Extension to Multiple Goal Functionals}
\label{sec_multigoal}
In Section~\ref{Section: DWR for Reduced System}, we discussed how the DWR
method works for one functional. However, for some problems, several functional
evaluations would be of interest. Let us consider $N$ goal
functionals $I_1, I_2,\ldots, I_N$ for some $N \in \mathbb{N}$. One possibility
would be to compute the error estimators  separately as described in Section~\ref{Section:
  DWR for Reduced System}. However, we would have to solve the adjoint problem
$N$ times, leading to high computational cost. There are several ways to
tackle this problem as for example discussed in
\cite{HaHou03,Ha08,PARDO20101953,AlvParBar2013} and more recently in 
\cite{KerPruChaLaf2017,EnWi17,EnLaWi18,EnLaWiPAMM18,EnLaWi18b}. 

Adopting the techniques presented in~\cite{EnLaWi18}, we try to combine the functionals to one, and apply the DWR method for one functional to it. In the following section, we consider $\overline{u}$, $\overline{q}$ as the solution of~\eqref{Problem: Abstract Optimization Problem}, and $\tilde{u}_h$, $\tilde{q}_h$ as some approximations. To construct the combination, we introduce a so called error weighting function:

\begin{definition}[Error weighting function~\cite{EnLaWi18}]
  Let $ M \subseteq \mathbb{R}^N$. 
  We say that $\mathfrak{E}: (\mathbb{R}^+_0)^N \times  M \mapsto \mathbb{R}^+_0$ is an \textit{error-weighting
    function} if  $\mathfrak{E}(\cdot,m) \in
  \mathcal{C}^1((\mathbb{R}^+_0)^N,\mathbb{R}^+_0)$ is strictly monotonically
  increasing in each component and $\mathfrak{E} (0,m)=0$ for all
  $m \in M$.
\end{definition}

As in~\cite{EnLaWi18}, let $\vec{I}(\cdot):=(I_1(\cdot),I_2(\cdot),\ldots,I_N(\cdot))$ mapping from $\bigcap_{i=1}^N\mathcal{D}(I_i) \subset U \times Q \mapsto \mathbb{R}^N$. Furthermore, we define $|\cdot|_N: \mathbb{R}^N \mapsto (\mathbb{R}_0^+)^N$ as the component-wise absolute value.
This allows us to construct the error function $I_\mathfrak{E}$ as follows: 
\begin{align}\label{Equation: ExactErrorweightingFunction}
	\tilde{I}_\mathfrak{E}(\cdot):=\mathfrak{E}(|\vec{I}(\overline{u},\overline{q})-\vec{I}(\cdot)|_N,\vec{I}(\tilde{u}_h,\tilde{q}_h)).
\end{align}
\begin{remark}
  The error functional $\tilde{I}_\mathfrak{E}$ is constructed in a way, that avoids error cancellation between two or more functionals. For a more detailed discussion, we refer the reader to~\cite{EnLaWi18,EnLaWi18b}.
\end{remark}
\begin{remark}
  The quantity~\eqref{Equation: ExactErrorweightingFunction} is not computable, since it depends on $\vec{I}(\overline{u},\overline{q})$, which is not known. However, we can use a higher order polynomial approximation to approximate this quantity, as done in~\cite{HaHou03,EnLaWi18,EnLaWi18b}, where  consequences of the replacement are discussed in~\cite{EnLaWi18b}. 
\end{remark}
The resulting error weighting functional is given by
\begin{equation}
  I_\mathfrak{E}(\cdot):=\mathfrak{E}(|\vec{I}(u_h^{(2)},q_h^{(2)})-\vec{I}(\cdot)|_N,\vec{I}(\tilde{u}_h,\tilde{q}_h)),
\end{equation}
where $u_h^{(2)}$, $q_h^{(2)}$ denote the solutions on enriched finite element spaces. 
\begin{remark}
  We notice that, for the choice $\mathfrak{E}(x,m):= \sum_{\ell=1}^{N}\frac{x_\ell}{|m_\ell|}$, we obtain the same combined functional as in~\cite{EnWi17} up to sign. The same holds for~\cite{HaHou03,Ha08} in the case of linear problems. This choice is used in our numerical examples.
\end{remark}

\begin{remark}
  Finally, the method explained in Section~\ref{sec_DWR} is applied to $I_\mathfrak{E}$ instead of $I$ to achieve a control of the errors in all functionals at once, as algorithmically illustrated in Section~\ref{sec_algo}.
\end{remark}

% ##############################################################################%
\section{Algorithmic Details}
\label{sec_algo}
In this section, we briefly recapitulate the algorithmic techniques to solve 
the optimal control problem with multiple goal functionals that we have
outlined in the previous sections. The algorithms for the forward 
problem including multiple goal functionals evaluations were derived 
in~\cite{EnLaWi18}. Therein, the goal functionals were estimated using 
the DWR method (thus an adjoint approach). Hence, the extension to optimal control 
using a gradient-based approach is straightforward. 
The implementation of the following algorithms is done in the open-source
library 
\verb|DOpElib|~\cite{dope,DOpElib}.
For a general overview 
of optimization algorithms, we refer to~\cite{NoWri06,DissMeidner07}.
First, we present the reduced Newton method described
in Algorithm~\ref{inexat_newton_algorithm_for_multiple_goal_functionals}.

\begin{algorithm}[H]
  \caption{Reduced Newton algorithm for multiple-goal
    functionals with adaptive stopping rule on level $l$ } \label{inexat_newton_algorithm_for_multiple_goal_functionals}
  \begin{algorithmic}[1]
    \State Start with  some initial guess $q^{l,0}_h \in Q_h^l$ and $k=0$.
    \State For $p^{l,0}_h$, solve $$ j_{h}''(q^{l,0}_h)(v_h,p^{l,0}_h)=(i_{\mathfrak{E},h}^{(0)})'(q^{l,0}_h)(v_h) \quad \forall v_h \in V_h^l,$$
    with  $(i_{\mathfrak{E},h}^{(0)})'$ constructed with $q^{l,(2)}_h$ and $q^{l,0}_h$.
    \While{$|j_{h}'(q^{l,k}_h)(p^{l,k}_h)|> \gamma \eta_h^{l-1,(2)}$}
    \State For $\delta q^{l,k}_h$, solve $$ j_{h}''(q^{l,k}_h)(\delta q^{l,k}_h,v_h)=-j_{h}'(q^{l,k}_h)(v_h)  \quad \forall v_h \in V_h^l.$$
    \State Update : $u^{l,k+1}_h=q^{l,k}_h+\alpha^k \delta q^{l,k}_h$ for some good choice $\alpha^k \in (0,1]$.
    % \State $u^{l,k,(2)}_h=I_h^{h_2}q^{l,k}_h$.
    \State $k=k+1.$
    \State For $p^{l,k}_h$, solve $$ j_{h}''(q^{l,k}_h)(v_h,p^{l,k}_h)=(i_{\mathfrak{E},h}^{(k)})'(q^{l,k}_h)(v_h) \quad \forall v_h \in U_h^l,$$
    $\quad$ with $(i_{\mathfrak{E},h}^{(k)})'$ constructed with $q^{l,(2)}_h$ and $q^{l,k}_h$.
    \EndWhile
  \end{algorithmic}
\end{algorithm}

\begin{remark}
  The parameter $\gamma$ is chosen as $10^{-2}$ in the numerical experiments.
\end{remark}

\begin{remark}
  In~\cite{DOpElib}, we specifically used 
  \verb| DOpE::ReducedNewtonAlgorithm::ReducedNewtonLineSearch| to
  obtain the line search parameter $\alpha^k$.
\end{remark}
\begin{remark}
  The arising linear problems $
  j_{h}''(q^{l,k}_h)(v_h,p^{l,k}_h)=j_{h}'(q^{l,k}_h)(v_h)$
  and \newline $
  j_{h}''(q^{l,k}_h)(v_h,p^{l,k}_h)=(i_{\mathfrak{E},h}^{(k)})'(q^{l,k}_h)(v_h)$
  were solved by using the
  algorithm  \newline\verb|DOpE::ReducedNewtonAlgorithm::SolveReducedLinearSystem|
  implemented in~\cite{DOpElib}.
\end{remark}

With the help of Algorithm
\ref{inexat_newton_algorithm_for_multiple_goal_functionals}, we can now state
the final Algorithm~\ref{final algorithm} used in this paper.

\begin{algorithm}[H]
  \caption{The final algorithm }\label{final algorithm}
  \begin{algorithmic}[1]
    \State Start with some initial guess 
    $q_{h}^{0,(2)}$,$q_h^{0}$, set $l=1$
    and set $TOL_{dis} > 0$.
    \State Solve  $j_h'(q_h^{l,(2)})=0$ for  $q_h^{l,(2)}$ using Newton Algorithm with the initial guess $q_h^{l-1,(2)}$ on the discrete space $Q_{h}^{l,(2)}$. \label{solve Uh2}
    \State Solve $j_h'(q_h^{l})=0$ and $j_h''(q^{l}_h)(\cdot,p^{l}_h)=(i_{\mathfrak{E},h}^{})'(q^{l}_h)(\cdot)$ using Reduced Adaptive Newton algorithm with the initial guess $q_h^{l-1}$ on the discrete space $Q_{h}^l$ . \label{final algorithm: solveprimal}
    \State Construct the combined functional $i_{\mathfrak{E},h}$.
    \State Solve the adjoint problem $j_h''(q_h^{l-1,(2)})(\cdot,p^{l,(2)}_h)=(i_{\mathfrak{E},h}^{(k)})'(q^{l,(2)}_h)(\cdot)$ on $V_h^{l,(2)}$. \label{final algorithm: solveadjoint}
    \State Recover $v_h^{l,(2)}$ and $y_h^{l,(2)}$ using the first an last row in~\eqref{Problem: Discrete Adjoint Optimality System} for the enriched spaces $U_h^{l,(2)}$ and $Q_h^{l,(2)}$.
    \State Compute  the local error estimator $\eta_{h,K}$  from element and face contributions following Section~\ref{Subsection: The error estimator parts}.
    \State Mark elements with some refinement strategy. \label{final algorithm: refinement}
    \State Refine marked elements: $\mathcal{T}_h^l \mapsto\mathcal{T }_h^{l+1}$ and $l=l+1$.
    \State If $|\eta_h| < TOL_{dis}$ stop, else go to~\ref{solve Uh2}.
  \end{algorithmic}
\end{algorithm}

\begin{remark}
  In Algorithm \ref{final algorithm} in Step 8, we use D\"orfler marking with $\theta=0.5$ as marking strategy~\cite{Dorfler:1998:ASE:301695.301697}.
\end{remark}
\begin{remark}
  The reduced discrete cost functional $j_h$ on the space $Q_h^{l,(2)}$ is constructed by means of the corresponding discrete solution operator on the enriched space.
\end{remark}

\begin{remark}
  To solve the linear systems arising form the forward state equation, we use the sparse direct solver UMFPACK~\cite{UMFPACK}. 
\end{remark}

% ##############################################################################%
% 
\section{Numerical examples}
\label{sec_num_tests}
In the current section, we provide some numerical examples demonstrating the performance of the
theoretical arguments and algorithms developed previously. The implementation 
is done in \verb|DOpElib|~\cite{dope,DOpElib} using the finite elements 
from deal.II~\cite{BangerthHartmannKanschat2007,dealII84}. 
However, large parts of the programming are new . For this reason,
we first present a linear example with a single goal functional, which has
been already studied in the literature. In the second example, we then 
consider the $p$-Laplacian and again the case of a single goal functional. In Example 3,
we study several nonlinear goal functionals that are simultaneously
controlled.
The quality of our results will be measured by effectivity index which is given by $$I_{\text{eff}}:= \frac{\eta_h^{(2)}}{I(\overline{u},\overline{q})-I(\tilde{u}_h,\tilde{q_h})},$$
whereas the primal and adjoint effectivity indices are defined by 
$$I_{\text{effp}}:= \frac{\rho_u(\tilde{\xi}_h,\tilde{\xi}_h^*)(y_h^{(2)}-\tilde{y}_h) + \rho_z(\tilde{\xi}_h,\tilde{\xi}_h^*)(v_h^{(2)}-\tilde{v}_h) +\rho_q(\tilde{\xi}_h,\tilde{\xi}_h^*)(p_h^{(2)}-\tilde{p}_h)}{I(\overline{u},\overline{q})-I(\tilde{u}_h,\tilde{q_h})},$$ and $$I_{\text{effa}}:= \frac{\rho_v(\tilde{\xi}_h,\tilde{\xi}_h^*)(z_h^{(2)}-\tilde{z}_h) + \rho_y(\tilde{\xi}_h,\tilde{\xi}_h^*)(u_h^{(2)}-\tilde{u}_h)	  +\rho_p(\tilde{\xi}_h,\tilde{\xi}_h^*)(q_h^{(2)}-\tilde{q}_h)}{I(\overline{u},\overline{q})-I(\tilde{u}_h,\tilde{q_h})}.$$
Notice that we do not apply the absolute value to the contributions. Hence, we also estimate the sign of the error.

\subsection{Example 1: linear Laplacian, single goal functional}
In this first numerical test, we consider a standard linear example, 
which is implemented, for instance, in 
\verb|DOpElib|\cite{dope,DOpElib}[\verb|OPT/StatPDE/Example1|, Section 6.1.1].
The main purpose is to validate our novel programming code against 
known findings. The domain is $\Omega := (0,1)^2$.
The right-hand side forces of the PDE are $f(x,y):=\big(20\pi^2
\text{sin}(4\pi x)-\alpha^{-1}\text{sin}( \pi x)\big)\text{sin}(2\pi y)$.
The given control is $q^d:=0$, and the desired state is 
$u^d:=\big(5\pi ^2 \text{sin}(\pi x)+ \text{sin}(4\pi x)\big)\text{sin}(2\pi y)$. 
The regularization is chosen as $\alpha =10^{-2}$.

The problem statement is as follows:
Find $(\overline{u},\overline{q}) \in H^1_0(\Omega) \times L^2(\Omega)$ such that it is a minimizer of 
\begin{align*}
  \min_{(u,q) \in H^1_0(\Omega) \times L^2(\Omega)} J(u,q):= \frac{1}{2}\Vert u-u^d \Vert_{L^2(\Omega)}^2 + \frac{\alpha}{2}\Vert q-q^d \Vert_{L^2(\Omega)}^2,
\end{align*}
with the constraints
\begin{align*}
  -\Delta u &=f+q \qquad &\text{ in } \Omega,\\
  u &= 0\qquad &\text{ on } \partial  \Omega.
\end{align*}

The exact minimizer of the problem is known, and given by $\overline{u}(x,y)=\text{sin}(4 \pi x) \text{sin}(2\pi y)$ and $\overline{q}(x,y)= \alpha^{-1}\text{sin}(\pi x) \text{sin}(2\pi y)$.
First of all, we use $I=J$, so the cost functional as quantity of interest. Here, the exact value is given by $J(\overline{u},\overline{q}) = \frac{1}{8}\big(25 \pi^4 + \alpha^{-1}\big)$.\newline

In the Figures~\ref{pic_ex_1_a} and~\ref{pic_ex_1_b}, the effectivity index 
$I_{\text{eff}}$ and the error are both shown against the number of degrees of
freedom (DOFs). For the single error parts, primal and adjoint estimators, 
the effectivity indices show significant differences from the asymptotically 
expected value. Combining both parts, then yields an optimal $I_{\text{eff}} = 1$.
Convergence of adaptive and uniform mesh refinement are shown 
in Figure~\ref{pic_ex_1_b}.

\begin{minipage}[t]{0.45 \textwidth}
  \ifMAKEPICS
  \begin{gnuplot}[terminal=epslatex]
    set output "Figures/Example1btex.tex"
    set key right
    #set key opaque
    set datafile separator "|"
    set logscale x
    # set yrange [-0.5:3]
    set grid ytics lc rgb "#bbbbbb" lw 1 lt 0
    set grid xtics lc rgb "#bbbbbb" lw 1 lt 0
    set xlabel '\text{DOFs}'
    set format '%g'
    plot  '< sqlite3 Data/SinglegoalLaplace/Costfunctional/datacost.db "SELECT DISTINCT DOFs, Ieff from data "' u 1:2 w  lp lw 3 title '  \footnotesize $I_{\text{eff}}$', \
    '< sqlite3 Data/SinglegoalLaplace/Costfunctional/datacost.db "SELECT DISTINCT DOFs, Ieffprimal from data "' u 1:2 w  lp lw 3 title '\footnotesize $I_{\text{effp}}$', \
    '< sqlite3 Data/SinglegoalLaplace/Costfunctional/datacost.db "SELECT DISTINCT DOFs, Ieffadjoint from data "' u 1:2 w  lp lw 3 title ' \footnotesize $I_{\text{effa}}$',\
    1 lw 3																	
  \end{gnuplot}
  \fi
  \scalebox{0.65}{% GNUPLOT: LaTeX picture with Postscript
\begingroup
  \makeatletter
  \providecommand\color[2][]{%
    \GenericError{(gnuplot) \space\space\space\@spaces}{%
      Package color not loaded in conjunction with
      terminal option `colourtext'%
    }{See the gnuplot documentation for explanation.%
    }{Either use 'blacktext' in gnuplot or load the package
      color.sty in LaTeX.}%
    \renewcommand\color[2][]{}%
  }%
  \providecommand\includegraphics[2][]{%
    \GenericError{(gnuplot) \space\space\space\@spaces}{%
      Package graphicx or graphics not loaded%
    }{See the gnuplot documentation for explanation.%
    }{The gnuplot epslatex terminal needs graphicx.sty or graphics.sty.}%
    \renewcommand\includegraphics[2][]{}%
  }%
  \providecommand\rotatebox[2]{#2}%
  \@ifundefined{ifGPcolor}{%
    \newif\ifGPcolor
    \GPcolorfalse
  }{}%
  \@ifundefined{ifGPblacktext}{%
    \newif\ifGPblacktext
    \GPblacktexttrue
  }{}%
  % define a \g@addto@macro without @ in the name:
  \let\gplgaddtomacro\g@addto@macro
  % define empty templates for all commands taking text:
  \gdef\gplbacktext{}%
  \gdef\gplfronttext{}%
  \makeatother
  \ifGPblacktext
    % no textcolor at all
    \def\colorrgb#1{}%
    \def\colorgray#1{}%
  \else
    % gray or color?
    \ifGPcolor
      \def\colorrgb#1{\color[rgb]{#1}}%
      \def\colorgray#1{\color[gray]{#1}}%
      \expandafter\def\csname LTw\endcsname{\color{white}}%
      \expandafter\def\csname LTb\endcsname{\color{black}}%
      \expandafter\def\csname LTa\endcsname{\color{black}}%
      \expandafter\def\csname LT0\endcsname{\color[rgb]{1,0,0}}%
      \expandafter\def\csname LT1\endcsname{\color[rgb]{0,1,0}}%
      \expandafter\def\csname LT2\endcsname{\color[rgb]{0,0,1}}%
      \expandafter\def\csname LT3\endcsname{\color[rgb]{1,0,1}}%
      \expandafter\def\csname LT4\endcsname{\color[rgb]{0,1,1}}%
      \expandafter\def\csname LT5\endcsname{\color[rgb]{1,1,0}}%
      \expandafter\def\csname LT6\endcsname{\color[rgb]{0,0,0}}%
      \expandafter\def\csname LT7\endcsname{\color[rgb]{1,0.3,0}}%
      \expandafter\def\csname LT8\endcsname{\color[rgb]{0.5,0.5,0.5}}%
    \else
      % gray
      \def\colorrgb#1{\color{black}}%
      \def\colorgray#1{\color[gray]{#1}}%
      \expandafter\def\csname LTw\endcsname{\color{white}}%
      \expandafter\def\csname LTb\endcsname{\color{black}}%
      \expandafter\def\csname LTa\endcsname{\color{black}}%
      \expandafter\def\csname LT0\endcsname{\color{black}}%
      \expandafter\def\csname LT1\endcsname{\color{black}}%
      \expandafter\def\csname LT2\endcsname{\color{black}}%
      \expandafter\def\csname LT3\endcsname{\color{black}}%
      \expandafter\def\csname LT4\endcsname{\color{black}}%
      \expandafter\def\csname LT5\endcsname{\color{black}}%
      \expandafter\def\csname LT6\endcsname{\color{black}}%
      \expandafter\def\csname LT7\endcsname{\color{black}}%
      \expandafter\def\csname LT8\endcsname{\color{black}}%
    \fi
  \fi
    \setlength{\unitlength}{0.0500bp}%
    \ifx\gptboxheight\undefined%
      \newlength{\gptboxheight}%
      \newlength{\gptboxwidth}%
      \newsavebox{\gptboxtext}%
    \fi%
    \setlength{\fboxrule}{0.5pt}%
    \setlength{\fboxsep}{1pt}%
\begin{picture}(7200.00,5040.00)%
    \gplgaddtomacro\gplbacktext{%
      \csname LTb\endcsname%
      \put(990,704){\makebox(0,0)[r]{\strut{}1e-05}}%
      \csname LTb\endcsname%
      \put(990,1286){\makebox(0,0)[r]{\strut{}0.0001}}%
      \csname LTb\endcsname%
      \put(990,1867){\makebox(0,0)[r]{\strut{}0.001}}%
      \csname LTb\endcsname%
      \put(990,2449){\makebox(0,0)[r]{\strut{}0.01}}%
      \csname LTb\endcsname%
      \put(990,3030){\makebox(0,0)[r]{\strut{}0.1}}%
      \csname LTb\endcsname%
      \put(990,3612){\makebox(0,0)[r]{\strut{}1}}%
      \csname LTb\endcsname%
      \put(990,4193){\makebox(0,0)[r]{\strut{}10}}%
      \csname LTb\endcsname%
      \put(990,4775){\makebox(0,0)[r]{\strut{}100}}%
      \csname LTb\endcsname%
      \put(1122,484){\makebox(0,0){\strut{}10}}%
      \csname LTb\endcsname%
      \put(2258,484){\makebox(0,0){\strut{}100}}%
      \csname LTb\endcsname%
      \put(3394,484){\makebox(0,0){\strut{}1000}}%
      \csname LTb\endcsname%
      \put(4531,484){\makebox(0,0){\strut{}10000}}%
      \csname LTb\endcsname%
      \put(5667,484){\makebox(0,0){\strut{}100000}}%
      \csname LTb\endcsname%
      \put(6803,484){\makebox(0,0){\strut{}1e+06}}%
    }%
    \gplgaddtomacro\gplfronttext{%
      \csname LTb\endcsname%
      \put(3962,154){\makebox(0,0){\strut{}\text{DOFs}}}%
      \csname LTb\endcsname%
      \put(5816,4602){\makebox(0,0)[r]{\strut{} \footnotesize Error in $I_\mathfrak{E}$ (adp.)}}%
      \csname LTb\endcsname%
      \put(5816,4382){\makebox(0,0)[r]{\strut{} \footnotesize Estimated Error}}%
      \csname LTb\endcsname%
      \put(5816,4162){\makebox(0,0)[r]{\strut{}  \footnotesize Error in $I_\mathfrak{E}$ (uni.)}}%
      \csname LTb\endcsname%
      \put(5816,3942){\makebox(0,0)[r]{\strut{}$\mathcal{O}(\text{DOFs}^{-1})$}}%
    }%
    \gplbacktext
    \put(0,0){\includegraphics{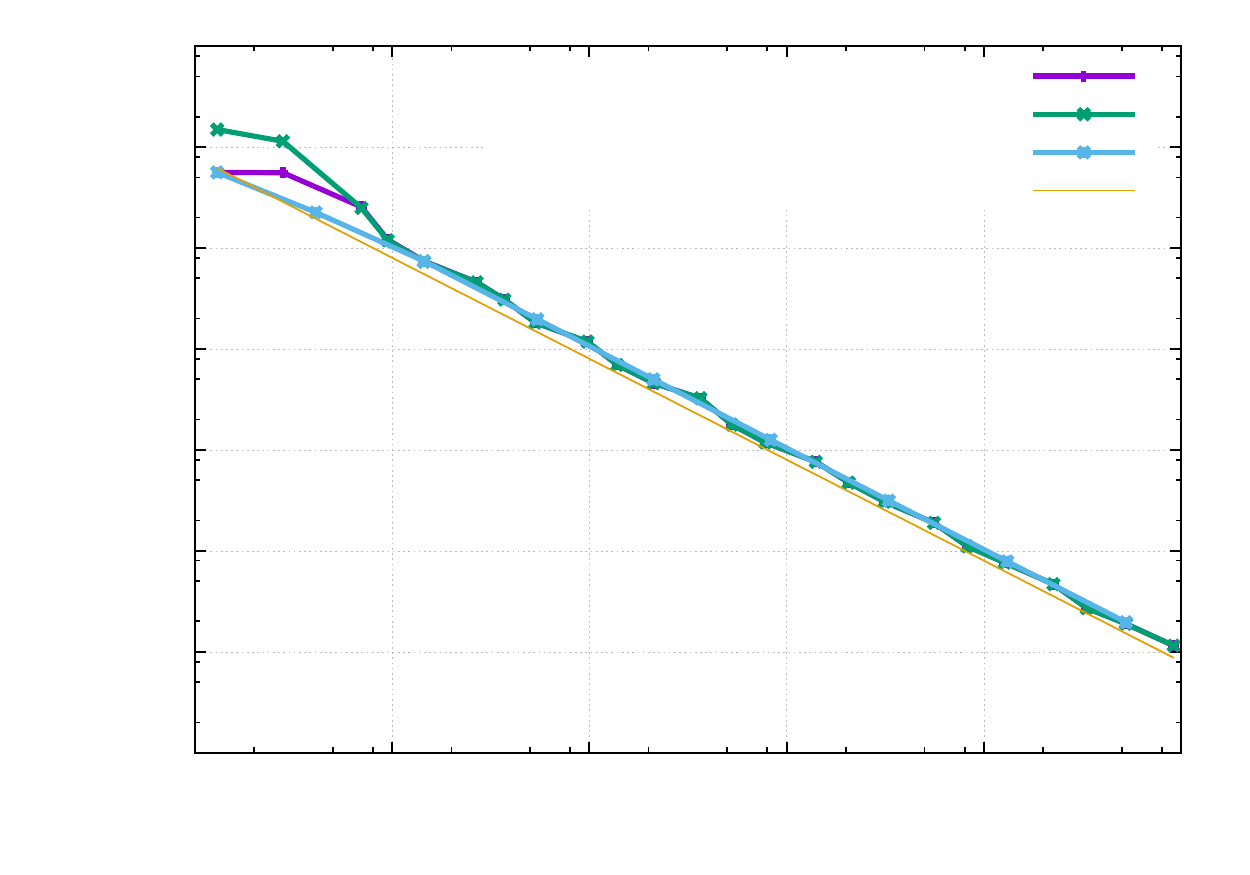}}%
    \gplfronttext
  \end{picture}%
\endgroup
}
  \captionof{figure}{Example 1. $I_{\text{eff}}$  vs DOFs for  the
    linear model problem. \label{pic_ex_1_a}}
  
\end{minipage}\hfill % 
\begin{minipage}[t]{0.45 \textwidth}
  \ifMAKEPICS
  \begin{gnuplot}[terminal=epslatex]
    set output "Figures/Example1atex.tex"
    set key opaque
    set datafile separator "|"
    set logscale x
    set logscale y
    set grid ytics lc rgb "#bbbbbb" lw 1 lt 0
    set grid xtics lc rgb "#bbbbbb" lw 1 lt 0
    set xlabel '\text{DOFs}'
    set format '%g'
    plot  '< sqlite3 Data/SinglegoalLaplace/Costfunctional/datacost.db "SELECT DISTINCT DOFs, abs(exacterror) from data "' u 1:2 w  lp lw 3 title ' \footnotesize Error in $I_\mathfrak{E}$ (adp.)', \
    '< sqlite3 Data/SinglegoalLaplace/Costfunctional/datacost.db "SELECT DISTINCT DOFs, abs(estimatederror) from data "' u 1:2 w  lp lw 3 title ' \footnotesize Estimated Error', \
    '< sqlite3 Data/SinglegoalLaplace/Costfunctional/datacostuniform.db "SELECT DISTINCT DOFs, abs(exacterror) from data "' u 1:2 w  lp lw 3 title '  \footnotesize Error in $I_\mathfrak{E}$ (uni.)', \
    80/x   lw  1	title '$\mathcal{O}(\text{DOFs}^{-1})$'
    #	'< sqlite3 Data/SinglegoalLaplace/Costfunctional/datacost.db "SELECT DISTINCT DOFs, abs(estimatederror0) from data "' u 1:2 w  lp lw 3 title ' \footnotesize Error (a)', \
    '< sqlite3 Data/SinglegoalLaplace/Costfunctional/datacost.db "SELECT DISTINCT DOFs, abs(estimatederror1) from data "' u 1:2 w  lp lw 3 title ' \footnotesize Error (a)', \	
    '< sqlite3 Data/SinglegoalLaplace/Costfunctional/datacost.db "SELECT DISTINCT DOFs, abs(estimatederror2) from data "' u 1:2 w  lp lw 3 title ' \footnotesize Error (a)', \															
    '< sqlite3 Data/SinglegoalLaplace/Costfunctional/datacost.db "SELECT DISTINCT DOFs, correctederror from data "' u 1:2 w  lp lw 3 title ' \footnotesize Error (a)',\	
  \end{gnuplot}
  \fi
  {		\scalebox{0.65}{% GNUPLOT: LaTeX picture with Postscript
\begingroup
  \makeatletter
  \providecommand\color[2][]{%
    \GenericError{(gnuplot) \space\space\space\@spaces}{%
      Package color not loaded in conjunction with
      terminal option `colourtext'%
    }{See the gnuplot documentation for explanation.%
    }{Either use 'blacktext' in gnuplot or load the package
      color.sty in LaTeX.}%
    \renewcommand\color[2][]{}%
  }%
  \providecommand\includegraphics[2][]{%
    \GenericError{(gnuplot) \space\space\space\@spaces}{%
      Package graphicx or graphics not loaded%
    }{See the gnuplot documentation for explanation.%
    }{The gnuplot epslatex terminal needs graphicx.sty or graphics.sty.}%
    \renewcommand\includegraphics[2][]{}%
  }%
  \providecommand\rotatebox[2]{#2}%
  \@ifundefined{ifGPcolor}{%
    \newif\ifGPcolor
    \GPcolorfalse
  }{}%
  \@ifundefined{ifGPblacktext}{%
    \newif\ifGPblacktext
    \GPblacktexttrue
  }{}%
  % define a \g@addto@macro without @ in the name:
  \let\gplgaddtomacro\g@addto@macro
  % define empty templates for all commands taking text:
  \gdef\gplbacktext{}%
  \gdef\gplfronttext{}%
  \makeatother
  \ifGPblacktext
    % no textcolor at all
    \def\colorrgb#1{}%
    \def\colorgray#1{}%
  \else
    % gray or color?
    \ifGPcolor
      \def\colorrgb#1{\color[rgb]{#1}}%
      \def\colorgray#1{\color[gray]{#1}}%
      \expandafter\def\csname LTw\endcsname{\color{white}}%
      \expandafter\def\csname LTb\endcsname{\color{black}}%
      \expandafter\def\csname LTa\endcsname{\color{black}}%
      \expandafter\def\csname LT0\endcsname{\color[rgb]{1,0,0}}%
      \expandafter\def\csname LT1\endcsname{\color[rgb]{0,1,0}}%
      \expandafter\def\csname LT2\endcsname{\color[rgb]{0,0,1}}%
      \expandafter\def\csname LT3\endcsname{\color[rgb]{1,0,1}}%
      \expandafter\def\csname LT4\endcsname{\color[rgb]{0,1,1}}%
      \expandafter\def\csname LT5\endcsname{\color[rgb]{1,1,0}}%
      \expandafter\def\csname LT6\endcsname{\color[rgb]{0,0,0}}%
      \expandafter\def\csname LT7\endcsname{\color[rgb]{1,0.3,0}}%
      \expandafter\def\csname LT8\endcsname{\color[rgb]{0.5,0.5,0.5}}%
    \else
      % gray
      \def\colorrgb#1{\color{black}}%
      \def\colorgray#1{\color[gray]{#1}}%
      \expandafter\def\csname LTw\endcsname{\color{white}}%
      \expandafter\def\csname LTb\endcsname{\color{black}}%
      \expandafter\def\csname LTa\endcsname{\color{black}}%
      \expandafter\def\csname LT0\endcsname{\color{black}}%
      \expandafter\def\csname LT1\endcsname{\color{black}}%
      \expandafter\def\csname LT2\endcsname{\color{black}}%
      \expandafter\def\csname LT3\endcsname{\color{black}}%
      \expandafter\def\csname LT4\endcsname{\color{black}}%
      \expandafter\def\csname LT5\endcsname{\color{black}}%
      \expandafter\def\csname LT6\endcsname{\color{black}}%
      \expandafter\def\csname LT7\endcsname{\color{black}}%
      \expandafter\def\csname LT8\endcsname{\color{black}}%
    \fi
  \fi
    \setlength{\unitlength}{0.0500bp}%
    \ifx\gptboxheight\undefined%
      \newlength{\gptboxheight}%
      \newlength{\gptboxwidth}%
      \newsavebox{\gptboxtext}%
    \fi%
    \setlength{\fboxrule}{0.5pt}%
    \setlength{\fboxsep}{1pt}%
\begin{picture}(7200.00,5040.00)%
    \gplgaddtomacro\gplbacktext{%
      \csname LTb\endcsname%
      \put(990,704){\makebox(0,0)[r]{\strut{}1e-05}}%
      \csname LTb\endcsname%
      \put(990,1286){\makebox(0,0)[r]{\strut{}0.0001}}%
      \csname LTb\endcsname%
      \put(990,1867){\makebox(0,0)[r]{\strut{}0.001}}%
      \csname LTb\endcsname%
      \put(990,2449){\makebox(0,0)[r]{\strut{}0.01}}%
      \csname LTb\endcsname%
      \put(990,3030){\makebox(0,0)[r]{\strut{}0.1}}%
      \csname LTb\endcsname%
      \put(990,3612){\makebox(0,0)[r]{\strut{}1}}%
      \csname LTb\endcsname%
      \put(990,4193){\makebox(0,0)[r]{\strut{}10}}%
      \csname LTb\endcsname%
      \put(990,4775){\makebox(0,0)[r]{\strut{}100}}%
      \csname LTb\endcsname%
      \put(1122,484){\makebox(0,0){\strut{}10}}%
      \csname LTb\endcsname%
      \put(2258,484){\makebox(0,0){\strut{}100}}%
      \csname LTb\endcsname%
      \put(3394,484){\makebox(0,0){\strut{}1000}}%
      \csname LTb\endcsname%
      \put(4531,484){\makebox(0,0){\strut{}10000}}%
      \csname LTb\endcsname%
      \put(5667,484){\makebox(0,0){\strut{}100000}}%
      \csname LTb\endcsname%
      \put(6803,484){\makebox(0,0){\strut{}1e+06}}%
    }%
    \gplgaddtomacro\gplfronttext{%
      \csname LTb\endcsname%
      \put(3962,154){\makebox(0,0){\strut{}\text{DOFs}}}%
      \csname LTb\endcsname%
      \put(5816,4602){\makebox(0,0)[r]{\strut{} \footnotesize Error in $I_\mathfrak{E}$ (adp.)}}%
      \csname LTb\endcsname%
      \put(5816,4382){\makebox(0,0)[r]{\strut{} \footnotesize Estimated Error}}%
      \csname LTb\endcsname%
      \put(5816,4162){\makebox(0,0)[r]{\strut{}  \footnotesize Error in $I_\mathfrak{E}$ (uni.)}}%
      \csname LTb\endcsname%
      \put(5816,3942){\makebox(0,0)[r]{\strut{}$\mathcal{O}(\text{DOFs}^{-1})$}}%
    }%
    \gplbacktext
    \put(0,0){\includegraphics{Figures/Example1a}}%
    \gplfronttext
  \end{picture}%
\endgroup
}
    \captionof{figure}{Example 1. Error vs DOFs for the linear model problem.\label{pic_ex_1_b}}
  }

\end{minipage}
\vspace*{1cm}

In this second part of the example, we apply the method to a quantity that is different
to the cost functional. We are interested in $I(u,q):=\Vert
u\Vert_{L^1(\Omega)}$. The exact value is given by
$I(\overline{u},\overline{q})= 4\pi^{-2}$.
The corresponding numerical findings are displayed in 
the Figures~\ref{pic_ex_1_c} and~\ref{pic_ex_1_d}. We observe excellent 
effectivity indices in Figure~\ref{pic_ex_1_c}. Optimal convergence rates 
also in comparison with uniform mesh refinement are observed in Figure~
\ref{pic_ex_1_d}.

\begin{minipage}[t]{0.45 \textwidth}
  \ifMAKEPICS
  \begin{gnuplot}[terminal=epslatex]
    set output "Figures/ExampleL1atex.tex"
    set key right
    #set key opaque
    set datafile separator "|"
    set logscale x
    # set yrange [-0.5:3]
    set grid ytics lc rgb "#bbbbbb" lw 1 lt 0
    set grid xtics lc rgb "#bbbbbb" lw 1 lt 0
    set xlabel '\text{DOFs}'
    set format '%g'
    plot  '< sqlite3 Data/SinglegoalLaplace/L1Norm/dataL1.db "SELECT DISTINCT DOFs, Ieff from data "' u 1:2 w  lp lw 3 title '  \footnotesize $I_{\text{eff}}$', \
    '< sqlite3 Data/SinglegoalLaplace/L1Norm/dataL1.db "SELECT DISTINCT DOFs, Ieffprimal from data "' u 1:2 w  lp lw 3 title '\footnotesize $I_{\text{effp}}$', \
    '< sqlite3 Data/SinglegoalLaplace/L1Norm/dataL1.db "SELECT DISTINCT DOFs, Ieffadjoint from data "' u 1:2 w  lp lw 3 title ' \footnotesize $I_{\text{effa}}$',\
    1 lw 3																	
  \end{gnuplot}
  \fi
  \scalebox{0.65}{% GNUPLOT: LaTeX picture with Postscript
\begingroup
  \makeatletter
  \providecommand\color[2][]{%
    \GenericError{(gnuplot) \space\space\space\@spaces}{%
      Package color not loaded in conjunction with
      terminal option `colourtext'%
    }{See the gnuplot documentation for explanation.%
    }{Either use 'blacktext' in gnuplot or load the package
      color.sty in LaTeX.}%
    \renewcommand\color[2][]{}%
  }%
  \providecommand\includegraphics[2][]{%
    \GenericError{(gnuplot) \space\space\space\@spaces}{%
      Package graphicx or graphics not loaded%
    }{See the gnuplot documentation for explanation.%
    }{The gnuplot epslatex terminal needs graphicx.sty or graphics.sty.}%
    \renewcommand\includegraphics[2][]{}%
  }%
  \providecommand\rotatebox[2]{#2}%
  \@ifundefined{ifGPcolor}{%
    \newif\ifGPcolor
    \GPcolorfalse
  }{}%
  \@ifundefined{ifGPblacktext}{%
    \newif\ifGPblacktext
    \GPblacktexttrue
  }{}%
  % define a \g@addto@macro without @ in the name:
  \let\gplgaddtomacro\g@addto@macro
  % define empty templates for all commands taking text:
  \gdef\gplbacktext{}%
  \gdef\gplfronttext{}%
  \makeatother
  \ifGPblacktext
    % no textcolor at all
    \def\colorrgb#1{}%
    \def\colorgray#1{}%
  \else
    % gray or color?
    \ifGPcolor
      \def\colorrgb#1{\color[rgb]{#1}}%
      \def\colorgray#1{\color[gray]{#1}}%
      \expandafter\def\csname LTw\endcsname{\color{white}}%
      \expandafter\def\csname LTb\endcsname{\color{black}}%
      \expandafter\def\csname LTa\endcsname{\color{black}}%
      \expandafter\def\csname LT0\endcsname{\color[rgb]{1,0,0}}%
      \expandafter\def\csname LT1\endcsname{\color[rgb]{0,1,0}}%
      \expandafter\def\csname LT2\endcsname{\color[rgb]{0,0,1}}%
      \expandafter\def\csname LT3\endcsname{\color[rgb]{1,0,1}}%
      \expandafter\def\csname LT4\endcsname{\color[rgb]{0,1,1}}%
      \expandafter\def\csname LT5\endcsname{\color[rgb]{1,1,0}}%
      \expandafter\def\csname LT6\endcsname{\color[rgb]{0,0,0}}%
      \expandafter\def\csname LT7\endcsname{\color[rgb]{1,0.3,0}}%
      \expandafter\def\csname LT8\endcsname{\color[rgb]{0.5,0.5,0.5}}%
    \else
      % gray
      \def\colorrgb#1{\color{black}}%
      \def\colorgray#1{\color[gray]{#1}}%
      \expandafter\def\csname LTw\endcsname{\color{white}}%
      \expandafter\def\csname LTb\endcsname{\color{black}}%
      \expandafter\def\csname LTa\endcsname{\color{black}}%
      \expandafter\def\csname LT0\endcsname{\color{black}}%
      \expandafter\def\csname LT1\endcsname{\color{black}}%
      \expandafter\def\csname LT2\endcsname{\color{black}}%
      \expandafter\def\csname LT3\endcsname{\color{black}}%
      \expandafter\def\csname LT4\endcsname{\color{black}}%
      \expandafter\def\csname LT5\endcsname{\color{black}}%
      \expandafter\def\csname LT6\endcsname{\color{black}}%
      \expandafter\def\csname LT7\endcsname{\color{black}}%
      \expandafter\def\csname LT8\endcsname{\color{black}}%
    \fi
  \fi
    \setlength{\unitlength}{0.0500bp}%
    \ifx\gptboxheight\undefined%
      \newlength{\gptboxheight}%
      \newlength{\gptboxwidth}%
      \newsavebox{\gptboxtext}%
    \fi%
    \setlength{\fboxrule}{0.5pt}%
    \setlength{\fboxsep}{1pt}%
\begin{picture}(7200.00,5040.00)%
    \gplgaddtomacro\gplbacktext{%
      \csname LTb\endcsname%
      \put(726,704){\makebox(0,0)[r]{\strut{}-2.5}}%
      \csname LTb\endcsname%
      \put(726,1111){\makebox(0,0)[r]{\strut{}-2}}%
      \csname LTb\endcsname%
      \put(726,1518){\makebox(0,0)[r]{\strut{}-1.5}}%
      \csname LTb\endcsname%
      \put(726,1925){\makebox(0,0)[r]{\strut{}-1}}%
      \csname LTb\endcsname%
      \put(726,2332){\makebox(0,0)[r]{\strut{}-0.5}}%
      \csname LTb\endcsname%
      \put(726,2740){\makebox(0,0)[r]{\strut{}0}}%
      \csname LTb\endcsname%
      \put(726,3147){\makebox(0,0)[r]{\strut{}0.5}}%
      \csname LTb\endcsname%
      \put(726,3554){\makebox(0,0)[r]{\strut{}1}}%
      \csname LTb\endcsname%
      \put(726,3961){\makebox(0,0)[r]{\strut{}1.5}}%
      \csname LTb\endcsname%
      \put(726,4368){\makebox(0,0)[r]{\strut{}2}}%
      \csname LTb\endcsname%
      \put(726,4775){\makebox(0,0)[r]{\strut{}2.5}}%
      \csname LTb\endcsname%
      \put(858,484){\makebox(0,0){\strut{}10}}%
      \csname LTb\endcsname%
      \put(2047,484){\makebox(0,0){\strut{}100}}%
      \csname LTb\endcsname%
      \put(3236,484){\makebox(0,0){\strut{}1000}}%
      \csname LTb\endcsname%
      \put(4425,484){\makebox(0,0){\strut{}10000}}%
      \csname LTb\endcsname%
      \put(5614,484){\makebox(0,0){\strut{}100000}}%
      \csname LTb\endcsname%
      \put(6803,484){\makebox(0,0){\strut{}1e+06}}%
    }%
    \gplgaddtomacro\gplfronttext{%
      \csname LTb\endcsname%
      \put(3830,154){\makebox(0,0){\strut{}\text{DOFs}}}%
      \csname LTb\endcsname%
      \put(5816,4602){\makebox(0,0)[r]{\strut{}  \footnotesize $I_{eff}$}}%
      \csname LTb\endcsname%
      \put(5816,4382){\makebox(0,0)[r]{\strut{}\footnotesize $I_{effp}$}}%
      \csname LTb\endcsname%
      \put(5816,4162){\makebox(0,0)[r]{\strut{} \footnotesize $I_{effa}$}}%
      \csname LTb\endcsname%
      \put(5816,3942){\makebox(0,0)[r]{\strut{}1}}%
    }%
    \gplbacktext
    \put(0,0){\includegraphics{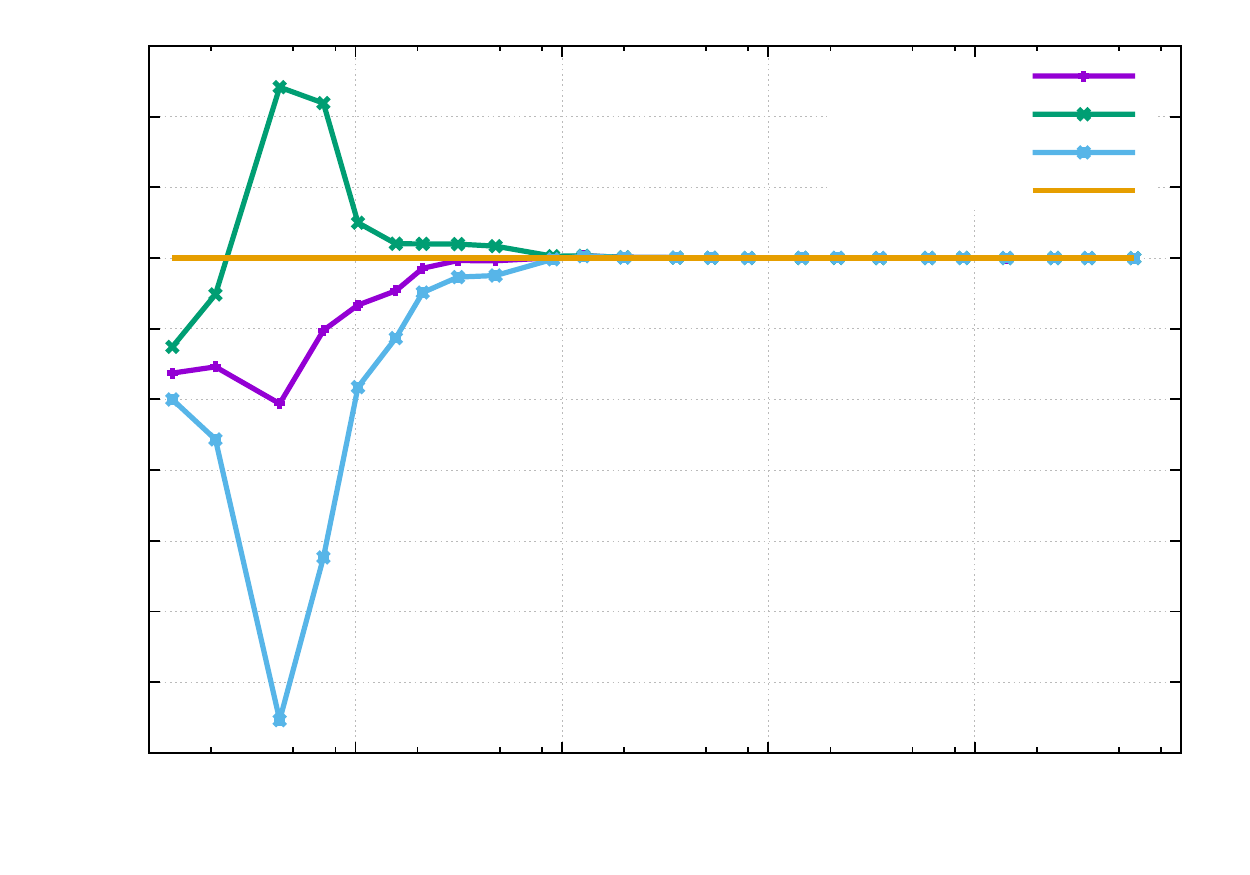}}%
    \gplfronttext
  \end{picture}%
\endgroup
}
  \captionof{figure}{Example 1. $I_{\text{eff}}$  vs DOFs for  the
    linear model problem.  \label{pic_ex_1_c}}
  
\end{minipage} % 
\hfill
\begin{minipage}[t]{0.45 \textwidth}
  \ifMAKEPICS
  \begin{gnuplot}[terminal=epslatex]
    set output "Figures/ExampleL1btex.tex"
    set key opaque
    set datafile separator "|"
    set logscale x
    set logscale y
    set grid ytics lc rgb "#bbbbbb" lw 1 lt 0
    set grid xtics lc rgb "#bbbbbb" lw 1 lt 0
    set xlabel '\text{DOFs}'
    set format '%g'
    plot  '< sqlite3 Data/SinglegoalLaplace/L1Norm/dataL1.db "SELECT DISTINCT DOFs, abs(exacterror) from data "' u 1:2 w  lp lw 3 title ' \footnotesize Error in $I_\mathfrak{E}$ (adp.)', \
    '< sqlite3 Data/SinglegoalLaplace/L1Norm/dataL1.db "SELECT DISTINCT DOFs, abs(estimatederror) from data "' u 1:2 w  lp lw 3 title ' \footnotesize Estimated Error', \
    '< sqlite3 Data/SinglegoalLaplace/L1Norm/dataL1uniform.db "SELECT DISTINCT DOFs, abs(exacterror) from data "' u 1:2 w  lp lw 3 title '  \footnotesize Error in $I_\mathfrak{E}$ (uni.)', \
    30/x   lw  1	title '$\mathcal{O}(\text{DOFs}^{-1})$'
    #	'< sqlite3 Data/SinglegoalLaplace/Costfunctional/datacost.db "SELECT DISTINCT DOFs, abs(estimatederror0) from data "' u 1:2 w  lp lw 3 title ' \footnotesize Error (a)', \
    '< sqlite3 Data/SinglegoalLaplace/Costfunctional/datacost.db "SELECT DISTINCT DOFs, abs(estimatederror1) from data "' u 1:2 w  lp lw 3 title ' \footnotesize Error (a)', \	
    '< sqlite3 Data/SinglegoalLaplace/Costfunctional/datacost.db "SELECT DISTINCT DOFs, abs(estimatederror2) from data "' u 1:2 w  lp lw 3 title ' \footnotesize Error (a)', \															
    '< sqlite3 Data/SinglegoalLaplace/Costfunctional/datacost.db "SELECT DISTINCT DOFs, correctederror from data "' u 1:2 w  lp lw 3 title ' \footnotesize Error (a)',\	
  \end{gnuplot}
  \fi
  {		\scalebox{0.65}{% GNUPLOT: LaTeX picture with Postscript
\begingroup
  \makeatletter
  \providecommand\color[2][]{%
    \GenericError{(gnuplot) \space\space\space\@spaces}{%
      Package color not loaded in conjunction with
      terminal option `colourtext'%
    }{See the gnuplot documentation for explanation.%
    }{Either use 'blacktext' in gnuplot or load the package
      color.sty in LaTeX.}%
    \renewcommand\color[2][]{}%
  }%
  \providecommand\includegraphics[2][]{%
    \GenericError{(gnuplot) \space\space\space\@spaces}{%
      Package graphicx or graphics not loaded%
    }{See the gnuplot documentation for explanation.%
    }{The gnuplot epslatex terminal needs graphicx.sty or graphics.sty.}%
    \renewcommand\includegraphics[2][]{}%
  }%
  \providecommand\rotatebox[2]{#2}%
  \@ifundefined{ifGPcolor}{%
    \newif\ifGPcolor
    \GPcolorfalse
  }{}%
  \@ifundefined{ifGPblacktext}{%
    \newif\ifGPblacktext
    \GPblacktexttrue
  }{}%
  % define a \g@addto@macro without @ in the name:
  \let\gplgaddtomacro\g@addto@macro
  % define empty templates for all commands taking text:
  \gdef\gplbacktext{}%
  \gdef\gplfronttext{}%
  \makeatother
  \ifGPblacktext
    % no textcolor at all
    \def\colorrgb#1{}%
    \def\colorgray#1{}%
  \else
    % gray or color?
    \ifGPcolor
      \def\colorrgb#1{\color[rgb]{#1}}%
      \def\colorgray#1{\color[gray]{#1}}%
      \expandafter\def\csname LTw\endcsname{\color{white}}%
      \expandafter\def\csname LTb\endcsname{\color{black}}%
      \expandafter\def\csname LTa\endcsname{\color{black}}%
      \expandafter\def\csname LT0\endcsname{\color[rgb]{1,0,0}}%
      \expandafter\def\csname LT1\endcsname{\color[rgb]{0,1,0}}%
      \expandafter\def\csname LT2\endcsname{\color[rgb]{0,0,1}}%
      \expandafter\def\csname LT3\endcsname{\color[rgb]{1,0,1}}%
      \expandafter\def\csname LT4\endcsname{\color[rgb]{0,1,1}}%
      \expandafter\def\csname LT5\endcsname{\color[rgb]{1,1,0}}%
      \expandafter\def\csname LT6\endcsname{\color[rgb]{0,0,0}}%
      \expandafter\def\csname LT7\endcsname{\color[rgb]{1,0.3,0}}%
      \expandafter\def\csname LT8\endcsname{\color[rgb]{0.5,0.5,0.5}}%
    \else
      % gray
      \def\colorrgb#1{\color{black}}%
      \def\colorgray#1{\color[gray]{#1}}%
      \expandafter\def\csname LTw\endcsname{\color{white}}%
      \expandafter\def\csname LTb\endcsname{\color{black}}%
      \expandafter\def\csname LTa\endcsname{\color{black}}%
      \expandafter\def\csname LT0\endcsname{\color{black}}%
      \expandafter\def\csname LT1\endcsname{\color{black}}%
      \expandafter\def\csname LT2\endcsname{\color{black}}%
      \expandafter\def\csname LT3\endcsname{\color{black}}%
      \expandafter\def\csname LT4\endcsname{\color{black}}%
      \expandafter\def\csname LT5\endcsname{\color{black}}%
      \expandafter\def\csname LT6\endcsname{\color{black}}%
      \expandafter\def\csname LT7\endcsname{\color{black}}%
      \expandafter\def\csname LT8\endcsname{\color{black}}%
    \fi
  \fi
    \setlength{\unitlength}{0.0500bp}%
    \ifx\gptboxheight\undefined%
      \newlength{\gptboxheight}%
      \newlength{\gptboxwidth}%
      \newsavebox{\gptboxtext}%
    \fi%
    \setlength{\fboxrule}{0.5pt}%
    \setlength{\fboxsep}{1pt}%
\begin{picture}(7200.00,5040.00)%
    \gplgaddtomacro\gplbacktext{%
      \csname LTb\endcsname%
      \put(990,704){\makebox(0,0)[r]{\strut{}1e-05}}%
      \csname LTb\endcsname%
      \put(990,1383){\makebox(0,0)[r]{\strut{}0.0001}}%
      \csname LTb\endcsname%
      \put(990,2061){\makebox(0,0)[r]{\strut{}0.001}}%
      \csname LTb\endcsname%
      \put(990,2740){\makebox(0,0)[r]{\strut{}0.01}}%
      \csname LTb\endcsname%
      \put(990,3418){\makebox(0,0)[r]{\strut{}0.1}}%
      \csname LTb\endcsname%
      \put(990,4097){\makebox(0,0)[r]{\strut{}1}}%
      \csname LTb\endcsname%
      \put(990,4775){\makebox(0,0)[r]{\strut{}10}}%
      \csname LTb\endcsname%
      \put(1122,484){\makebox(0,0){\strut{}10}}%
      \csname LTb\endcsname%
      \put(2258,484){\makebox(0,0){\strut{}100}}%
      \csname LTb\endcsname%
      \put(3394,484){\makebox(0,0){\strut{}1000}}%
      \csname LTb\endcsname%
      \put(4531,484){\makebox(0,0){\strut{}10000}}%
      \csname LTb\endcsname%
      \put(5667,484){\makebox(0,0){\strut{}100000}}%
      \csname LTb\endcsname%
      \put(6803,484){\makebox(0,0){\strut{}1e+06}}%
    }%
    \gplgaddtomacro\gplfronttext{%
      \csname LTb\endcsname%
      \put(3962,154){\makebox(0,0){\strut{}\text{DOFs}}}%
      \csname LTb\endcsname%
      \put(5816,4602){\makebox(0,0)[r]{\strut{} \footnotesize Error in $I_\mathfrak{E}$ (adp.)}}%
      \csname LTb\endcsname%
      \put(5816,4382){\makebox(0,0)[r]{\strut{} \footnotesize Estimated Error}}%
      \csname LTb\endcsname%
      \put(5816,4162){\makebox(0,0)[r]{\strut{}  \footnotesize Error in $I_\mathfrak{E}$ (uni.)}}%
      \csname LTb\endcsname%
      \put(5816,3942){\makebox(0,0)[r]{\strut{}$\mathcal{O}(\text{DOFs}^{-1})$}}%
    }%
    \gplbacktext
    \put(0,0){\includegraphics{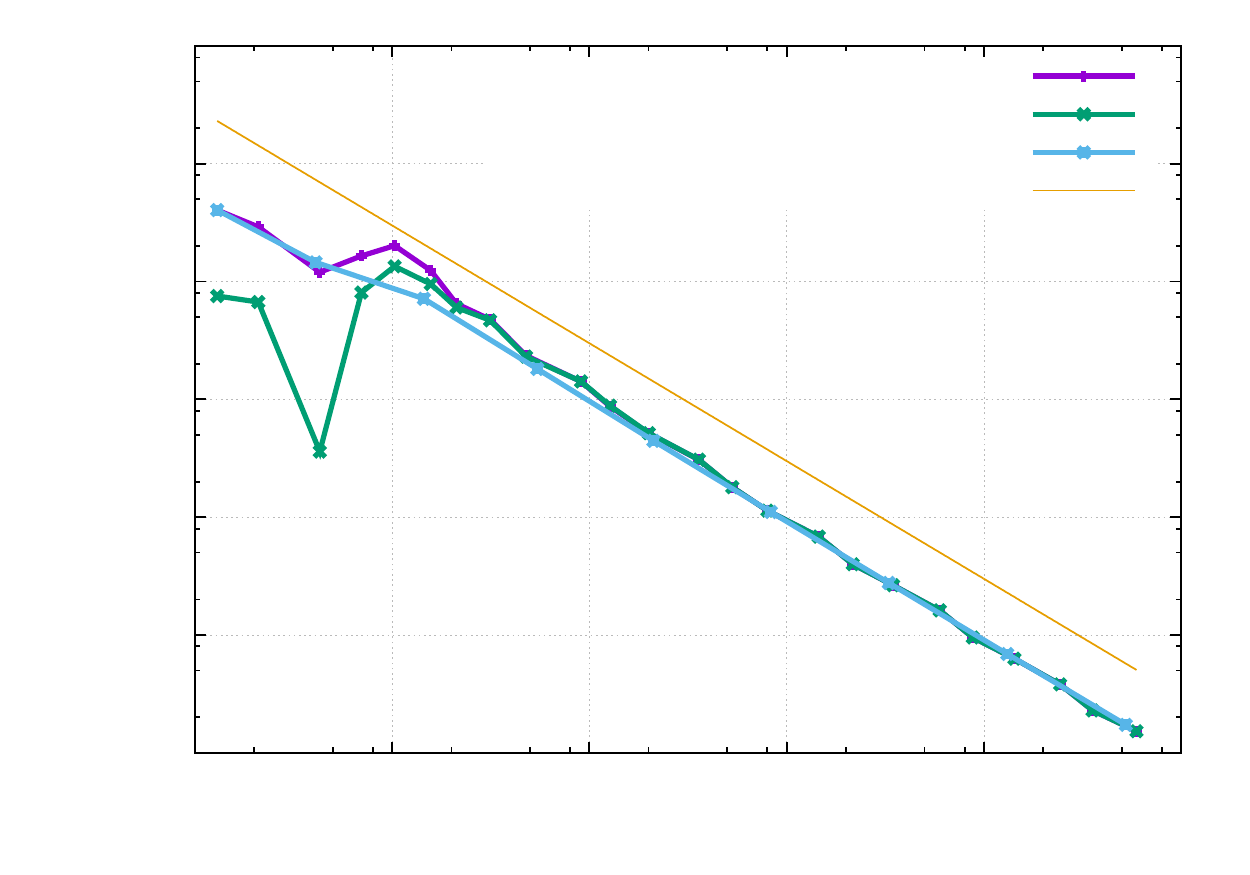}}%
    \gplfronttext
  \end{picture}%
\endgroup
}
    \captionof{figure}{Example 1. Error vs DOFs for the linear model problem.\label{pic_ex_1_d}}}

\end{minipage}

%%%%%%%%%%%%%%%%%%%%%%%%%%%%%%%%%%%%%%%%%%%%%%%%%%%%%%%%%%%%%%%%%%%%%%%% 
\subsection{Example 2: $p$-Laplacian, single goal functional} 
\begin{figure}[H]
  \centering
  \definecolor{zzttqq}{rgb}{0.6,0.2,0}
  \definecolor{ududff}{rgb}{0.30196078431372547,0.30196078431372547,1}
  \definecolor{xdxdff}{rgb}{0.49019607843137253,0.49019607843137253,1}
  \definecolor{cqcqcq}{rgb}{0.7529411764705882,0.7529411764705882,0.7529411764705882}
  \begin{tikzpicture}[line cap=round,line join=round,>=triangle 45,x=1cm,y=1cm]
    \draw [color=cqcqcq,, xstep=0.5cm,ystep=0.5cm] (-0.0,-0.0) grid (7.0,5.0);
    \clip(0,0) rectangle (7,5);
    \fill[line width=2pt,color=zzttqq,fill=zzttqq,fill opacity=0.10000000149011612] (0,0) -- (0,5) -- (7,5) -- (7,0) -- cycle;
    \fill[line width=2pt,color=blue,fill=blue,fill opacity=0.40000000149011612] (2.5,4.5) -- (4.5,4.5) -- (4.5,2.5) -- (2.5,2.5) -- cycle;
    \fill[line width=2pt,color=white,fill=white,fill opacity=1.00000000149011612] (1,4) -- (1,3) -- (2,3) -- (2,4) -- cycle;
    \fill[line width=2pt,color=white,fill=white,fill opacity=1.0000000049011612] (1,2) -- (1,1) -- (2,1) -- (2,2) -- cycle;
    \fill[line width=2pt,color=white,fill=white,fill opacity=1.0000000000149011612] (3,4) -- (3,3) -- (4,3) -- (4,4) -- cycle;
    \fill[line width=2pt,color=white,fill=white,fill opacity=1.0000000000000149011612] (3,2) -- (3,1) -- (4,1) -- (4,2) -- cycle;
    \fill[line width=2pt,color=white,fill=white,fill opacity=1.00000000000149011612] (5,4) -- (5,3) -- (6,3) -- (6,4) -- cycle;
    \fill[line width=2pt,color=white,fill=white,fill opacity=1.00000000000149011612] (5,2) -- (5,1) -- (6,1) -- (6,2) -- cycle;
    
    \draw [line width=2pt,color=black] (0,0)-- (0,5);
    \draw [line width=2pt,color=black] (0,5)-- (7,5);
    \draw [line width=2pt,color=black] (7,5)-- (7,0);
    \draw [line width=2pt,color=black] (7,0)-- (0,0);
    \draw [line width=2pt,color=black] (1,4)-- (1,3);
    \draw [line width=2pt,color=black] (1,3)-- (2,3);
    \draw [line width=2pt,color=black] (2,3)-- (2,4);
    \draw [line width=2pt,color=black] (2,4)-- (1,4);
    \draw [line width=2pt,color=black] (1,2)-- (1,1);
    \draw [line width=2pt,color=black] (1,1)-- (2,1);
    \draw [line width=2pt,color=black] (2,1)-- (2,2);
    \draw [line width=2pt,color=black] (2,2)-- (1,2);
    \draw [line width=2pt,color=black] (3,4)-- (3,3);
    \draw [line width=2pt,color=black] (3,3)-- (4,3);
    \draw [line width=2pt,color=black] (4,3)-- (4,4);
    \draw [line width=2pt,color=black] (4,4)-- (3,4);
    \draw [line width=2pt,color=black] (3,2)-- (3,1);
    \draw [line width=2pt,color=black] (3,1)-- (4,1);
    \draw [line width=2pt,color=black] (4,1)-- (4,2);
    \draw [line width=2pt,color=black] (4,2)-- (3,2);
    \draw [line width=2pt,color=black] (5,4)-- (5,3);
    \draw [line width=2pt,color=black] (5,3)-- (6,3);
    \draw [line width=2pt,color=black] (6,3)-- (6,4);
    \draw [line width=2pt,color=black] (6,4)-- (5,4);
    \draw [line width=2pt,color=black] (5,2)-- (5,1);
    \draw [line width=2pt,color=black] (5,1)-- (6,1);
    \draw [line width=2pt,color=black] (6,1)-- (6,2);
    \draw [line width=2pt,color=black] (6,2)-- (5,2);
    \draw [line width=1pt,color=black] (2.5,4.5)-- (4.5,4.5);
    \draw [line width=1pt,color=black] (4.5,4.5)-- (4.5,2.5);
    \draw [line width=1pt,color=black] (4.5,2.5)-- (2.5,2.5);
    \draw [line width=1pt,color=black] (2.5,2.5)-- (2.5,4.5);
    
    \begin{scriptsize}
      \draw[color=white] (3.5,2.75) node {$u^d=-1$};
      \draw[color=black] (3.5,2.25) node {$u^d=0$};
      \draw[color=black] (0.35,0.25) node {$(0,0)$};
      \draw[color=black] (6.65,4.75) node {$(7,5)$};
    \end{scriptsize}
  \end{tikzpicture}
  \caption{Example 2: The domain $\Omega$ with initial mesh and values of $u^d$.\label{pic_ex_2_a}}

\end{figure}
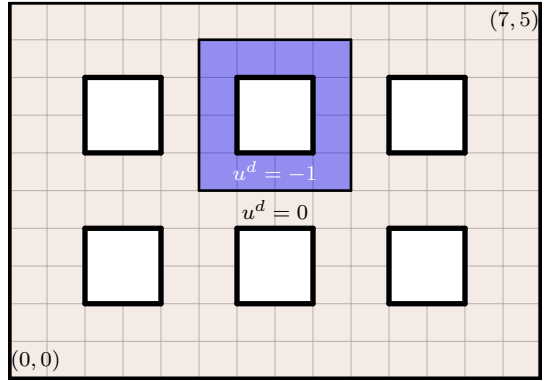

We now proceed to nonlinear state equations  
and consider the example PDE provided in Section~\ref{Example: p-Laplace}. Here,
$\Omega$ (and the initial mesh) and $u^d$ are given in
Figure~\ref{pic_ex_2_a}. Furthermore, $q^d=1$, $p=4$,
$\varepsilon=1$ and $f=0$. In particular, we investigate various
regularization 
parameters $\alpha$. The goal functional $I(u,q)$ is given by 
$I(u,q):= \int_{\Omega} u(x)^2 q(x)^2 dx$.

% \begin{minipage}[t]{0,45\textwidth}
% \end{minipage} \hfill%

\vspace*{1cm}

\begin{table}[H]
  \caption{\label{Table:Differentalphas}Example 1: $I_{\text{effs}}$ for $I(u,q):= \int_{\Omega} u(x)^2 q(x)^2 dx$.}
  \begin{tabular}{|r|r|r|r|r|r|r|r|r|}
    \hline
    \multicolumn{1}{|c|}{$\alpha$} & \multicolumn{2}{c|}{0.01}                             & \multicolumn{2}{c|}{0.1}                              & \multicolumn{2}{c|}{1}                                & \multicolumn{2}{c|}{10}                               \\ \hline
    \multicolumn{1}{|c|}{$l$} & \multicolumn{1}{c|}{$I_{\text{eff}}$} & \multicolumn{1}{c|}{DOFs} & \multicolumn{1}{c|}{$I_{\text{eff}}$} & \multicolumn{1}{c|}{DOFs} & \multicolumn{1}{c|}{$I_{\text{eff}}$} & \multicolumn{1}{c|}{DOFs} & \multicolumn{1}{c|}{$I_{\text{eff}}$} & \multicolumn{1}{c|}{DOFs} \\ \hline
    0     & 0.88 & 275    & 0.69      & 275             & 0.89   & 275    & 0.90       & 275       \\ \hline
    1     & 0.94 & 326    & 0.79      & 506           & 0.91    & 565     & 0.93      & 568       \\ \hline
    2     & 0.94 & 381    & 0.68      & 759           & 0.91  & 832       & 0.92      & 845       \\ \hline
    3     & 0.99 & 561    & 0.66      & 1 266          & 0.92  & 1 367      & 0.93      & 1 451      \\ \hline
    4     & 1.05 & 719    & 0.63      & 2 084          & 0.91    & 2 246    & 0.93      & 2 385      \\ \hline
    5     & 1.06 & 1 151   & 0.50       & 3 013          & 0.89  & 3 115      & 0.93      & 3 263      \\ \hline
    6     & 1.13 & 1 856   & 0.59      & 5 031          & 0.92   & 5 072     & 0.95      & 5 444      \\ \hline
    7     & 1.05 & 2 419   & 0.55      & 8 137          & 0.94   & 8 367     & 0.97      & 8 865      \\ \hline
    8     & 1.11 & 3 363   & 0.36      & 12 498         & 0.94   & 11 880    & 0.97      & 12 479     \\ \hline
    9     & 1.12 & 5 691   & 0.56      & 20 690          & 0.95  & 17 591    & 0.98      & 19 357     \\ \hline
    10    & 1.15 & 7 852   & 0.47      & 33 247         & 0.95  & 31 035    & 0.99      & 32 970     \\ \hline
    11    & 1.13 & 10 752  & 0.38      & 50 864         & 0.95    & 45 721   & 0.99      & 47 850     \\ \hline
    12    & 1.14 & 17 094  & 0.56      & 84 368          & 0.96    & 72 636  & 0.99      & 78 502     \\ \hline
    13    & 1.19 & 25 916  & 0.44      & 135 166        & 0.96   & 126 711   & 0.99      & 133 541    \\ \hline
    14    & 1.14 & 35 482  & 0.39      & 207 466      & 0.96     & 184 754   & 1.00      & 192 946    \\ \hline
    % 15    & 1.11 & 51 194  & $--$ & $--$ & $--$ & $--$ & $--$ & $--$ \\ \hline
    % 16    & 1.14 & 84 411  & $--$ & $--$ & $--$ & $--$ & $--$ & $--$ \\ \hline
    % 17    & 1.30 & 120 485 & $--$ & $--$ & $--$ & $--$ & $--$ & $--$ \\ \hline
    \multicolumn{1}{|c|}{} & \multicolumn{1}{c|}{$I(u,q)$} & \multicolumn{1}{c|}{DOFs} & \multicolumn{1}{c|}{$I(u,q)$} & \multicolumn{1}{c|}{DOFs} & \multicolumn{1}{c|}{$I(u,q)$} & \multicolumn{1}{c|}{DOFs} & \multicolumn{1}{c|}{$I(u,q)$} &
    \multicolumn{1}{c|}{DOFs} \\ \hline
    $\infty$
    & 0.2316036                      &          1 326 503         & 0.07069658                     & 2 127 499 
    & 0.1502366                      & 1 996 755                   & 0.1635741                     & 2 107 007                    \\ \hline
  \end{tabular}
\end{table}
In Table~\ref{Table:Differentalphas}, we obtain, for $\alpha = 0.01,
\ldots, 10$,  effectivity indices in the range of  $0.88$ to $1.30$, which are
excellent findings in view of the nonlinear behavior of the state equation and 
the geometric singularities introduced by the domain.
In the case of $\alpha = 0.1$, we obtain a $I_{\text{eff}}$ in the range of $0.36$ to $0.79$, which might be affected by cancellation effects from adding the different  contributions to the error estimator. The exact value of the functionals was approximated by one additional $p$ and $h$ refinement, and is given in the last line of Table~\ref{Table:Differentalphas} corresponding to $l = \infty$, with additional information on the number of DOFs used to compute this values.

\begin{figure}[H]
  \centering
  \includegraphics[scale
  =0.18]{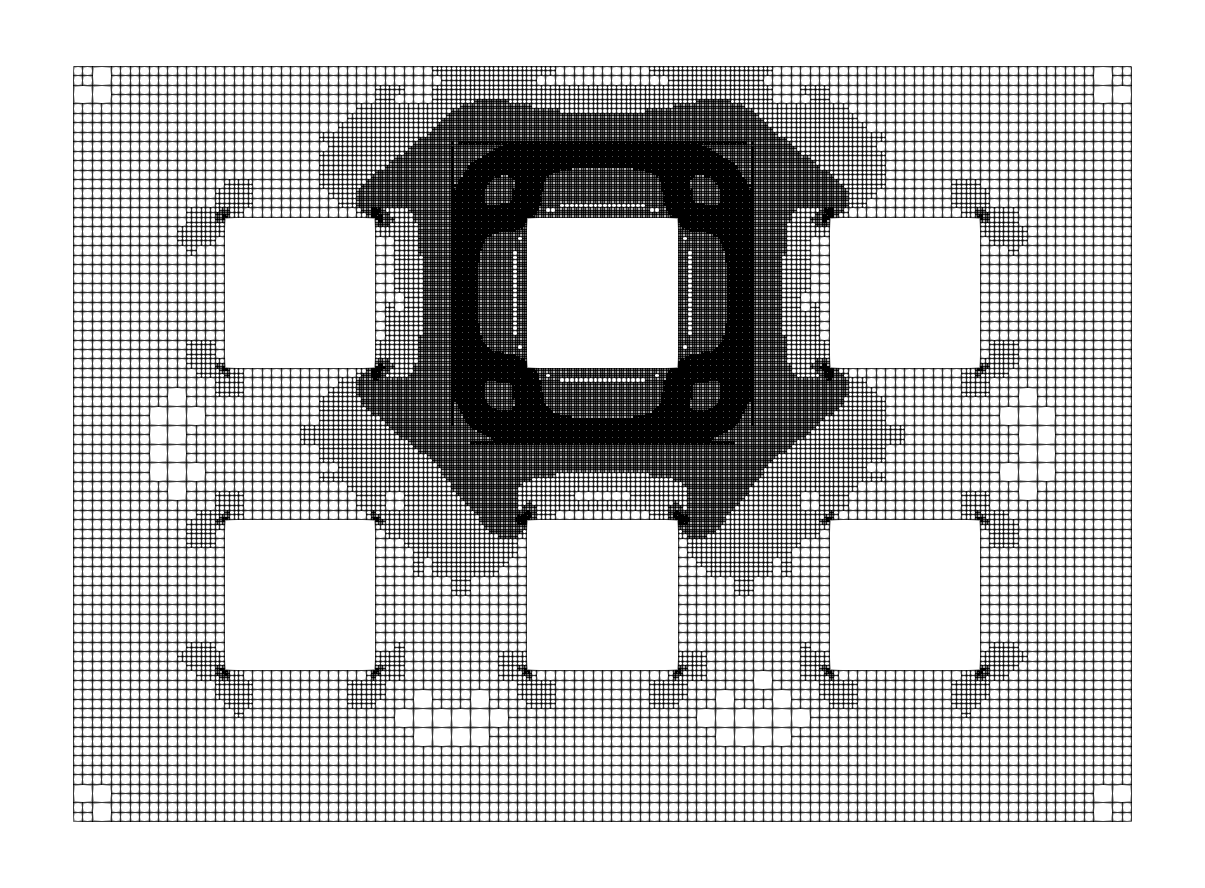}
  \includegraphics[scale =0.18]{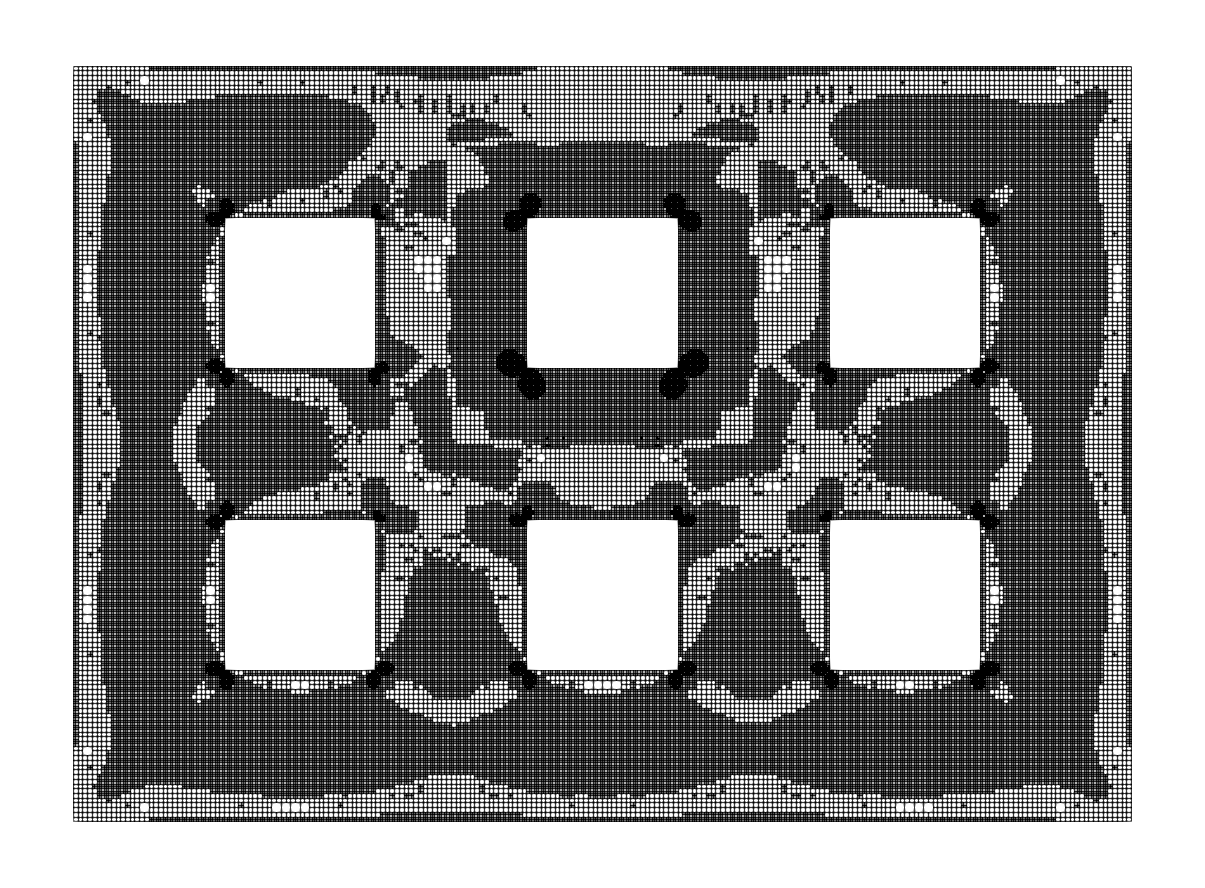}
  \caption{Example 2: Final meshes for $\alpha=10^{-2}$ and $\alpha=10^{-1}$.}
  \label{pic_ex_2_mesh_a}
\end{figure}

\begin{figure}[H]
  \centering
  \includegraphics[scale
  =0.18]{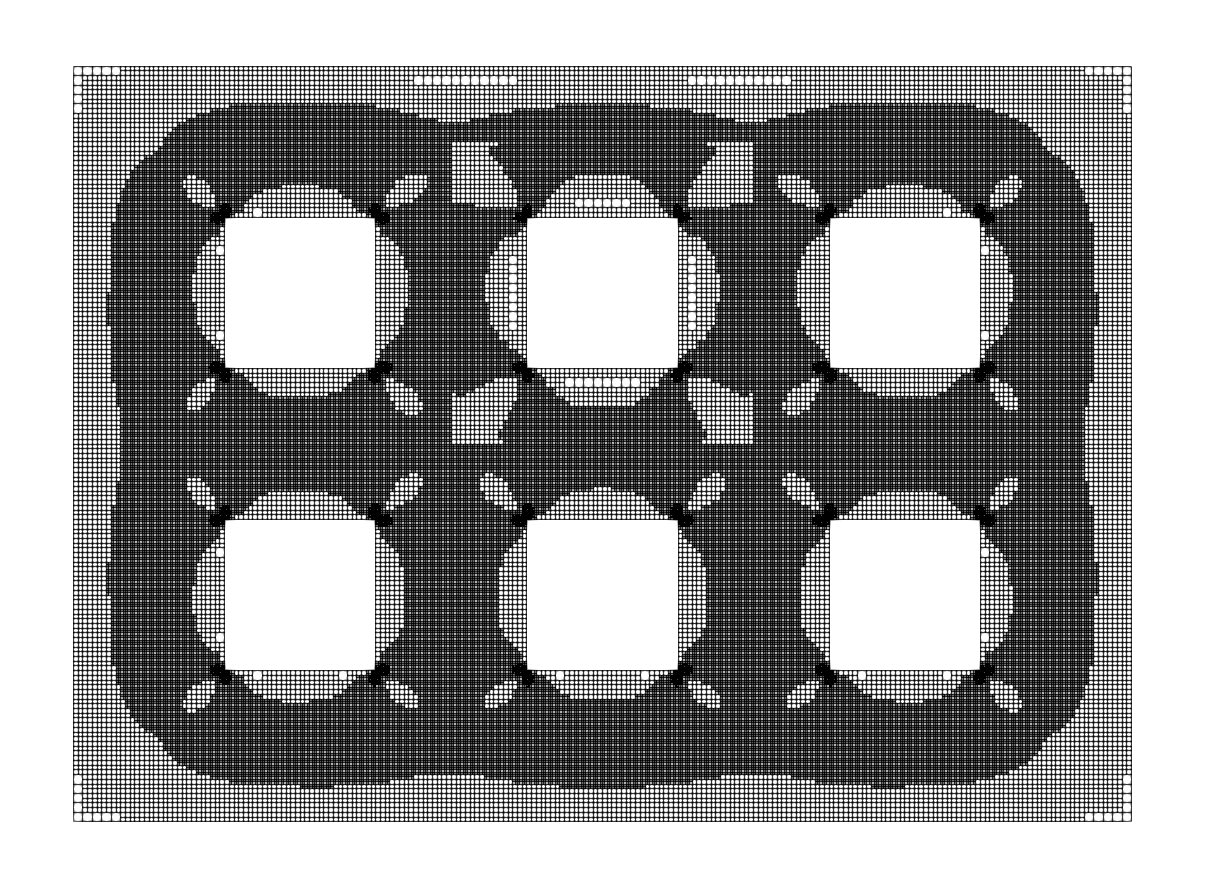}
  \includegraphics[scale =0.18]{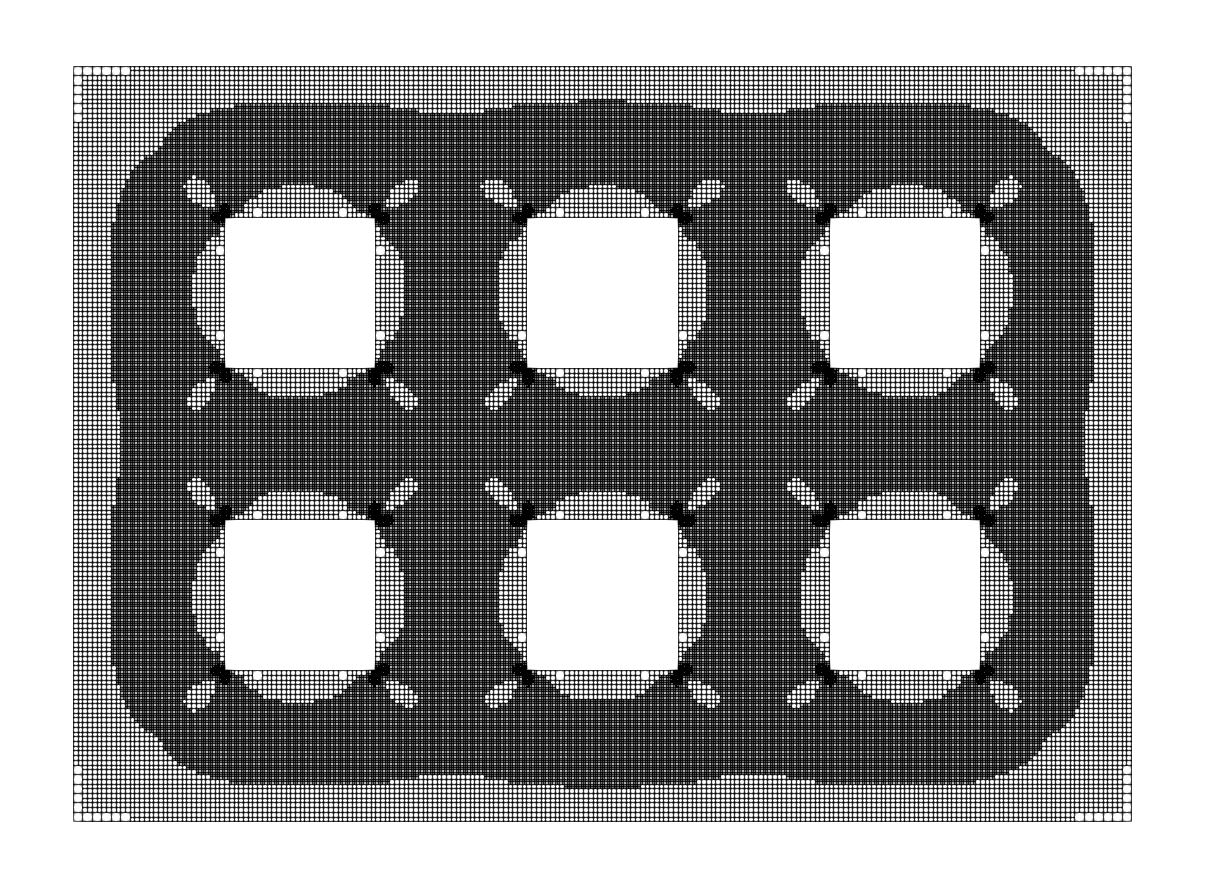}
  \caption{Example 2: Final meshes for $\alpha=10^{-0}$ and $\alpha=10^{1}$.}
  \label{pic_ex_2_mesh_b}
\end{figure}

\begin{figure}[H]
  \centering
  \includegraphics[scale
  =0.18]{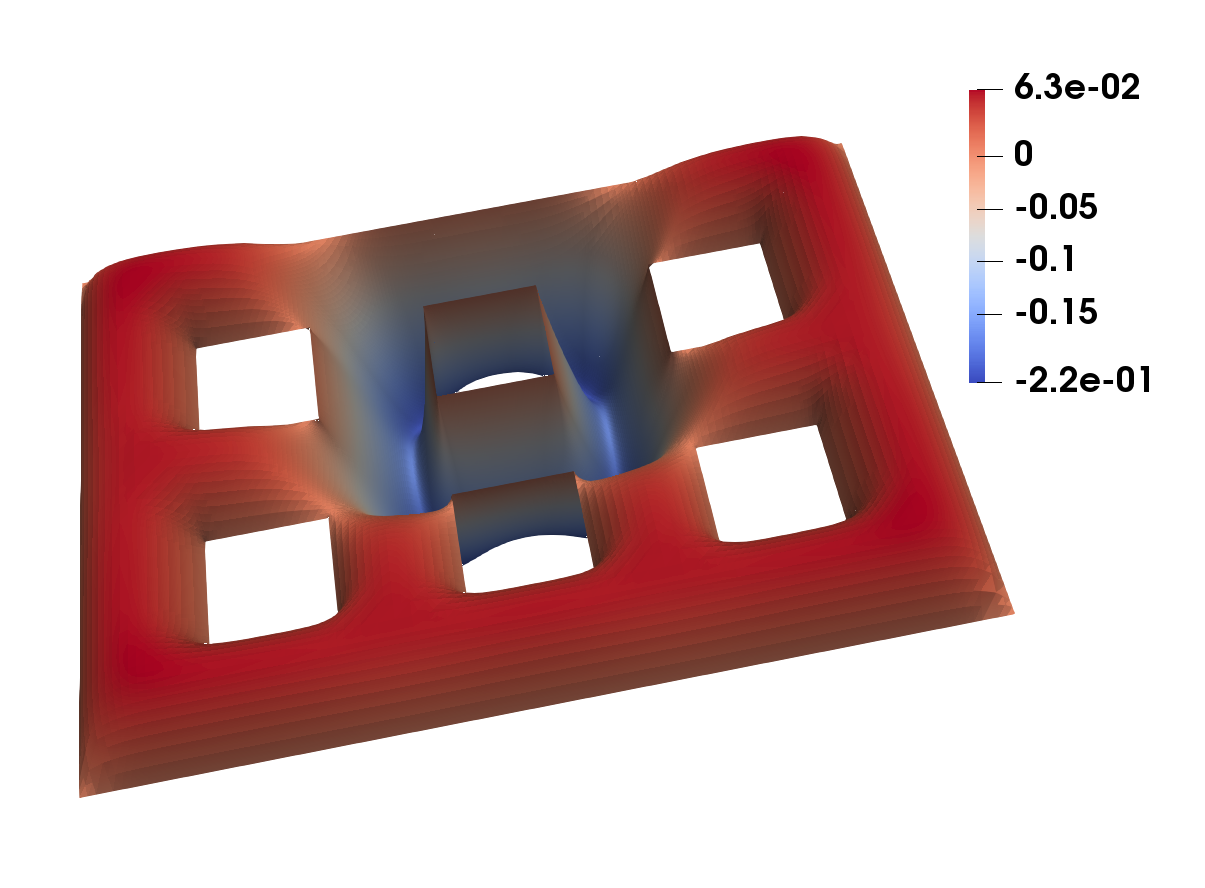}
  \includegraphics[scale =0.18]{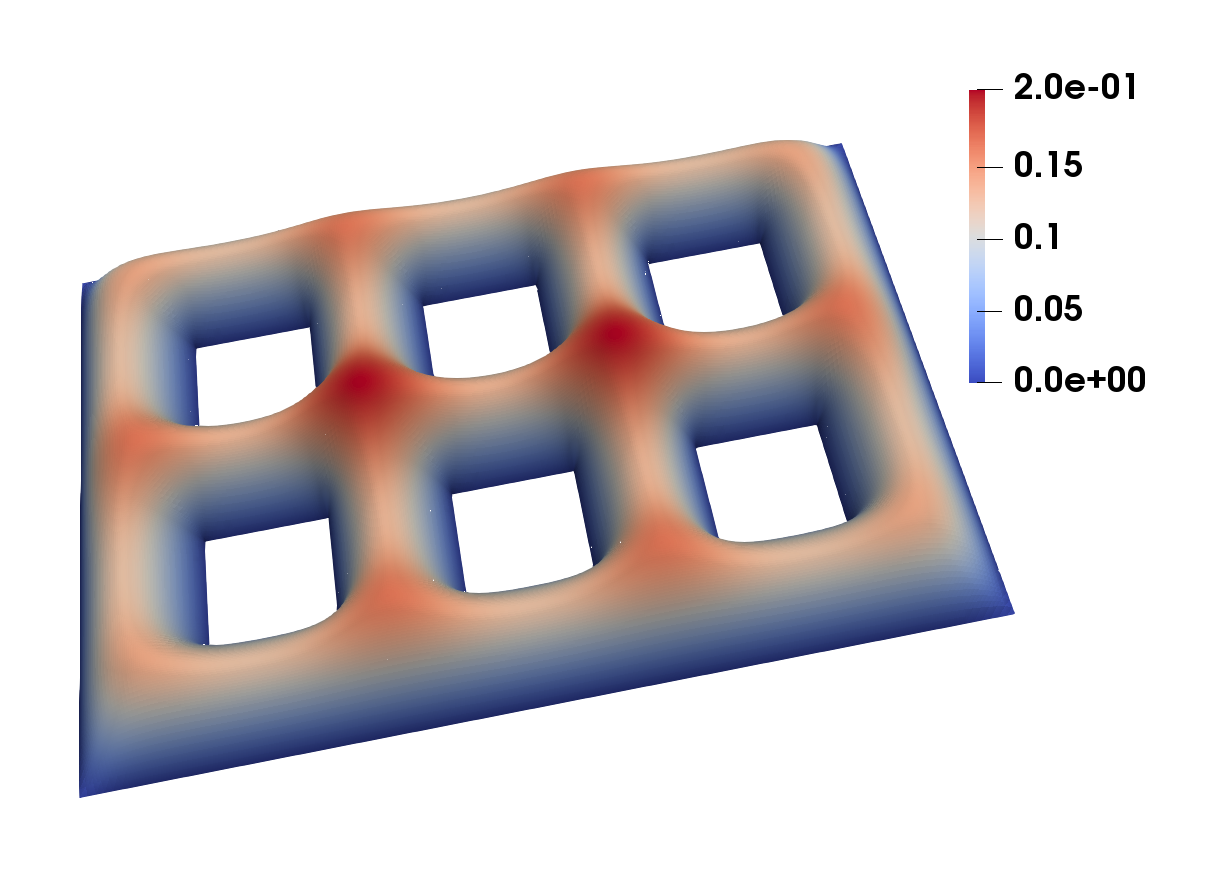}
  \caption{Example 2: State on finest grid for $\alpha=10^{-2}$ and $\alpha=10^{1}$.}
  \label{pic_ex_2_mesh_c}
\end{figure}

\begin{figure}[H]
  \centering
  \includegraphics[scale
  =0.18]{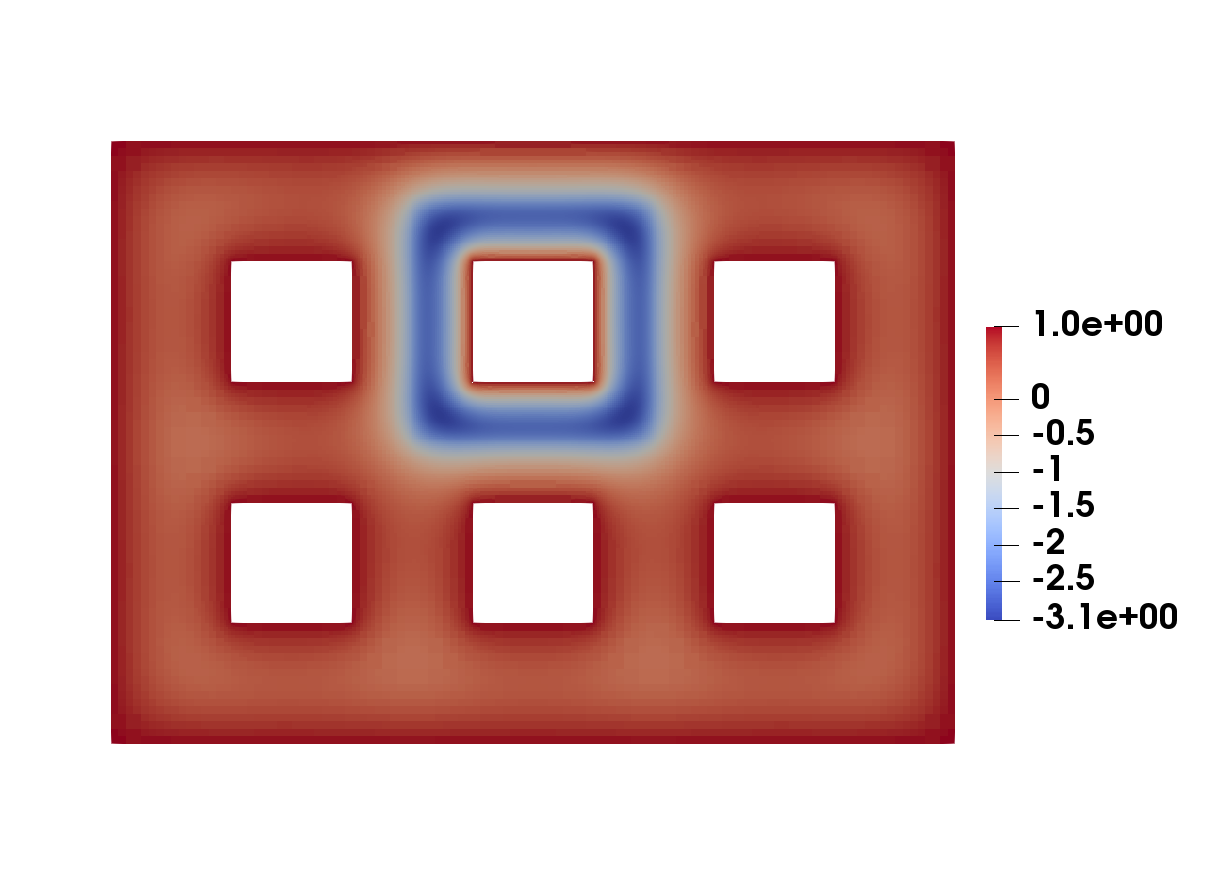}
  \includegraphics[scale =0.18]{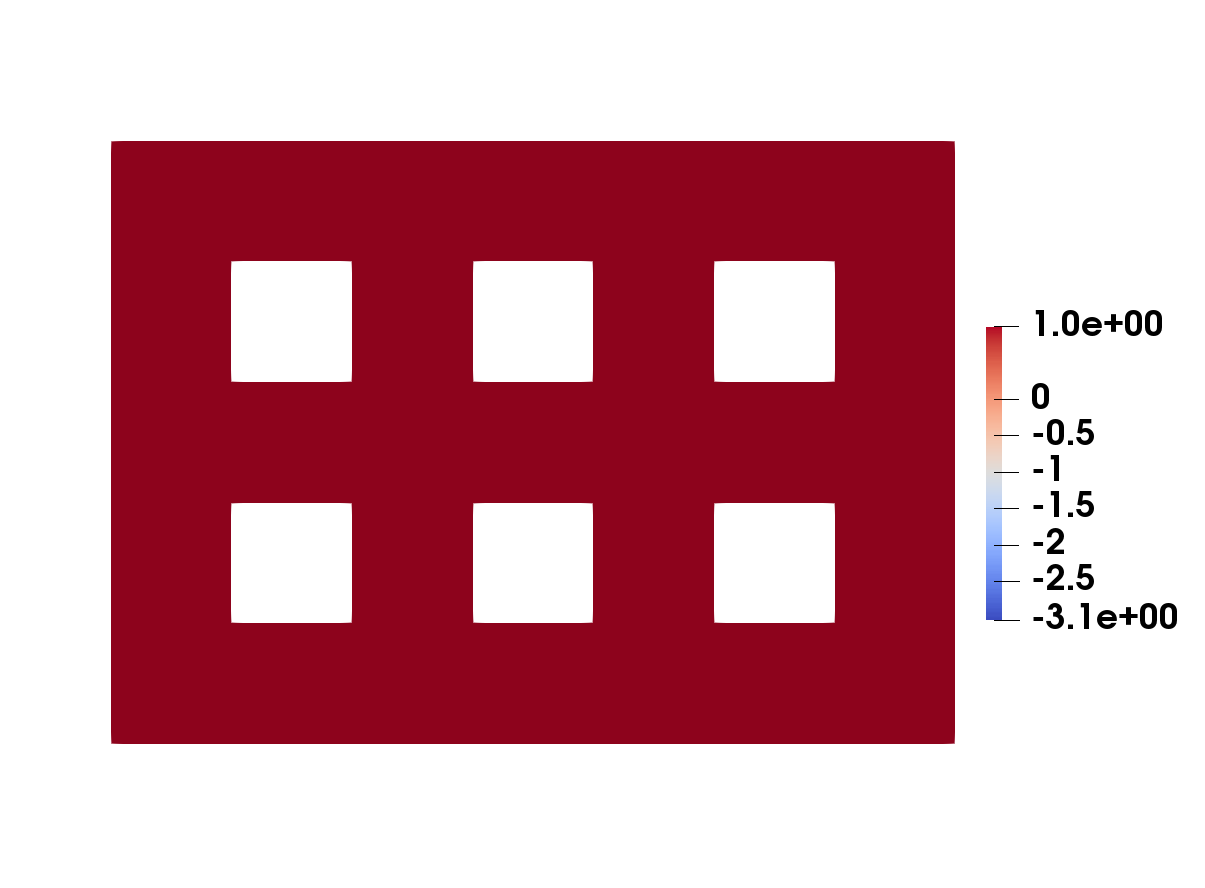}
  \caption{Example 2: Control on finest grid for $\alpha=10^{-2}$ and $\alpha=10^{1}$.}
  \label{pic_ex_2_mesh_d}
\end{figure}

In the Figures~\ref{pic_ex_2_mesh_a} and~\ref{pic_ex_2_mesh_b}, the final
meshes for different $\alpha$ are shown. For $\alpha = 10^{-2}$, we observe 
very localized mesh refinement, while, for larger $\alpha$, the mesh is still
locally refined, but in a somewhat uniform behavior.
The states and controls 
on these final meshes are displayed in the Figures~\ref{pic_ex_2_mesh_c}
and~\ref{pic_ex_2_mesh_d}. 

%\newpage
%%%%%%%%%%%%%%%%%%%%%%%%%%%%%%%%%%%%%%%%%%%%%%%%%%%%%%%%%%%%%%%%%%%%%%%%%%%%%%%%%%%%%% 
\subsection{Example 3: $p$-Laplacian, multiple goal functionals}
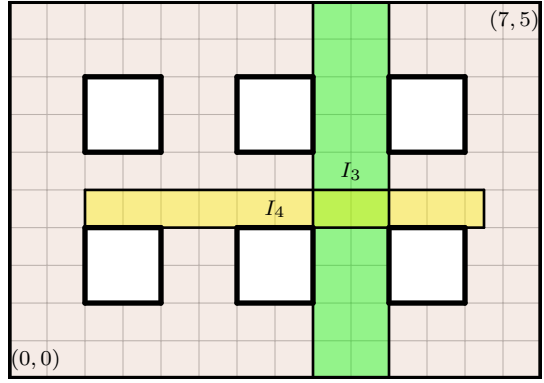
\begin{figure}[H]
  \centering
  \definecolor{zzttqq}{rgb}{0.6,0.2,0}
  \definecolor{ududff}{rgb}{0.30196078431372547,0.30196078431372547,1}
  \definecolor{xdxdff}{rgb}{0.49019607843137253,0.49019607843137253,1}
  \definecolor{cqcqcq}{rgb}{0.7529411764705882,0.7529411764705882,0.7529411764705882}
  \begin{tikzpicture}[line cap=round,line join=round,>=triangle 45,x=1cm,y=1cm]
    \draw [color=cqcqcq,, xstep=0.5cm,ystep=0.5cm] (-0.0,-0.0) grid (7.0,5.0);
    \clip(0,0) rectangle (7,5);
    \fill[line width=2pt,color=zzttqq,fill=zzttqq,fill opacity=0.10000000149011612] (0,0) -- (0,5) -- (7,5) -- (7,0) -- cycle;
    \fill[line width=2pt,color=green,fill=green,fill opacity=0.40000000149011612] (4,0) -- (5,0) -- (5,5) -- (4,5) -- cycle;
    \fill[line width=2pt,color=yellow,fill=yellow,fill opacity=0.40000000149011612] (6.25,2) -- (6.25,2.5) -- (1,2.5) -- (1,2) -- cycle;
    \fill[line width=2pt,color=white,fill=white,fill opacity=1.00000000149011612] (1,4) -- (1,3) -- (2,3) -- (2,4) -- cycle;
    \fill[line width=2pt,color=white,fill=white,fill opacity=1.0000000049011612] (1,2) -- (1,1) -- (2,1) -- (2,2) -- cycle;
    \fill[line width=2pt,color=white,fill=white,fill opacity=1.0000000000149011612] (3,4) -- (3,3) -- (4,3) -- (4,4) -- cycle;
    \fill[line width=2pt,color=white,fill=white,fill opacity=1.0000000000000149011612] (3,2) -- (3,1) -- (4,1) -- (4,2) -- cycle;
    \fill[line width=2pt,color=white,fill=white,fill opacity=1.00000000000149011612] (5,4) -- (5,3) -- (6,3) -- (6,4) -- cycle;
    \fill[line width=2pt,color=white,fill=white,fill opacity=1.00000000000149011612] (5,2) -- (5,1) -- (6,1) -- (6,2) -- cycle;
    
    \draw [line width=2pt,color=black] (0,0)-- (0,5);
    \draw [line width=2pt,color=black] (0,5)-- (7,5);
    \draw [line width=2pt,color=black] (7,5)-- (7,0);
    \draw [line width=2pt,color=black] (7,0)-- (0,0);
    \draw [line width=2pt,color=black] (1,4)-- (1,3);
    \draw [line width=2pt,color=black] (1,3)-- (2,3);
    \draw [line width=2pt,color=black] (2,3)-- (2,4);
    \draw [line width=2pt,color=black] (2,4)-- (1,4);
    \draw [line width=2pt,color=black] (1,2)-- (1,1);
    \draw [line width=2pt,color=black] (1,1)-- (2,1);
    \draw [line width=2pt,color=black] (2,1)-- (2,2);
    \draw [line width=2pt,color=black] (2,2)-- (1,2);
    \draw [line width=2pt,color=black] (3,4)-- (3,3);
    \draw [line width=2pt,color=black] (3,3)-- (4,3);
    \draw [line width=2pt,color=black] (4,3)-- (4,4);
    \draw [line width=2pt,color=black] (4,4)-- (3,4);
    \draw [line width=2pt,color=black] (3,2)-- (3,1);
    \draw [line width=2pt,color=black] (3,1)-- (4,1);
    \draw [line width=2pt,color=black] (4,1)-- (4,2);
    \draw [line width=2pt,color=black] (4,2)-- (3,2);
    \draw [line width=2pt,color=black] (5,4)-- (5,3);
    \draw [line width=2pt,color=black] (5,3)-- (6,3);
    \draw [line width=2pt,color=black] (6,3)-- (6,4);
    \draw [line width=2pt,color=black] (6,4)-- (5,4);
    \draw [line width=2pt,color=black] (5,2)-- (5,1);
    \draw [line width=2pt,color=black] (5,1)-- (6,1);
    \draw [line width=2pt,color=black] (6,1)-- (6,2);
    \draw [line width=2pt,color=black] (6,2)-- (5,2);
    \draw [line width=1pt,color=black] (4,0) -- (5,0) ;
    \draw [line width=1pt,color=black]  (5,0) -- (5,5);
    \draw [line width=1pt,color=black] (5,5) -- (4,5) ;
    \draw [line width=1pt,color=black] (4,0) -- (4,5) ;
    (6.25,2) -- (6.25,2.5) -- (1,2.5) -- (1,2)
    \draw [line width=1pt,color=black] (6.25,2) -- (6.25,2.5) ;
    \draw [line width=1pt,color=black]  (6.25,2.5) -- (1,2.5);
    \draw [line width=1pt,color=black]  (1,2.5) -- (1,2) ;
    \draw [line width=1pt,color=black]  (6.25,2) -- (1,2) ;
    
    \begin{scriptsize}
      \draw[color=black] (4.5,2.75) node {$I_3$};
      \draw[color=black] (3.5,2.25) node {$I_4$};
      \draw[color=black] (0.35,0.25) node {$(0,0)$};
      \draw[color=black] (6.65,4.75) node {$(7,5)$};
    \end{scriptsize}
  \end{tikzpicture}
  \caption{Example 3: The domain $\Omega$ with initial mesh and the domains of Integration for $I_3$ and $I_4$ \label{pic_ex_3_mesh}.}
\end{figure}

In this third example, we proceed to multiple goal functionals. The setup is
the same as in Example 2, but with a single $\alpha=0.01$ and 
multiple goal functionals:
\begin{itemize}
\item $I_1(u,q)=\frac{1}{2}\int_{\Omega}(u-\overline{u}^d)^2 dx \approx 1.15760  $,
\item $I_2(u,q)=\frac{1}{2}\int_{\Omega}(q-\overline{q}^d)^2 dx \approx 21.3305 $ ,
\item $I_3(u,q)=\int_{([4,5]\times \mathbb{R})\cap \Omega}u dx \approx  -0.236288$ ,
\item $I_4(u,q)=\int_{[1,\frac{25}{4}]\times [2,\frac{5}{2}]}q dx \approx 0.328042$ ,
\item $I_5(u,q)=\frac{1}{2}\int_{\Omega}u^2q^2 dx \approx 0.231615$.
\end{itemize}
The geometry alongside with the goal functionals $I_3$ and $I_4$ is
illustrated in Figure~\ref{pic_ex_3_mesh}.

\begin{minipage}[t]{0.45 \textwidth}
  \ifMAKEPICS
  \begin{gnuplot}[terminal=epslatex]
    set output "Figures/Example3btex.tex"
    set key left
    set key opaque
    set datafile separator "|"
    set logscale x
    # set yrange [-0.5:3]
    set xrange [50:100000]
    set grid ytics lc rgb "#bbbbbb" lw 1 lt 0
    set grid xtics lc rgb "#bbbbbb" lw 1 lt 0
    set xlabel '\text{DOFs}'
    set format '%g'
    plot  '< sqlite3 Data/Multigoalp4/Higher_Order/dataHigherOrderJE.db "SELECT DISTINCT DOFs, Ieff from data "' u 1:2 w  lp lw 3 title '  \footnotesize $I_{\text{eff}}$', \
    '< sqlite3 Data/Multigoalp4/Higher_Order/dataHigherOrderJE.db "SELECT DISTINCT DOFs, Ieffprimal from data "' u 1:2 w  lp lw 3 title '\footnotesize $I_{\text{effp}}$', \
    '< sqlite3 Data/Multigoalp4/Higher_Order/dataHigherOrderJE.db "SELECT DISTINCT DOFs, Ieffadjoint from data "' u 1:2 w  lp lw 3 title ' \footnotesize $I_{\text{effa}}$',\
    1 lw 3																	
  \end{gnuplot}
  \fi
  \scalebox{0.65}{% GNUPLOT: LaTeX picture with Postscript
\begingroup
  \makeatletter
  \providecommand\color[2][]{%
    \GenericError{(gnuplot) \space\space\space\@spaces}{%
      Package color not loaded in conjunction with
      terminal option `colourtext'%
    }{See the gnuplot documentation for explanation.%
    }{Either use 'blacktext' in gnuplot or load the package
      color.sty in LaTeX.}%
    \renewcommand\color[2][]{}%
  }%
  \providecommand\includegraphics[2][]{%
    \GenericError{(gnuplot) \space\space\space\@spaces}{%
      Package graphicx or graphics not loaded%
    }{See the gnuplot documentation for explanation.%
    }{The gnuplot epslatex terminal needs graphicx.sty or graphics.sty.}%
    \renewcommand\includegraphics[2][]{}%
  }%
  \providecommand\rotatebox[2]{#2}%
  \@ifundefined{ifGPcolor}{%
    \newif\ifGPcolor
    \GPcolorfalse
  }{}%
  \@ifundefined{ifGPblacktext}{%
    \newif\ifGPblacktext
    \GPblacktexttrue
  }{}%
  % define a \g@addto@macro without @ in the name:
  \let\gplgaddtomacro\g@addto@macro
  % define empty templates for all commands taking text:
  \gdef\gplbacktext{}%
  \gdef\gplfronttext{}%
  \makeatother
  \ifGPblacktext
    % no textcolor at all
    \def\colorrgb#1{}%
    \def\colorgray#1{}%
  \else
    % gray or color?
    \ifGPcolor
      \def\colorrgb#1{\color[rgb]{#1}}%
      \def\colorgray#1{\color[gray]{#1}}%
      \expandafter\def\csname LTw\endcsname{\color{white}}%
      \expandafter\def\csname LTb\endcsname{\color{black}}%
      \expandafter\def\csname LTa\endcsname{\color{black}}%
      \expandafter\def\csname LT0\endcsname{\color[rgb]{1,0,0}}%
      \expandafter\def\csname LT1\endcsname{\color[rgb]{0,1,0}}%
      \expandafter\def\csname LT2\endcsname{\color[rgb]{0,0,1}}%
      \expandafter\def\csname LT3\endcsname{\color[rgb]{1,0,1}}%
      \expandafter\def\csname LT4\endcsname{\color[rgb]{0,1,1}}%
      \expandafter\def\csname LT5\endcsname{\color[rgb]{1,1,0}}%
      \expandafter\def\csname LT6\endcsname{\color[rgb]{0,0,0}}%
      \expandafter\def\csname LT7\endcsname{\color[rgb]{1,0.3,0}}%
      \expandafter\def\csname LT8\endcsname{\color[rgb]{0.5,0.5,0.5}}%
    \else
      % gray
      \def\colorrgb#1{\color{black}}%
      \def\colorgray#1{\color[gray]{#1}}%
      \expandafter\def\csname LTw\endcsname{\color{white}}%
      \expandafter\def\csname LTb\endcsname{\color{black}}%
      \expandafter\def\csname LTa\endcsname{\color{black}}%
      \expandafter\def\csname LT0\endcsname{\color{black}}%
      \expandafter\def\csname LT1\endcsname{\color{black}}%
      \expandafter\def\csname LT2\endcsname{\color{black}}%
      \expandafter\def\csname LT3\endcsname{\color{black}}%
      \expandafter\def\csname LT4\endcsname{\color{black}}%
      \expandafter\def\csname LT5\endcsname{\color{black}}%
      \expandafter\def\csname LT6\endcsname{\color{black}}%
      \expandafter\def\csname LT7\endcsname{\color{black}}%
      \expandafter\def\csname LT8\endcsname{\color{black}}%
    \fi
  \fi
    \setlength{\unitlength}{0.0500bp}%
    \ifx\gptboxheight\undefined%
      \newlength{\gptboxheight}%
      \newlength{\gptboxwidth}%
      \newsavebox{\gptboxtext}%
    \fi%
    \setlength{\fboxrule}{0.5pt}%
    \setlength{\fboxsep}{1pt}%
\begin{picture}(7200.00,5040.00)%
    \gplgaddtomacro\gplbacktext{%
      \csname LTb\endcsname%
      \put(594,704){\makebox(0,0)[r]{\strut{}0.6}}%
      \csname LTb\endcsname%
      \put(594,1213){\makebox(0,0)[r]{\strut{}0.7}}%
      \csname LTb\endcsname%
      \put(594,1722){\makebox(0,0)[r]{\strut{}0.8}}%
      \csname LTb\endcsname%
      \put(594,2231){\makebox(0,0)[r]{\strut{}0.9}}%
      \csname LTb\endcsname%
      \put(594,2739){\makebox(0,0)[r]{\strut{}1}}%
      \csname LTb\endcsname%
      \put(594,3248){\makebox(0,0)[r]{\strut{}1.1}}%
      \csname LTb\endcsname%
      \put(594,3757){\makebox(0,0)[r]{\strut{}1.2}}%
      \csname LTb\endcsname%
      \put(594,4266){\makebox(0,0)[r]{\strut{}1.3}}%
      \csname LTb\endcsname%
      \put(594,4775){\makebox(0,0)[r]{\strut{}1.4}}%
      \csname LTb\endcsname%
      \put(1280,484){\makebox(0,0){\strut{}100}}%
      \csname LTb\endcsname%
      \put(3121,484){\makebox(0,0){\strut{}1000}}%
      \csname LTb\endcsname%
      \put(4962,484){\makebox(0,0){\strut{}10000}}%
      \csname LTb\endcsname%
      \put(6803,484){\makebox(0,0){\strut{}100000}}%
    }%
    \gplgaddtomacro\gplfronttext{%
      \csname LTb\endcsname%
      \put(3764,154){\makebox(0,0){\strut{}\text{DOFs}}}%
      \csname LTb\endcsname%
      \put(1914,4602){\makebox(0,0)[r]{\strut{}  \footnotesize $I_{eff}$}}%
      \csname LTb\endcsname%
      \put(1914,4382){\makebox(0,0)[r]{\strut{}\footnotesize $I_{effp}$}}%
      \csname LTb\endcsname%
      \put(1914,4162){\makebox(0,0)[r]{\strut{} \footnotesize $I_{effa}$}}%
      \csname LTb\endcsname%
      \put(1914,3942){\makebox(0,0)[r]{\strut{}1}}%
    }%
    \gplbacktext
    \put(0,0){\includegraphics{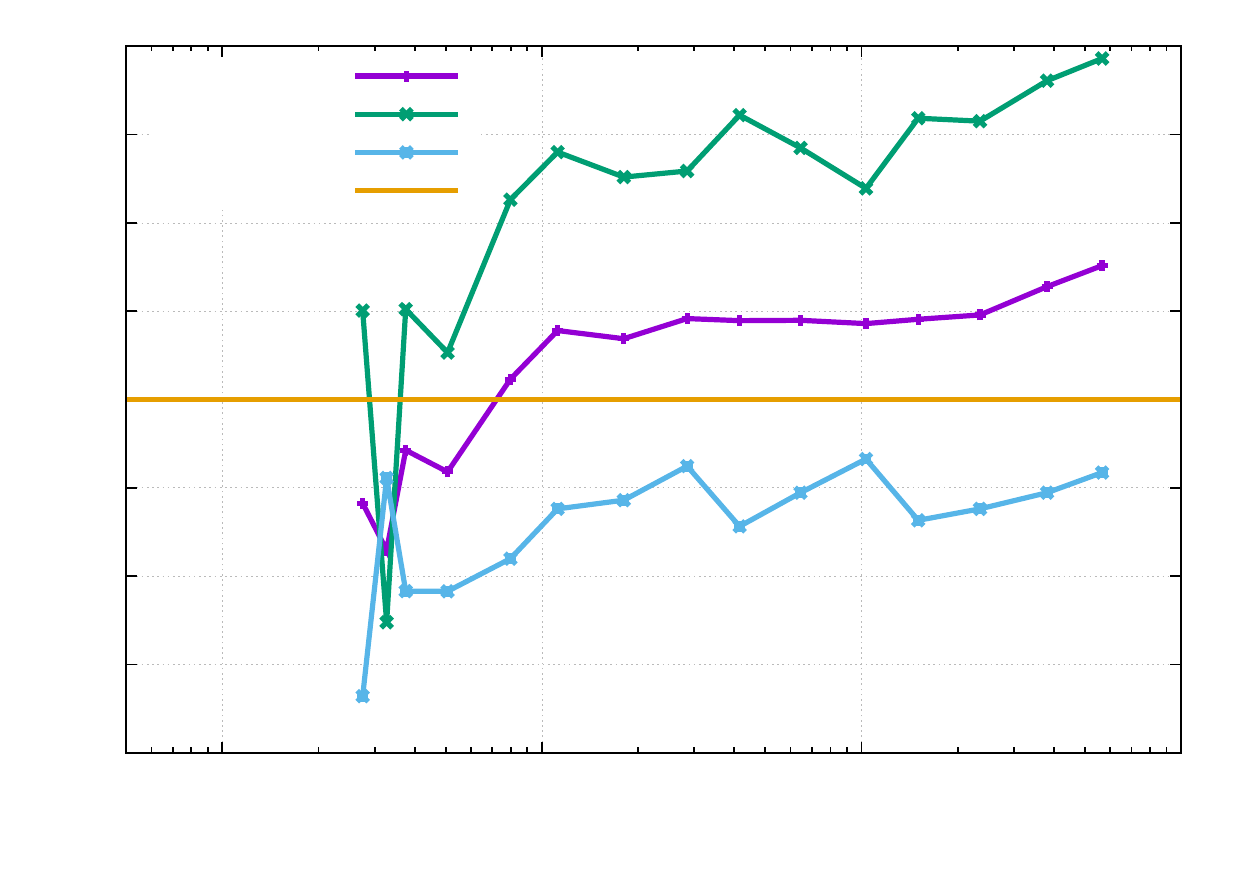}}%
    \gplfronttext
  \end{picture}%
\endgroup
}
  \captionof{figure}{Example 3. $I_{\text{eff}}$  vs DOFs for $p=4$, $\varepsilon=10^{-0}$.\label{pic_ex_3_a}}
\end{minipage} % 
\hfill
\begin{minipage}[t]{0.45 \textwidth}
  \ifMAKEPICS
  \begin{gnuplot}[terminal=epslatex]
    set output "Figures/Example3atex.tex"
    set key opaque
    set datafile separator "|"
    set logscale x
    set logscale y
    set grid ytics lc rgb "#bbbbbb" lw 1 lt 0
    set grid xtics lc rgb "#bbbbbb" lw 1 lt 0
    set xlabel '\text{DOFs}'
    set format '%g'
    plot  '< sqlite3 Data/Multigoalp4/Higher_Order/dataHigherOrderJE.db "SELECT DISTINCT DOFs, abs(exacterror) from data "' u 1:2 w  lp lw 3 title ' \footnotesize Error in $I_\mathfrak{E}$ (adp.)', \
    '< sqlite3 Data/Multigoalp4/Higher_Order/dataHigherOrderJE.db "SELECT DISTINCT DOFs, abs(estimatederror) from data "' u 1:2 w  lp lw 3 title ' \footnotesize Estimated Error', \
    '< sqlite3 Data/Multigoalp4/Higher_Order/dataHigherOrderJEuniform.db "SELECT DISTINCT DOFs, abs(exacterror) from data "' u 1:2 w  lp lw 3 title '  \footnotesize Error in $I_\mathfrak{E}$ (uni.)', \
    80/x   lw  1	title '$\mathcal{O}(\text{DOFs}^{-1})$'
    #	'< sqlite3 Data/Multigoalp4/Higher_Order/dataHigherOrderJE.db "SELECT DISTINCT DOFs, abs(estimatederror0) from data "' u 1:2 w  lp lw 3 title ' \footnotesize Error (a)', \
    '< sqlite3 Data/Multigoalp4/Higher_Order/dataHigherOrderJE.db "SELECT DISTINCT DOFs, abs(estimatederror1) from data "' u 1:2 w  lp lw 3 title ' \footnotesize Error (a)', \	
    '< sqlite3 Data/Multigoalp4/Higher_Order/dataHigherOrderJE.db "SELECT DISTINCT DOFs, abs(estimatederror2) from data "' u 1:2 w  lp lw 3 title ' \footnotesize Error (a)', \															
    '< sqlite3 Data/Multigoalp4/Higher_Order/dataHigherOrderJE.db "SELECT DISTINCT DOFs, correctederror from data "' u 1:2 w  lp lw 3 title ' \footnotesize Error (a)',\	
  \end{gnuplot}
  \fi
  {		\scalebox{0.65}{% GNUPLOT: LaTeX picture with Postscript
\begingroup
  \makeatletter
  \providecommand\color[2][]{%
    \GenericError{(gnuplot) \space\space\space\@spaces}{%
      Package color not loaded in conjunction with
      terminal option `colourtext'%
    }{See the gnuplot documentation for explanation.%
    }{Either use 'blacktext' in gnuplot or load the package
      color.sty in LaTeX.}%
    \renewcommand\color[2][]{}%
  }%
  \providecommand\includegraphics[2][]{%
    \GenericError{(gnuplot) \space\space\space\@spaces}{%
      Package graphicx or graphics not loaded%
    }{See the gnuplot documentation for explanation.%
    }{The gnuplot epslatex terminal needs graphicx.sty or graphics.sty.}%
    \renewcommand\includegraphics[2][]{}%
  }%
  \providecommand\rotatebox[2]{#2}%
  \@ifundefined{ifGPcolor}{%
    \newif\ifGPcolor
    \GPcolorfalse
  }{}%
  \@ifundefined{ifGPblacktext}{%
    \newif\ifGPblacktext
    \GPblacktexttrue
  }{}%
  % define a \g@addto@macro without @ in the name:
  \let\gplgaddtomacro\g@addto@macro
  % define empty templates for all commands taking text:
  \gdef\gplbacktext{}%
  \gdef\gplfronttext{}%
  \makeatother
  \ifGPblacktext
    % no textcolor at all
    \def\colorrgb#1{}%
    \def\colorgray#1{}%
  \else
    % gray or color?
    \ifGPcolor
      \def\colorrgb#1{\color[rgb]{#1}}%
      \def\colorgray#1{\color[gray]{#1}}%
      \expandafter\def\csname LTw\endcsname{\color{white}}%
      \expandafter\def\csname LTb\endcsname{\color{black}}%
      \expandafter\def\csname LTa\endcsname{\color{black}}%
      \expandafter\def\csname LT0\endcsname{\color[rgb]{1,0,0}}%
      \expandafter\def\csname LT1\endcsname{\color[rgb]{0,1,0}}%
      \expandafter\def\csname LT2\endcsname{\color[rgb]{0,0,1}}%
      \expandafter\def\csname LT3\endcsname{\color[rgb]{1,0,1}}%
      \expandafter\def\csname LT4\endcsname{\color[rgb]{0,1,1}}%
      \expandafter\def\csname LT5\endcsname{\color[rgb]{1,1,0}}%
      \expandafter\def\csname LT6\endcsname{\color[rgb]{0,0,0}}%
      \expandafter\def\csname LT7\endcsname{\color[rgb]{1,0.3,0}}%
      \expandafter\def\csname LT8\endcsname{\color[rgb]{0.5,0.5,0.5}}%
    \else
      % gray
      \def\colorrgb#1{\color{black}}%
      \def\colorgray#1{\color[gray]{#1}}%
      \expandafter\def\csname LTw\endcsname{\color{white}}%
      \expandafter\def\csname LTb\endcsname{\color{black}}%
      \expandafter\def\csname LTa\endcsname{\color{black}}%
      \expandafter\def\csname LT0\endcsname{\color{black}}%
      \expandafter\def\csname LT1\endcsname{\color{black}}%
      \expandafter\def\csname LT2\endcsname{\color{black}}%
      \expandafter\def\csname LT3\endcsname{\color{black}}%
      \expandafter\def\csname LT4\endcsname{\color{black}}%
      \expandafter\def\csname LT5\endcsname{\color{black}}%
      \expandafter\def\csname LT6\endcsname{\color{black}}%
      \expandafter\def\csname LT7\endcsname{\color{black}}%
      \expandafter\def\csname LT8\endcsname{\color{black}}%
    \fi
  \fi
    \setlength{\unitlength}{0.0500bp}%
    \ifx\gptboxheight\undefined%
      \newlength{\gptboxheight}%
      \newlength{\gptboxwidth}%
      \newsavebox{\gptboxtext}%
    \fi%
    \setlength{\fboxrule}{0.5pt}%
    \setlength{\fboxsep}{1pt}%
\begin{picture}(7200.00,5040.00)%
    \gplgaddtomacro\gplbacktext{%
      \csname LTb\endcsname%
      \put(990,704){\makebox(0,0)[r]{\strut{}0.0001}}%
      \csname LTb\endcsname%
      \put(990,1286){\makebox(0,0)[r]{\strut{}0.001}}%
      \csname LTb\endcsname%
      \put(990,1867){\makebox(0,0)[r]{\strut{}0.01}}%
      \csname LTb\endcsname%
      \put(990,2449){\makebox(0,0)[r]{\strut{}0.1}}%
      \csname LTb\endcsname%
      \put(990,3030){\makebox(0,0)[r]{\strut{}1}}%
      \csname LTb\endcsname%
      \put(990,3612){\makebox(0,0)[r]{\strut{}10}}%
      \csname LTb\endcsname%
      \put(990,4193){\makebox(0,0)[r]{\strut{}100}}%
      \csname LTb\endcsname%
      \put(990,4775){\makebox(0,0)[r]{\strut{}1000}}%
      \csname LTb\endcsname%
      \put(1122,484){\makebox(0,0){\strut{}100}}%
      \csname LTb\endcsname%
      \put(2542,484){\makebox(0,0){\strut{}1000}}%
      \csname LTb\endcsname%
      \put(3963,484){\makebox(0,0){\strut{}10000}}%
      \csname LTb\endcsname%
      \put(5383,484){\makebox(0,0){\strut{}100000}}%
      \csname LTb\endcsname%
      \put(6803,484){\makebox(0,0){\strut{}1e+06}}%
    }%
    \gplgaddtomacro\gplfronttext{%
      \csname LTb\endcsname%
      \put(3962,154){\makebox(0,0){\strut{}\text{DOFs}}}%
      \csname LTb\endcsname%
      \put(5816,4602){\makebox(0,0)[r]{\strut{} \footnotesize Error in $I_\mathfrak{E}$ (adp.)}}%
      \csname LTb\endcsname%
      \put(5816,4382){\makebox(0,0)[r]{\strut{} \footnotesize Estimated Error}}%
      \csname LTb\endcsname%
      \put(5816,4162){\makebox(0,0)[r]{\strut{}  \footnotesize Error in $I_\mathfrak{E}$ (uni.)}}%
      \csname LTb\endcsname%
      \put(5816,3942){\makebox(0,0)[r]{\strut{}$\mathcal{O}(\text{DOFs}^{-1})$}}%
    }%
    \gplbacktext
    \put(0,0){\includegraphics{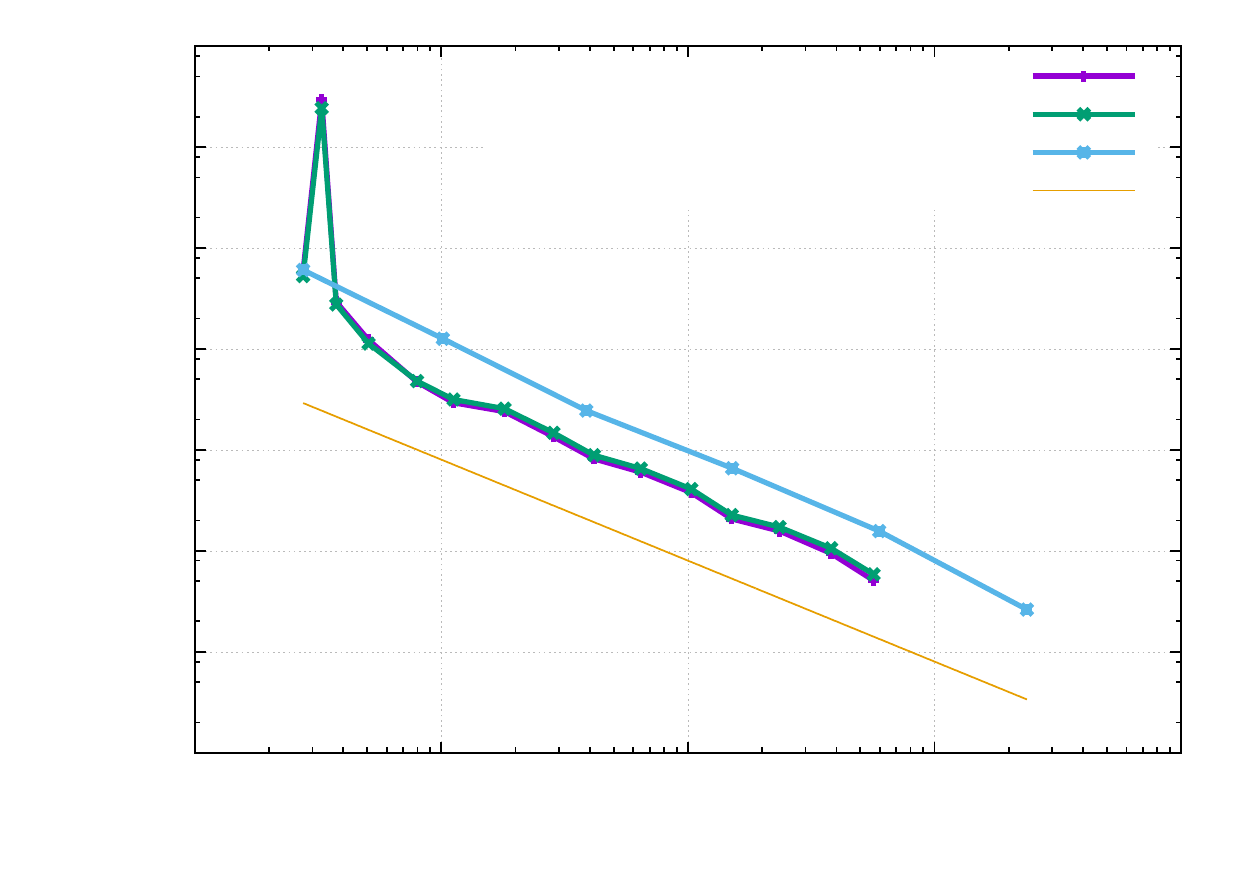}}%
    \gplfronttext
  \end{picture}%
\endgroup
}
    \captionof{figure}{Example 3. Error vs DOFs for $p=4$, $\varepsilon=10^{-0}$.\label{pic_ex_3_b}}}

\end{minipage}

\vspace*{1cm}

\begin{minipage}[t]{0.45 \textwidth}
  \ifMAKEPICS
  \begin{gnuplot}[terminal=epslatex]
    set output "Figures/Example3ctex.tex"
    set key right
    set key opaque
    set datafile separator "|"
    set yrange [0.00001:1000]
    set logscale 
    set grid ytics lc rgb "#bbbbbb" lw 1 lt 0
    set grid xtics lc rgb "#bbbbbb" lw 1 lt 0
    set xlabel '\text{DOFs}'
    set format '%g'
    plot  '< sqlite3 Data/Multigoalp4/Higher_Order/dataHigherOrderJE.db "SELECT DISTINCT DOFs, \"RelativeError(u)0\" from data "' u 1:2 w  lp lw 3 title '  \footnotesize $I_1$', \
    '< sqlite3 Data/Multigoalp4/Higher_Order/dataHigherOrderJE.db "SELECT DISTINCT DOFs, \"RelativeError(u)1\" from data "' u 1:2 w  lp lw 3 title '\footnotesize $I_2$', \
    '< sqlite3 Data/Multigoalp4/Higher_Order/dataHigherOrderJE.db "SELECT DISTINCT DOFs, \"RelativeError(u)2\" from data "' u 1:2 w  lp lw 3 title ' \footnotesize $I_3$',\
    '< sqlite3 Data/Multigoalp4/Higher_Order/dataHigherOrderJE.db "SELECT DISTINCT DOFs, \"RelativeError(u)3\" from data "' u 1:2 w  lp lw 3 title '  \footnotesize $I_4$', \
    '< sqlite3 Data/Multigoalp4/Higher_Order/dataHigherOrderJE.db "SELECT DISTINCT DOFs, \"RelativeError(u)4\" from data "' u 1:2 w  lp lw 3 title '\footnotesize $I_5$', \
    '< sqlite3 Data/Multigoalp4/Higher_Order/dataHigherOrderJE.db "SELECT DISTINCT DOFs, exacterror from data "' u 1:2 w  lp lw 3 title ' \footnotesize $I_{\mathfrak{E}}$',\													
  \end{gnuplot}
  \fi
  \scalebox{0.65}{% GNUPLOT: LaTeX picture with Postscript
\begingroup
  \makeatletter
  \providecommand\color[2][]{%
    \GenericError{(gnuplot) \space\space\space\@spaces}{%
      Package color not loaded in conjunction with
      terminal option `colourtext'%
    }{See the gnuplot documentation for explanation.%
    }{Either use 'blacktext' in gnuplot or load the package
      color.sty in LaTeX.}%
    \renewcommand\color[2][]{}%
  }%
  \providecommand\includegraphics[2][]{%
    \GenericError{(gnuplot) \space\space\space\@spaces}{%
      Package graphicx or graphics not loaded%
    }{See the gnuplot documentation for explanation.%
    }{The gnuplot epslatex terminal needs graphicx.sty or graphics.sty.}%
    \renewcommand\includegraphics[2][]{}%
  }%
  \providecommand\rotatebox[2]{#2}%
  \@ifundefined{ifGPcolor}{%
    \newif\ifGPcolor
    \GPcolorfalse
  }{}%
  \@ifundefined{ifGPblacktext}{%
    \newif\ifGPblacktext
    \GPblacktexttrue
  }{}%
  % define a \g@addto@macro without @ in the name:
  \let\gplgaddtomacro\g@addto@macro
  % define empty templates for all commands taking text:
  \gdef\gplbacktext{}%
  \gdef\gplfronttext{}%
  \makeatother
  \ifGPblacktext
    % no textcolor at all
    \def\colorrgb#1{}%
    \def\colorgray#1{}%
  \else
    % gray or color?
    \ifGPcolor
      \def\colorrgb#1{\color[rgb]{#1}}%
      \def\colorgray#1{\color[gray]{#1}}%
      \expandafter\def\csname LTw\endcsname{\color{white}}%
      \expandafter\def\csname LTb\endcsname{\color{black}}%
      \expandafter\def\csname LTa\endcsname{\color{black}}%
      \expandafter\def\csname LT0\endcsname{\color[rgb]{1,0,0}}%
      \expandafter\def\csname LT1\endcsname{\color[rgb]{0,1,0}}%
      \expandafter\def\csname LT2\endcsname{\color[rgb]{0,0,1}}%
      \expandafter\def\csname LT3\endcsname{\color[rgb]{1,0,1}}%
      \expandafter\def\csname LT4\endcsname{\color[rgb]{0,1,1}}%
      \expandafter\def\csname LT5\endcsname{\color[rgb]{1,1,0}}%
      \expandafter\def\csname LT6\endcsname{\color[rgb]{0,0,0}}%
      \expandafter\def\csname LT7\endcsname{\color[rgb]{1,0.3,0}}%
      \expandafter\def\csname LT8\endcsname{\color[rgb]{0.5,0.5,0.5}}%
    \else
      % gray
      \def\colorrgb#1{\color{black}}%
      \def\colorgray#1{\color[gray]{#1}}%
      \expandafter\def\csname LTw\endcsname{\color{white}}%
      \expandafter\def\csname LTb\endcsname{\color{black}}%
      \expandafter\def\csname LTa\endcsname{\color{black}}%
      \expandafter\def\csname LT0\endcsname{\color{black}}%
      \expandafter\def\csname LT1\endcsname{\color{black}}%
      \expandafter\def\csname LT2\endcsname{\color{black}}%
      \expandafter\def\csname LT3\endcsname{\color{black}}%
      \expandafter\def\csname LT4\endcsname{\color{black}}%
      \expandafter\def\csname LT5\endcsname{\color{black}}%
      \expandafter\def\csname LT6\endcsname{\color{black}}%
      \expandafter\def\csname LT7\endcsname{\color{black}}%
      \expandafter\def\csname LT8\endcsname{\color{black}}%
    \fi
  \fi
    \setlength{\unitlength}{0.0500bp}%
    \ifx\gptboxheight\undefined%
      \newlength{\gptboxheight}%
      \newlength{\gptboxwidth}%
      \newsavebox{\gptboxtext}%
    \fi%
    \setlength{\fboxrule}{0.5pt}%
    \setlength{\fboxsep}{1pt}%
\begin{picture}(7200.00,5040.00)%
    \gplgaddtomacro\gplbacktext{%
      \csname LTb\endcsname%
      \put(990,704){\makebox(0,0)[r]{\strut{}1e-05}}%
      \csname LTb\endcsname%
      \put(990,1213){\makebox(0,0)[r]{\strut{}0.0001}}%
      \csname LTb\endcsname%
      \put(990,1722){\makebox(0,0)[r]{\strut{}0.001}}%
      \csname LTb\endcsname%
      \put(990,2231){\makebox(0,0)[r]{\strut{}0.01}}%
      \csname LTb\endcsname%
      \put(990,2740){\makebox(0,0)[r]{\strut{}0.1}}%
      \csname LTb\endcsname%
      \put(990,3248){\makebox(0,0)[r]{\strut{}1}}%
      \csname LTb\endcsname%
      \put(990,3757){\makebox(0,0)[r]{\strut{}10}}%
      \csname LTb\endcsname%
      \put(990,4266){\makebox(0,0)[r]{\strut{}100}}%
      \csname LTb\endcsname%
      \put(990,4775){\makebox(0,0)[r]{\strut{}1000}}%
      \csname LTb\endcsname%
      \put(1122,484){\makebox(0,0){\strut{}100}}%
      \csname LTb\endcsname%
      \put(3016,484){\makebox(0,0){\strut{}1000}}%
      \csname LTb\endcsname%
      \put(4909,484){\makebox(0,0){\strut{}10000}}%
      \csname LTb\endcsname%
      \put(6803,484){\makebox(0,0){\strut{}100000}}%
    }%
    \gplgaddtomacro\gplfronttext{%
      \csname LTb\endcsname%
      \put(3962,154){\makebox(0,0){\strut{}\text{DOFs}}}%
      \csname LTb\endcsname%
      \put(5816,4602){\makebox(0,0)[r]{\strut{}  \footnotesize $I_1$}}%
      \csname LTb\endcsname%
      \put(5816,4382){\makebox(0,0)[r]{\strut{}\footnotesize $I_2$}}%
      \csname LTb\endcsname%
      \put(5816,4162){\makebox(0,0)[r]{\strut{} \footnotesize $I_3$}}%
      \csname LTb\endcsname%
      \put(5816,3942){\makebox(0,0)[r]{\strut{}  \footnotesize $I_4$}}%
      \csname LTb\endcsname%
      \put(5816,3722){\makebox(0,0)[r]{\strut{}\footnotesize $I_5$}}%
      \csname LTb\endcsname%
      \put(5816,3502){\makebox(0,0)[r]{\strut{} \footnotesize $I_{\mathfrak{E}}$}}%
    }%
    \gplbacktext
    \put(0,0){\includegraphics{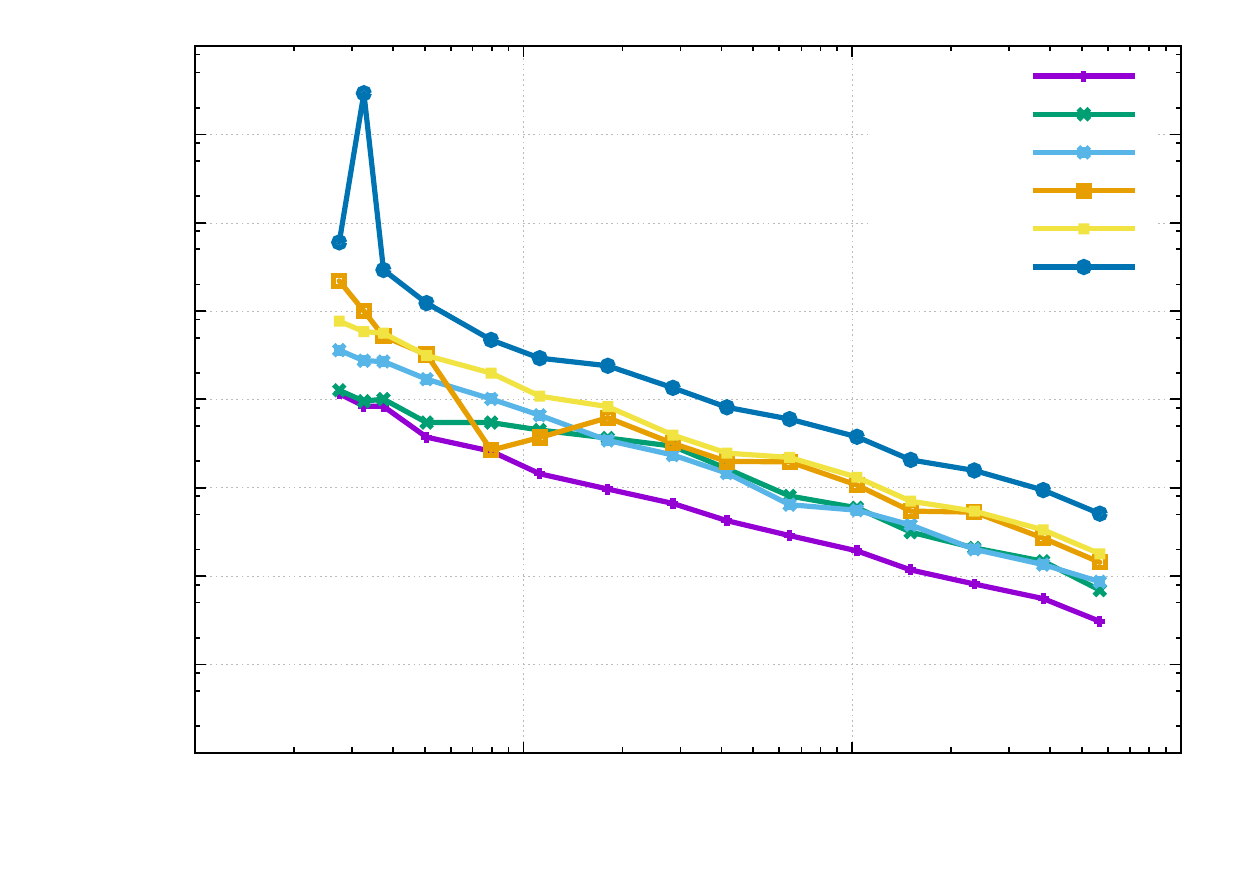}}%
    \gplfronttext
  \end{picture}%
\endgroup
}
  \captionof{figure}{Example 3. Relative error vs DOFs for $p=4$, $\varepsilon=10^{-0}$.\label{pic_ex_3_c}}
\end{minipage} % 
\hfill
\begin{minipage}[t]{0.45 \textwidth}
  \ifMAKEPICS
  \begin{gnuplot}[terminal=epslatex]
    set output "Figures/Example3dtex.tex"
    set key bottom left
    set key opaque
    set datafile separator "|"
    set logscale x
    set logscale y
    set yrange [0.00001:1000]
    set grid ytics lc rgb "#bbbbbb" lw 1 lt 0
    set grid xtics lc rgb "#bbbbbb" lw 1 lt 0
    set xlabel '\text{DOFs}'
    set format '%g'
    plot \
    '< sqlite3 Data/Multigoalp4/Higher_Order/dataHigherOrderJE.db "SELECT DISTINCT DOFs, abs(estimatederror) from data "' u 1:2 w  lp lw 3 title ' \footnotesize $\eta_h^{(2)}$', \
    '< sqlite3 Data/Multigoalp4/Higher_Order/dataHigherOrderJE.db "SELECT DISTINCT DOFs, abs(estimatederror1) from data "' u 1:2 w  lp lw 3 title ' \footnotesize $\rho_u$', \
    '< sqlite3 Data/Multigoalp4/Higher_Order/dataHigherOrderJE.db "SELECT DISTINCT DOFs, abs(estimatederror2) from data "' u 1:2 w  lp lw 3 title ' \footnotesize $\rho_p$', \	
    '< sqlite3 Data/Multigoalp4/Higher_Order/dataHigherOrderJE.db "SELECT DISTINCT DOFs, abs(estimatederror3) from data "' u 1:2 w  lp lw 3 title ' \footnotesize $\rho_z$', \			
    '< sqlite3 Data/Multigoalp4/Higher_Order/dataHigherOrderJE.db "SELECT DISTINCT DOFs, abs(estimatederror4) from data "' u 1:2 w  lp lw 3 title ' \footnotesize $\rho_v$', \					
    '< sqlite3 Data/Multigoalp4/Higher_Order/dataHigherOrderJE.db "SELECT DISTINCT DOFs, abs(estimatederror5) from data "' u 1:2 w  lp lw 3 title ' \footnotesize $\rho_q$', \					
    '< sqlite3 Data/Multigoalp4/Higher_Order/dataHigherOrderJE.db "SELECT DISTINCT DOFs, abs(estimatederror6) from data "' u 1:2 w  lp lw 3 title ' \footnotesize $\rho_y$',				
    #				10/x   lw  1											
    #					 '< sqlite3 Data/Multigoalp4/Higher_Order/dataHigherOrderJE.db "SELECT DISTINCT DOFs, abs(exacterror) from data "' u 1:2 w  lp lw 3 title ' \footnotesize Error in $I_\mathfrak{E}$', \
  \end{gnuplot}
  \fi
  {		\scalebox{0.65}{% GNUPLOT: LaTeX picture with Postscript
\begingroup
  \makeatletter
  \providecommand\color[2][]{%
    \GenericError{(gnuplot) \space\space\space\@spaces}{%
      Package color not loaded in conjunction with
      terminal option `colourtext'%
    }{See the gnuplot documentation for explanation.%
    }{Either use 'blacktext' in gnuplot or load the package
      color.sty in LaTeX.}%
    \renewcommand\color[2][]{}%
  }%
  \providecommand\includegraphics[2][]{%
    \GenericError{(gnuplot) \space\space\space\@spaces}{%
      Package graphicx or graphics not loaded%
    }{See the gnuplot documentation for explanation.%
    }{The gnuplot epslatex terminal needs graphicx.sty or graphics.sty.}%
    \renewcommand\includegraphics[2][]{}%
  }%
  \providecommand\rotatebox[2]{#2}%
  \@ifundefined{ifGPcolor}{%
    \newif\ifGPcolor
    \GPcolorfalse
  }{}%
  \@ifundefined{ifGPblacktext}{%
    \newif\ifGPblacktext
    \GPblacktexttrue
  }{}%
  % define a \g@addto@macro without @ in the name:
  \let\gplgaddtomacro\g@addto@macro
  % define empty templates for all commands taking text:
  \gdef\gplbacktext{}%
  \gdef\gplfronttext{}%
  \makeatother
  \ifGPblacktext
    % no textcolor at all
    \def\colorrgb#1{}%
    \def\colorgray#1{}%
  \else
    % gray or color?
    \ifGPcolor
      \def\colorrgb#1{\color[rgb]{#1}}%
      \def\colorgray#1{\color[gray]{#1}}%
      \expandafter\def\csname LTw\endcsname{\color{white}}%
      \expandafter\def\csname LTb\endcsname{\color{black}}%
      \expandafter\def\csname LTa\endcsname{\color{black}}%
      \expandafter\def\csname LT0\endcsname{\color[rgb]{1,0,0}}%
      \expandafter\def\csname LT1\endcsname{\color[rgb]{0,1,0}}%
      \expandafter\def\csname LT2\endcsname{\color[rgb]{0,0,1}}%
      \expandafter\def\csname LT3\endcsname{\color[rgb]{1,0,1}}%
      \expandafter\def\csname LT4\endcsname{\color[rgb]{0,1,1}}%
      \expandafter\def\csname LT5\endcsname{\color[rgb]{1,1,0}}%
      \expandafter\def\csname LT6\endcsname{\color[rgb]{0,0,0}}%
      \expandafter\def\csname LT7\endcsname{\color[rgb]{1,0.3,0}}%
      \expandafter\def\csname LT8\endcsname{\color[rgb]{0.5,0.5,0.5}}%
    \else
      % gray
      \def\colorrgb#1{\color{black}}%
      \def\colorgray#1{\color[gray]{#1}}%
      \expandafter\def\csname LTw\endcsname{\color{white}}%
      \expandafter\def\csname LTb\endcsname{\color{black}}%
      \expandafter\def\csname LTa\endcsname{\color{black}}%
      \expandafter\def\csname LT0\endcsname{\color{black}}%
      \expandafter\def\csname LT1\endcsname{\color{black}}%
      \expandafter\def\csname LT2\endcsname{\color{black}}%
      \expandafter\def\csname LT3\endcsname{\color{black}}%
      \expandafter\def\csname LT4\endcsname{\color{black}}%
      \expandafter\def\csname LT5\endcsname{\color{black}}%
      \expandafter\def\csname LT6\endcsname{\color{black}}%
      \expandafter\def\csname LT7\endcsname{\color{black}}%
      \expandafter\def\csname LT8\endcsname{\color{black}}%
    \fi
  \fi
    \setlength{\unitlength}{0.0500bp}%
    \ifx\gptboxheight\undefined%
      \newlength{\gptboxheight}%
      \newlength{\gptboxwidth}%
      \newsavebox{\gptboxtext}%
    \fi%
    \setlength{\fboxrule}{0.5pt}%
    \setlength{\fboxsep}{1pt}%
\begin{picture}(7200.00,5040.00)%
    \gplgaddtomacro\gplbacktext{%
      \csname LTb\endcsname%
      \put(990,704){\makebox(0,0)[r]{\strut{}1e-05}}%
      \csname LTb\endcsname%
      \put(990,1213){\makebox(0,0)[r]{\strut{}0.0001}}%
      \csname LTb\endcsname%
      \put(990,1722){\makebox(0,0)[r]{\strut{}0.001}}%
      \csname LTb\endcsname%
      \put(990,2231){\makebox(0,0)[r]{\strut{}0.01}}%
      \csname LTb\endcsname%
      \put(990,2740){\makebox(0,0)[r]{\strut{}0.1}}%
      \csname LTb\endcsname%
      \put(990,3248){\makebox(0,0)[r]{\strut{}1}}%
      \csname LTb\endcsname%
      \put(990,3757){\makebox(0,0)[r]{\strut{}10}}%
      \csname LTb\endcsname%
      \put(990,4266){\makebox(0,0)[r]{\strut{}100}}%
      \csname LTb\endcsname%
      \put(990,4775){\makebox(0,0)[r]{\strut{}1000}}%
      \csname LTb\endcsname%
      \put(1122,484){\makebox(0,0){\strut{}100}}%
      \csname LTb\endcsname%
      \put(3016,484){\makebox(0,0){\strut{}1000}}%
      \csname LTb\endcsname%
      \put(4909,484){\makebox(0,0){\strut{}10000}}%
      \csname LTb\endcsname%
      \put(6803,484){\makebox(0,0){\strut{}100000}}%
    }%
    \gplgaddtomacro\gplfronttext{%
      \csname LTb\endcsname%
      \put(3962,154){\makebox(0,0){\strut{}\text{DOFs}}}%
      \csname LTb\endcsname%
      \put(2310,2197){\makebox(0,0)[r]{\strut{} \footnotesize $\eta_h^{(2)}$}}%
      \csname LTb\endcsname%
      \put(2310,1977){\makebox(0,0)[r]{\strut{} \footnotesize $\rho_u$}}%
      \csname LTb\endcsname%
      \put(2310,1757){\makebox(0,0)[r]{\strut{} \footnotesize $\rho_p$}}%
      \csname LTb\endcsname%
      \put(2310,1537){\makebox(0,0)[r]{\strut{} \footnotesize $\rho_z$}}%
      \csname LTb\endcsname%
      \put(2310,1317){\makebox(0,0)[r]{\strut{} \footnotesize $\rho_v$}}%
      \csname LTb\endcsname%
      \put(2310,1097){\makebox(0,0)[r]{\strut{} \footnotesize $\rho_q$}}%
      \csname LTb\endcsname%
      \put(2310,877){\makebox(0,0)[r]{\strut{} \footnotesize $\rho_y$}}%
    }%
    \gplbacktext
    \put(0,0){\includegraphics{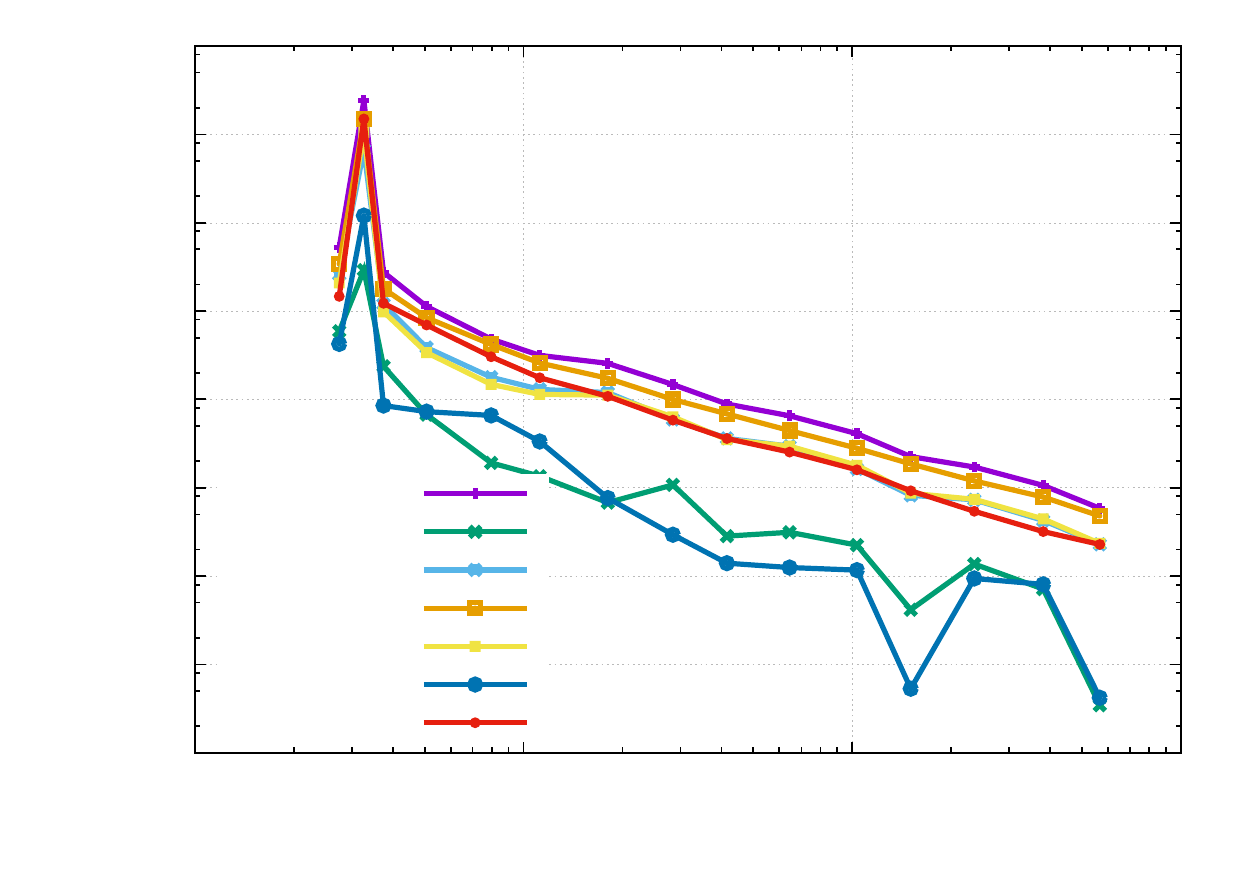}}%
    \gplfronttext
  \end{picture}%
\endgroup
}
    \captionof{figure}{Example 3. All error estimator parts for $p=4$, $\varepsilon=10^{-0}$.\label{pic_ex_3_d}	}}
  
\end{minipage}

%\vspace*{1cm}

\newpage
The reference values are computed on a fine grid ($716\, 792$ DOFs for $q$ +
$730\, 199$ DOFs for $u$) , which is obtained by 12 adaptive refinements for
$J_\mathfrak{E}$ followed by two uniform h-refinements and one uniform
p-refinement. Our findings are displayed in the 
Figures~\ref{pic_ex_3_a},~\ref{pic_ex_3_b},~\ref{pic_ex_3_c} and~\ref{pic_ex_3_d}.
In Figure~\ref{pic_ex_3_a} the calculated effectivity indices are excellent 
in view of the nonlinearities of the domain, state equation and multiple goal
functionals.
Curves of the errors and estimators are shown in the Figures~\ref{pic_ex_3_b},~\ref{pic_ex_3_c}
and~\ref{pic_ex_3_d}. Here, the combined functional (as expected) bounds all
single functionals. In Figure~\ref{pic_ex_3_b}, we observe that adaptive 
refinement pays off in delivering the same error as uniform mesh refinement,
but 
with a lower computational cost. The convergence rates are the same, which
lies in the fact that the control is chosen in such a way that a sufficiently smooth final 
solution is obtained.  
Finally, we compare the adaptive stopping rule used in Algorithm
\ref{inexat_newton_algorithm_for_multiple_goal_functionals} with the standard
stopping rule, which is used in the \verb|DOpElib|~\cite{dope,DOpElib}
algorithm  \verb|DOpE::ReducedNewtonAlgorithm::Solve| with the absolute
residual \verb|nonlinear_global_tol = 1.e-7| and relative residual
\verb|nonlinear_tol= 8.e-5|. Since the discretization error estimate is not
given for $l=0$ in
Algorithm~\ref{inexat_newton_algorithm_for_multiple_goal_functionals}, we use
$\eta_h^{l-1}=10^{-5}$.

We abbreviate the first algorithm with AN (Adaptive
Newton) and the second 
algorithm with FN (Full Newton). In Table~\ref{Table: Example3Ieff}, we
monitor that the $I_{\text{eff}}$ show a pretty similar behavior even for the
adaptive stopping rule. Even though we need $1-3$ iterations in case of the
adaptive stopping rule compared to $2-17$ iterations for the standard stopping
rule, which is illustrated in Table~\ref{Table: Example3Ieff} as
well. Furthermore, we want to notice that the refined meshes for both algorithms
coincide exactly up to $l=7$. For $l=8$, it is exactly 
one element, which is refined additionally in the case of FN. If we compare
the corrected effectivity indices 
\[
  I_{\text{eff,c}}:= \frac{\eta_h^{(2)}+
    \eta_k}{I(\overline{u},\overline{q})-I(\tilde{u}_h,\tilde{q_h})}
\]
for the
two stopping rules, we observe that they coincide even more 
after the correction.

In Table~\ref{Table: Example3Iteration},
the comparison between the estimated iteration error and the real error in the
combined functional is shown. The ratio between $\eta_k$ and the error mimics
the choice of $\gamma$ in Algorithm
\ref{inexat_newton_algorithm_for_multiple_goal_functionals} 
for our adaptive stopping rule, whereas there is almost no 
correlation for the standard stopping rule.

\begin{table}[H]
  \caption{Example 3: Comparison of Newton's Method with adaptive stopping rule
    (AN) and the classical, non-adaptive, Newton method (FN); $It_F$:
    number of iterations for FN, $It_A$: number of iterations for AN,
    $|\mathcal{T}_{h,F}|$: number of elements in the adaptive mesh
    resulting from using FN,  $|\mathcal{T}_{h,A}|$: number of elements
    in the adaptive mesh resulting from using AN, $I_{\text{eff,F}}$: $I_{\text{eff}}$
    for FN,  $I_{\text{eff,A}}$: $I_{\text{eff}}$ for AN, $I_{\text{eff,c,F}}$: $I_{\text{eff,c}}$
    for FN,  $I_{\text{eff,c,A}}$: $I_{\text{eff,c}}$ for AN \label{Table: Example3Ieff}}
  \renewcommand{\arraystretch}{1.3}
  \centering
  \begin{filecontents*}{Data/Multigoalp4/Higher_Order/IterationErrorEstimates/Differences.csv}
refinementcycle,NewtonIterationsf,NewtonIterationsa,DOFscontrol,DOFscontrola,Additionally Refined Cells,Ieff,Ieffa,Ieffcorrected,Ieffcorrecteda,exacterror,exacterrora,IterationError,IterationErrora
0,4,3,116,116,0,0.882,0.882,0.882,0.882,5.98E+00,5.98E+00,-1.99E-05,-3.20E-06
1,5,3,137,137,0,0.829,0.830,0.829,0.830,2.93E+02,2.75E+02,1.22E-02,-4.19E-02
2,5,1,158,158,0,0.942,0.933,0.942,0.941,2.92E+00,2.98E+00,-5.18E-05,2.46E-02
3,6,2,215,215,0,0.918,0.912,0.918,0.918,1.23E+00,1.24E+00,6.72E-05,7.71E-03
4,6,2,347,347,0,1.023,1.019,1.023,1.022,4.70E-01,4.73E-01,9.75E-06,1.67E-03
5,6,2,494,494,0,1.078,1.068,1.078,1.076,2.93E-01,2.96E-01,4.44E-06,2.59E-03
6,7,2,800,800,0,1.069,1.067,1.069,1.069,2.40E-01,2.40E-01,1.99E-06,3.48E-04
7,12,2,1283,1283,0,1.091,1.087,1.091,1.090,1.35E-01,1.36E-01,-2.41E-06,5.21E-04
8,17,2,1898,1895,1,1.089,1.093,1.089,1.090,8.13E-02,8.09E-02,-3.38E-06,-2.73E-04
9,6,2,2966,2957,3,1.089,1.090,1.090,1.090,6.00E-02,6.00E-02,1.94E-05,-2.64E-05
10,9,2,4802,4790,4,1.085,1.088,1.085,1.083,3.77E-02,3.77E-02,3.65E-06,-2.03E-04
11,5,2,7097,7091,2,1.089,1.077,1.089,1.093,2.07E-02,2.11E-02,-8.56E-07,3.45E-04
12,4,2,11153,11099,18,1.096,1.088,1.097,1.095,1.57E-02,1.58E-02,2.61E-05,1.17E-04
13,2,2,18206,18140,22,1.128,1.129,1.122,1.123,9.41E-03,9.45E-03,-5.67E-05,-5.96E-05
14,2,2,27341,27269,24,1.151,1.152,1.136,1.137,5.09E-03,5.11E-03,-7.89E-05,-7.85E-05
\end{filecontents*}
\pgfplotstabletypeset[
  columns={[index]0,[index]1,[index]2,[index]3,[index]4,[index]6,[index]7,[index]8,[index]9},
  col sep=comma,
  columns/refinementcycle/.style={column name=$l$},
  columns/NewtonIterationsf/.style={column name=$It_F$,column type={|r}},
  columns/NewtonIterationsa/.style={column name=$It_A$,column type={|r}},
  columns/DOFscontrol/.style={column name=$|\mathcal{T}_{h,F}|$,,column type={|r},1000 sep={\,}},
  columns/DOFscontrola/.style={column name=$|\mathcal{T}_{h,A}|$,column type={|r},1000 sep={\,}},
  columns/Ieff/.style={column name=$I_{\text{eff,F}}$,fixed zerofill,column type={|r},	precision=3},
  columns/Ieffa/.style={column name=$I_{\text{eff,A}}$,fixed zerofill,column type={|r},	precision=3},
  columns/Ieffcorrected/.style={column name=$I_{\text{eff,c,F}}$,column type={|r}, fixed zerofill,
    precision=3},
  columns/Ieffcorrecteda/.style={column name=$I_{\text{eff,c,A}}$,column type={|r|},fixed zerofill,
    precision=3},
  column type/.add={|}{},
  every head row/.style={},
  every last row/.style={after row= 	\hline},
  before row={
    \hline
  }
  ]{Data/Multigoalp4/Higher_Order/IterationErrorEstimates/Differences.csv}
\end{table}
\begin{table}[H]
  \caption{Example 3: Comparison of Newton's Method with adaptive stopping rule
    (AN) and the classical, non-adaptive, Newton method (FN);  $\eta_{k,F}$: iteration
    error estimate for FN, $\eta_{k,F}$: iteration error estimate for
    AN. 	\label{Table: Example3Iteration}}
  \renewcommand{\arraystretch}{1.3}
  \centering
  \pgfplotstabletypeset[
  columns={[index]0,[index]10,[index]11,[index]12,[index]13},
  col sep=comma,
  columns/refinementcycle/.style={column name=$l$},
  columns/NewtonIterationsf/.style={column name=$It_F$,column type={|r}},
  columns/NewtonIterationsa/.style={column name=$It_A$,column type={|r}},
  columns/DOFscontrol/.style={column name=$|\mathcal{T}_{h,F}|$,,column type={|r},1000 sep={\,}},
  columns/DOFscontrola/.style={column name=$|\mathcal{T}_{h,A}|$,column type={|r},1000 sep={\,}},
  columns/Ieff/.style={column name=$I_{\text{eff,F}}$,fixed zerofill,},
  columns/Ieffa/.style={column name=$I_{\text{eff,A}}$,fixed zerofill,},
  columns/Ieffcorrected/.style={column name=$I_{\text{eff,c,F}}$, fixed zerofill,
    precision=3},
  columns/Ieffcorrecteda/.style={column name=$I_{\text{eff,c,A}}$,fixed zerofill,
    precision=3},
  columns/exacterror/.style={column name=Error in $I_{\mathfrak{E},F}$, fixed zerofill,
    precision=2,sci,sci zerofill,column type={|c}},
  columns/exacterrora/.style={column name=Error in $I_{\mathfrak{E},A}$, fixed zerofill,
    precision=2,sci,,sci zerofill,column type={|c}},
  columns/IterationError/.style={column name=$\eta_{k,F}$, fixed zerofill,
    precision=3,sci,column type={|c}},
  columns/IterationErrora/.style={column name=$\eta_{k,A}$, fixed zerofill,
    precision=3,sci,column type={|c|}},
  column type/.add={|}{},
  every last row/.style={after row= 	\hline},
  before row={
    \hline
  }
  ]{Data/Multigoalp4/Higher_Order/IterationErrorEstimates/Differences.csv}
\end{table}

%%%%%%%%%%%%%%%%%%%%%%%%%%%%%%%%%%%%%%%%%%%%%%%%%%%%%%%%%%%%%%%%%%%%%%%%%%%%% 
\section{Conclusions}
\label{sec_conclusions}
In this work, we developed a novel a posteriori multiple goal-oriented error estimation 
for optimal control problems subject to a nonlinear state equation.
The error estimator also serves for balancing the discretization and nonlinear iteration error.
The overall optimization problem is solved via a reduced approach in which
the state equation is eliminated by a control-to-state solution operator. In
Section~\ref{sec_error_reduced_system}, the theoretical results yield 
an a posteriori estimate for a single goal functional.  The extension 
to multiple goal functionals was made in Section~\ref{sec_multigoal}. 
Based on these theoretical aspects, the algorithmic details were worked 
out in the following section. Three numerical examples were investigated.
In the first example, our approach  was tested against configurations 
known in the literature. The Examples 2 and 3 are more advanced by considering 
the regularized $p$-Laplacian as nonlinear state equation. The main criterion 
whether the  proposed error estimator works sufficiently well is given 
by the effectivity index. In the numerical examples,  values around one
were obtained. These are excellent findings in view of the 
challenging nature of the underlying problem configuration; namely 
domain (corner) singularities, quasi-linear state equations within an optimal control setting, 
and finally multiple nonlinear goal functionals. Ongoing work considers the
extension 
to elasticity and more practical applications.

%%%%%%%%%%%%%%%%%%%%%%%%%%%%%%%%%%%%%%%%%%%%%%%%%%%%%%%%%%%%%%%%%%%%%%%%%%%%% 
\section{Acknowledgments}
This work has been supported by the Austrian Science Fund (FWF) under the grant
P 29181
`Goal-Oriented Error Control for Phase-Field Fracture Coupled to
Multiphysics Problems'
and the DFG-SPP 1962 
`Non-smooth and Complementarity-based Distributed Parameter Systems: Simulation
and Hierarchical Optimization'
within the project `Optimizing Fracture Propagation Using a
Phase-Field Approach' under grant numbers NE1941/1-1 and WO1936/4-1.

\bibliography{./lit}
\bibliographystyle{abbrv}
% ##############################################################################%
\end{document}